\documentclass[a4paper,11pt]{article}
\usepackage[a4paper, total={6.45in, 9.7in}]{geometry}
\usepackage[utf8]{inputenc} % Required for inputting international characters
\usepackage[T1]{fontenc} % Output font encoding for international characters
\usepackage{amsmath}
\usepackage{amsthm}
\usepackage{verbatim}
\usepackage{mathtools} 
\usepackage{extarrows} 
\usepackage{amsthm}
\usepackage{xparse}
\usepackage{hyperref}
\usepackage{cleveref}
\usepackage{setspace}
\usepackage{mathpazo} % Use the Palatino font by default
\usepackage{amssymb}
\usepackage{tikz-cd}
\usepackage{float}
\usepackage{tabstackengine}
\setstackgap{L}{3.5mm}
\usepackage{graphicx}
\usepackage{tikz}
\usepackage[T1]{fontenc}
\usepackage{imakeidx}

\usepackage{xcolor}
\usepackage{amssymb,latexsym} 
 
\numberwithin{equation}{section} 
 
\newtheorem{theorem}{Theorem}[section] 
\newtheorem*{teo*}{Theorem}
\newtheorem{proposition}[theorem]{Proposition}

\newtheorem{lemma}[theorem]{Lemma} 
 
\theoremstyle{definition} 
 
\newtheorem{remark}[theorem]{Remark} 
\newtheorem{example}[theorem]{Example} 
\newtheorem{definition}[theorem]{Definition}

\def\PP{\mathbb{P}}

\def\Poly{\mathbf{P}} 

\usepackage[maxbibnames=99,style=alphabetic]{biblatex}
\usepackage{ytableau}
\usepackage{hyperref}
\usepackage[autostyle=true]{csquotes}

\makeindex

\begin{document}

\begin{center}
\textbf{\Large A categorification of cluster algebras of type B and C through symmetric quivers} 
\vspace{5mm}
 \\Azzurra Ciliberti\footnote{DIPARTIMENTO DI MATEMATICA “GUIDO CASTELNUOVO”, SAPIENZA UNIVERSIT\`{A} DI ROMA.\\
 \textit{Email address}: \textbf{azzurra.ciliberti@uniroma1.it}}

 \end{center}
 \begin{abstract}
  \noindent 
  \begin{comment}

 We associate cluster algebras of type $B_n$ and $C_n$ with principal coefficients to a regular polygon with $2n+2$ vertices. We set a correspondence between cluster variables and symmetric indecomposable representations of cluster-tilted bound symmetric quivers of type $A_{2n-1}$. We find a Caldero-Chapoton map in this setting. We also give a categorical interpretation of the cluster expansion formula in the case of acyclic quivers.
       
  \end{comment}
We express cluster variables of type $B_n$ and $C_n$ in terms of cluster variables of type $A_n$. Then we associate a cluster tilted bound symmetric quiver $Q$ of type $A_{2n-1}$ to any seed of a cluster algebra of type $B_n$ and $C_n$. Under this correspondence, cluster variables of type $B_n$ (resp. $C_n$) correspond to orthogonal (resp. symplectic) indecomposable representations of $Q$. We find a Caldero-Chapoton map in this setting. We also give a categorical interpretation of the cluster expansion formula in the case of acyclic quivers.
\end{abstract}
\tableofcontents
\section{Introduction}

Let $\Poly_{n+3}$ be the regular polygon with $n+3$ vertices. It is well known that clusters of cluster algebras of type $A_n$ correspond to triangulations of $\Poly_{n+3}$, while cluster variables correspond to diagonals. On the other hand, let $\Poly_{2n+2}$ be the regular polygon with $2n+2$ vertices, and let $\theta$ be the rotation of $180^\circ$. Fomin and Zelevinsky showed in \cite{CAII} that $\theta$-invariant triangulations of $\Poly_{2n+2}$ are in bijection with the clusters of cluster algebras of type $B_n$ and $C_n$. Furthermore, cluster variables correspond to the orbits of the action of $\theta$ on the diagonals of $\Poly_{2n+2}$, which can be either diameters or pairs of centrally symmetric non diameter diagonals. 

In this paper, given a $\theta$-invariant triangulation $T$, we define cluster algebras $\mathcal{A}^B(T)$ of type $B_n$, and $\mathcal{A}^C(T)$ of type $C_n$, with principal coefficients in $T$ (cf. Definition \ref{ca_def}), and we find an expansion formula for the cluster variable $x_{ab}$ corresponding to the $\theta$-orbit $[a,b]$ of the diagonal $(a,b)$ which connects the vertices $a$ and $b$. The formula we present is given in a combinatorial way. On the one hand, it expresses each cluster variable of type $B_n$ and $C_n$ in terms of cluster variables of type $A_n$, on the other hand, it allows one to get its expansion in terms of the cluster variables of the initial seed. In particular, we give a combinatorial description of the $F$-polynomial $F_{ab}$ and the $\bold{g}$-vector $\bold{g}_{ab}$ of $x_{ab}$. 

To state the result we need to define a simple operation on sets $\mathcal{D}$ of diagonals, the \emph{restriction}, denoted by $\text{Res}(\mathcal{D})$, which consists essentially of taking the diagonals obtained after identifying $n$ particular vertices of the polygon, see Definition \ref{def_restriction}. For a diagonal $\gamma$ of $\Poly_{n+3}$, we denote by $F_\gamma$ the $F$-polynomial of the cluster variable of type $A_n$ which corresponds to $\gamma$ in the cluster algebra with principal coefficients in the triangulation $\text{Res}(T)$ (cf. Definition \ref{def_ca_triang}). $F_\gamma$ has an explicit description, for example in terms of perfect matchings of the snake graph associated with $\gamma$. See \cite{MS,CSI} for details.

\begin{teo*} [\ref{theorem1}]
     Let $\mathcal{A}^B(T)$ be the cluster algebra of type $B_n$ with principal coefficients in a $\theta$-invariant
     triangulation $T=\{ \tau_1, \dots, \tau_{2n-1} \}$ of $\Poly_{2n+2}$. Then the $F$-polynomial $F_{ab}$ of $x_{ab}$ is given in the following way:

                       \begin{itemize}
                           \item[(i)] if $\text{Res}([a,b])$ contains only one diagonal $\gamma$, $ F_{ab}=F_\gamma$;
                           \item[(ii)] otherwise, $\text{Res}([a,b])=\{ \gamma_1, \gamma_2\}$, and
                           \begin{center}

                           $F_{ab}=F_{\gamma_1}F_{\gamma_2}- \bold{y}^cF_{(a,\theta(b))}$, \end{center}
 where $c \in \{0,1\}^n$ is such that $c_i=1$ if and only if the elementary lamination associated to $\tau_i$ crosses both $\gamma_1$ and $\gamma_2$, $i=1,\dots,n$.
                             \end{itemize}

\end{teo*}
We have analogous results for the $\bold{g}$-vector of $\mathcal{A}^B(T)$, and for the $F$-polynomial and the $\bold{g}$-vector of $\mathcal{A}^C(T)$. See Theorem \ref{theorem1} and Theorem \ref{theorem 2}.

Another cluster expansion formula for cluster algebras of type $B$ and $C$ has been given by Musiker in \cite{M} in terms of perfect matchings of certain labeled modified snake graphs. This formula holds only for the initial bipartite seed. In \cite{5} we use the results of the present paper to extend the work of Musiker to every seed. 

Moreover, Nakanishi and Stella provide in \cite{NS} a diagrammatic description of the $\bold{g}$-vectors of cluster algebras of type $B$ and $C$, while Reading studies them in \cite{R} using ring homomorphisms between cluster algebras of type $B$ and $C$, and cluster algebras of type $A$, induced by the fact that exchange matrices of type $B_n$ and $C_n$ ``dominate'' exchange matrices of type $A_n$. Furthermore, a cluster algebra of type $B_n$ (resp. $C_n$) can be realized as a disk with one orbifold point of weight 2 (resp. $\frac{1}{2}$), and $n+1$ boundary marked points \cite{FeST}. In \cite{FeliksonTumarkin}, Felikson and Tumarkin compute $\bold{g}$-vectors for cluster algebras from orbifolds, including type $B$ and $C$, in terms of laminations on the orbifolds. Finally, a relation between skew-symmetric and skew-symmetrizable cluster algebras has been investigated in \cite{FST_finite_type,Dupont} via folding.

\vspace{0.5cm}

On the other hand, the representation theory of symmetric quivers was developed by Derksen and Weyman in \cite{DW}, as well as Boos and Cerulli Irelli in \cite{boos2021degenerations}. A $symmetric$ $quiver$ is a quiver $Q$ with an involution $\sigma$ of vertices and arrows which reverses the orientation of arrows. A $symmetric$ $representation$ is an ordinary representation equipped with some extra data that forces each dual pair $(\alpha,\sigma(\alpha))$ of arrows of $Q$ to act anti-adjointly, see Section \ref{symmetric_quivers}. Symmetric representations are of two types: $orthogonal$ and $symplectic$. They form an additive category which is not abelian. Moreover, it was shown in \cite{DW,boos2021degenerations} that every indecomposable symmetric representation $M$ is uniquely determined by the $\nabla$-orbit of an indecomposable (ordinary) representation $L$ in one of the following forms:
\begin{center}
    \begin{itemize}
        \item [(I)] $M=L$ for $L\cong \nabla L$ ($indecomposable$ type);
        \item [(S)] $M=L\oplus \nabla L$ for $L \ncong \nabla L$ ($split$ type);
        \item [(R)] $M=L \oplus \nabla L$ for $L \cong \nabla L$ ($ramified$ type).
    \end{itemize}
\end{center}
Conversely, every indecomposable representation $L$ gives rise to exactly one of these indecomposable symmetric representations.

Derksen and Weyman in \cite{DW} stated the correspondence between positive roots of a root system of type $B_n$ (resp. $C_n$) and orthogonal (resp. symplectic) indecomposable representations of symmetric quivers of type $A_{2n-1}$. On the other hand, from the classification of finite type cluster algebras \cite{CAII}, we know that positive roots of type $B_n$ and $C_n$ correspond to non-initial cluster variables of type $B_n$ and $C_n$. Therefore, there is a one-to-one correspondence between non-initial cluster variables of type $B_n$ (resp. $C_n$) and orthogonal (resp. symplectic) indecomposable representations of symmetric quivers of type $A_{2n-1}$. The second objective of this work is to find explicitly this bijection. In the process of doing this, we extend it to any symmetric quiver in the mutation class of a symmetric quiver of type $A_{2n-1}$. 

Let $T$ be a $\theta$-invariant triangulation of $\Poly_{2n+2}$ with oriented diameter $d$. The quiver naturally associated to it (see Section \ref{quiver_to_triang}) is not symmetric. In order to get a symmetric quiver, we apply to the polygon an involution that we call $F_d$. It consists of cutting $\Poly_{2n+2}$ along $d$, then reflecting the right part with respect to the axis of symmetry of $d$, and finally gluing it again along $d$. $F_d$ induces an action on isotopy classes of diagonals of the polygon. Let $\rho$ denote the reflection of the polygon along $d$. Under the bijection $F_d$, $\theta$-orbits correspond to $\rho$-orbits. In particular, diameters correspond to $\rho$-invariant diagonals, while pairs of centrally symmetric diagonals correspond to $\rho$-invariant pairs of diagonals which are not orthogonal to $d$. Let $T'$ be the element in the isotopy class of $F_d(T)$ which is also a triangulation. Then $T'$ is a $\rho$-invariant triangulation of $\Poly_{2n+2}$. Therefore, the quiver $Q(T')$ associated to $T'$ is a cluster-tilted bound symmetric quiver of type $A_{2n-1}$ (\cite{S_Quiver_rep}, 3.4.1). Furthermore, indecomposable representations $L_\gamma$ of $Q(T')$ correspond to diagonals $\gamma$ of $\Poly_{2n+2}$ which are not in $T'$, and indecomposable symmetric representations correspond to their $\rho$-orbits.

In particular, let $\mathcal{A}^B (T)$ be the cluster algebra of type $B_n$ with principal coefficients in $T$. Let $[a,b]$ be a $\theta$-orbit and let $x_{ab}$ be the cluster variable which corresponds to $[a,b]$. If $F_d([a,b])=\{\alpha\}$ consists of only one $\rho$-invariant diagonal, then $x_{ab}$ corresponds to the orthogonal indecomposable $Q(T')$-representation $L_\alpha$ of type I. Otherwise, $F_d([a,b])=\{\alpha_1,\alpha_2\}$, and $L_{\alpha_2}=\nabla L_{\alpha_1}$. In this case, $x_{ab}$ corresponds to the orthogonal indecomposable $Q(T')$-representation  $L_{\alpha_1} \oplus \nabla L_{\alpha_1}$ of type S. 

On the other hand, for $\mathcal{A}^C(T)$, if $F_d([a,b])=\{\alpha\}$ consists of only one $\rho$-invariant diagonal, then $x_{ab}$ corresponds to the symplectic indecomposable $Q(T')$-representation $L_\alpha \oplus L_\alpha$ of type R. Otherwise, $F_d([a,b])=\{\alpha_1,\alpha_2\}$ with $L_{\alpha_2}=\nabla L_{\alpha_1}$. As before, $x_{ab}$ corresponds to the symplectic indecomposable $Q(T')$-representation $L_{\alpha_1} \oplus \nabla L_{\alpha_1}$ of type S.

Formulas of Theorem \ref{theorem1} (type $B$) and Theorem \ref{theorem 2} (type $C$) give the cluster expansion of each cluster variable associated to a $\theta$-orbit, on the one hand in terms of the cluster variables of the initial seed, on the other hand in terms of cluster variables of type $A_n$. It follows from the above correspondence that, given a cluster-tilted bound symmetric quiver $Q$ of type $A_{2n-1}$, they allow us to express the type $B_n$ (resp. type $C_n$) cluster variable that corresponds to an orthogonal (resp. symplectic) indecomposable representation of $Q$, on the one hand in terms of the initial cluster variables, on the other hand in terms of (ordinary) representations of $Q$. In other words, we get a Caldero-Chapoton like map (see \cite{CC}) from the category of symmetric representations of cluster tilted bound symmetric quivers of type $A_{2n-1}$ to cluster algebras of type $B_n$ and $C_n$.

This approach could be used to produce a categorification of other classes of non skew-symmetric cluster algebras through the representation theory of symmetric quivers. For example, they could provide an alternative categorification of non skew-symmetric cluster algebras associated by Felikson, Shapiro and Tumarkin \cite{FeST} to surfaces with marked points and order-2 orbifold points. These algebras have been categorified in the work of Geuenich and Labardini-Fragoso \cite{GeuLF1,GeuLF2} by species with potential.

To conclude, we give a categorical interpretation of Theorem \ref{theorem1} in the case where $Q(T')$ has no oriented cycles. To do this, we use the cluster multiplication formula of \cite{CEHR}. If $M$ is an orthogonal indecomposable representation of $Q(T')$, we denote by $\text{Res}(M)$ the representation of $Q(T')$ which corresponds to the restriction of the $\theta$-orbit corresponding to $M$. Moreover, we denote by $F_M$ the $F$-polynomial of the cluster variable $x_M$ of $\mathcal{A}^B(T)$ corresponding to $M$, and by $F_{\text{Res}(M)}$ the $F$-polynomial of the $Q(T')$-representation $\text{Res}(M)$ (see Section \ref{cat_interpretation}).

\begin{teo*}[\ref{cat_interpr}]
Let $M$ be an orthogonal indecomposable $Q(T')$-representation. If $\text{Res}(M)$ is indecomposable as $Q(T')$-representation, then 
    \begin{equation}
        \text{$F_M=F_{\text{Res}(M)}$.}
    \end{equation}
    
      Otherwise, $M=L\oplus \nabla L$ with $\text{dim Ext}^1(\nabla L, L)=1$, and there exists a non-split short exact sequence
    \begin{center}
        $0 \to L \to G_1 \oplus G_2 \to \nabla L \to 0 $,
    \end{center}
    where $G_1$ and $G_2$ are $\nabla$-invariant $Q(T')$-representations of type I. Then
   \begin{equation}
        F_{M}= F_{\text{Res}(M)} - \bold{y}^{\text{Res}(\textbf{dim}\nabla L^L)}F_{\text{Res}(L_{\nabla L}\oplus \nabla L / \nabla L^L)},
    \end{equation}
    %\begin{equation}
        %=F_{\text{Res}(G_1)}F_{\text{Res}(G_2)}-\bold{y}^{\text{Res}(\text{dim}\nabla L^L)}F_{\text{Res}(L_{\nabla L}\oplus \nabla L / \nabla L^L)},
    %\end{equation}
  where $L_{\nabla L} = \text{ker} (L \to \tau(\nabla L))$, $\nabla L^L=\text{im}(\tau^{-1}(L)\to \nabla L)$, with $\tau$ the Auslander-Reiten translation.
\end{teo*}

In literature there are other different categorifications of cluster algebras of type $B$ and $C$. In \cite{GLS1} Geiss, Leclerc and Schr\"{o}er use categories of locally free modules over certain Iwanaga-Gorenstein algebras; in \cite{Demonet} Demonet uses exact stably 2-Calabi-Yau categories endowed with the action of a finite group; in \cite{GeuLF1,GeuLF2} Geuenich and Labardini-Fragoso use species with potential.

The paper is organized as follows. Section \ref{s_background} is devoted to a quick overview of cluster algebras of geometric type with a particular focus on the geometric model for cluster algebras of type $A_n$, $B_n$ and $C_n$, that will be used throughout the paper. In Section \ref{s_formulas}, we give the definition of cluster algebras of type $B$ and $C$ with principal coefficient in a $\theta$-invariant triangulation of the polygon. Moreover, we state and prove the cluster expansion formulas for these algebras. Finally, in Section \ref{s_categ}, after a recollection on symmetric representation theory, we establish a correspondence between orthogonal (resp. symplectic) indecomposable representations of cluster-tilted bound symmetric quivers of type $A_{2n-1}$ and cluster variables of type $B_n$ (resp. $C_n$). Moreover, we give a categorical interpretation of Theorem \ref{theorem1}.
\section{Background}\label{s_background}

\subsection{Cluster algebras of geometric type}
Cluster algebras, introduced by Fomin and Zelevinsky in \cite{CAI}, are commutative algebras with a distinguished set of generators, the $cluster$ $variables$. Cluster variables are grouped into overlapping sets of constant cardinality $n$, the $clusters$, and the integer $n$ is called the rank of the cluster algebra. They are obtained combinatorially starting from an initial cluster $\bold{u}$, together with an integer $n\times n$ $exchange$ $matrix$ $B = (b_{ij})$ with the property that there exists a $symmetrizer$ $D=\text{diag}(d_1,\dots,d_n)$, with $d_i \in \mathbb{Z}_{>0}$ such that $DB$ is skew-symmetric, i.e. $B$ is skew symmetrizable, and a coefficient vector $\bold{y} = (y_i)$, whose entries are elements of a torsion-free abelian group $\mathbb{P}$. The triple $\Sigma=(\bold{u},\bold{y},B)$ is called the $initial$ $seed$. 
The set of cluster variables is obtained by repeated applications of the so called $mutations$ to the initial seed. To be more precise, let $\mathcal{F}$ be the field of rational functions in the indeterminates $u_1,\dots,u_n$ over the quotient field of the integer group ring $\mathbb{Z}\mathbb{P}$. Thus $\bold{u} = \{u_1,\dots, u_n\}$ is a transcendence basis for $\mathcal{F}$. For every $1\leq k \leq n$, the mutation $\mu_k(\bold{u})$ of the cluster $\bold{u}=\{u_1,\dots,u_n\}$ is a new cluster $\mu_k(\bold{u})=\bold{u} \setminus \{u_k\} \cup \{u_{k}'\}$ obtained from $\bold{u}$ by replacing the cluster variable $u_k$ by the new cluster variable $u_{k}'$ such that
\begin{equation}\label{intro_ex_rel}
    u_ku_k'= p_k^+ \displaystyle\prod_{b_{ik} >0 }u_i^{b_{ik}}+p_k^- \displaystyle\prod_{b_{ik}<0}u_i^{-b_{ik}}
\end{equation}
in $\mathcal{F}$, where $p_k^+$, $p_k^-$ are certain monomials in $y_1,\dots,y_n$. Equation \ref{intro_ex_rel} is the $exchange$ $relation$ between the cluster variables $u_k$ and $u_k'$. Each mutation also changes the coefficient vector $\bold{y}$, as well as the attached matrix $B$, but it does not change the symmetrizer which is the same for any matrix in the mutation class of $B$ (\cite{CAI}, Proposition 4.5). This combinatorics is encoded in the $cluster$ $complex$, which is the simplicial complex whose maximal faces are the clusters, and whose edges correspond to mutations.

The set $\mathcal{X}$ of all cluster variables is the union of all clusters obtained from the initial cluster $\bold{u}$ by repeated mutations. The $cluster$ $algebra$ $\mathcal{A}(\bold{u},\bold{y},B)$ is defined as the $\mathbb{ZP}$-subalgebra of $\mathcal{F}$ generated by $\mathcal{X}$. A cluster algebra is said to be of $finite$ $type$ if it has a finite number of cluster variables. Cluster algebras of finite type are classified by Dynkin diagrams, in the same way as semisimple Lie algebras and finite root systems \cite{CAII}.  

It is clear from the construction that every cluster variable is a rational function in the initial cluster variables $u_1,\dots,u_n$. In \cite{CAI} it is shown that every cluster variable $x$ is actually a Laurent polynomial in the $u_i$, that is, $x$ can be written as a reduced fraction
\begin{equation}\label{cluster_exp}
    x=\frac{f(u_1,\dots,u_n)}{\displaystyle\prod_{i=1}^n u_i^{d_i}},
\end{equation}
where $f \in \mathbb{Z}\mathbb{P}[u_1,\dots u_n]$ and $d_i \in \mathbb{Z}_{\geq 0}$. This is known as the $Laurent$ $phenomenon$. The right hand side of equation \ref{cluster_exp} is
called the cluster expansion of $x$ in $\bold{u}$.

The cluster algebra $\mathcal{A}(\bold{u},\bold{y},B)$ is determined by the initial matrix $B$ and the choice of a coefficient vector. If the coefficient group $\mathbb{P}$ is chosen to be the free abelian group on $m$ generators $y_1,\dots,y_m$, then the cluster algebra is said $of$ $geometric$ $type$. If $\Sigma=(\bold{x},\bold{y},B)$ is a seed of a cluster algebra of geometric type, then the datum of the pair $(\bold{y},B)$ is equivalent to the datum of an $extended$ $exchange$ $matrix$ $\Tilde{B}$, i.e. an $m \times n$ matrix whose top square matrix is $B$, and such that coefficient vectors can be recovered from the bottom part. 
A canonical choice in this setting is the $principal$ $coefficient$ $system$, introduced in \cite{CAIV}, which means that the coefficient group $\mathbb{P}$ is the free abelian group on $n$ generators $y_1,\dots,y_n$, and the initial coefficient tuple $\bold{y}=(y_1,\dots,y_n)$ consists of these $n$ generators. This is equivalent to taking in the initial seed the extended exchange matrix $\Tilde{B}=\left[\begin{matrix}
   B\\ I
\end{matrix} \right]$, where $I$ is the $n \times n$ identity matrix. The columns of the bottom part of the extended exchange matrices of any seed are called $c$-$vectors$. In \cite{CAIV}, the authors show that knowing the expansion formulas in the case where the cluster algebra has principal coefficients allows one to compute the expansion formulas for arbitrary coefficients. Moreover with this choice of coefficients, for each cluster variable $x$, a polynomial $F_x \in \mathbb{Z}[y_1,\dots,y_n]$ and an integer vector $\bold{g}_x \in \mathbb{Z}^n$ are defined. $F_x$ is called the $F$-$polynomial$ of $x$, and it is obtained by setting all $u_i=1$ in $x$. On the other hand, $\bold{g}_x$ is called the $g$-$vector$ of $x$, and it is the multi-degree of $x$ with respect to the $\mathbb{Z}^n$-grading in $\mathbb{Z}[u_1^{\pm 1},\dots,u_n^{\pm 1},y_1,\dots,y_n]$ given by $\textbf{deg}(u_i)=\bold{e}_i$ and $\textbf{deg}(y_j)=-\bold{b}_j$, where $\bold{e}_i$ is the standard basis vector of $\mathbb{Z}^n$ and $\bold{b}_j$ is the $j$-th column of $B$. Knowing the cluster expansion of $x$ in $\bold{u}$ is equivalent to knowing $F_x$ and $\bold{g}_x$. In fact,
\begin{equation}
    x=F_x(\hat{y}_1,\dots,\hat{y}_n)\bold{u}^{\bold{g}_x},
\end{equation}

where $\hat{y}_j=y_j\displaystyle\prod_{i=1}^n u_i^{b_{ij}}$, and $\bold{u}^{\bold{g}_x}$ is the monomial $u_1^{g_1}\cdots u_n^{g_n}$, if $\bold{g}_x=(g_1,\dots,g_n)$.

Fomin, Shapiro and Thurston in \cite{FST,FT}, and Labardini-Fragoso in \cite{LF}, initiated the study of cluster algebras, and quivers with potential, arising from triangulations of surfaces with boundary and marked points. In their approach, cluster variables correspond to arcs in the surface, and clusters correspond to triangulations. Musiker and Schiffler in \cite{MS}, and Musiker, Schiffler and Williams in \cite{MSW11}, gave an expansion formula for the cluster variables in terms of perfect matchings of some labeled graphs, called \emph{snake graphs}, that are recursively constructed from the surface by gluing together elementary pieces called tiles.

\subsection{Combinatorial description of the cluster complex of type $A_n$}
\label{sec:ca-type-A}
In this section, we recall the geometric model for cluster algebras of type $A$. 

Let $n$ be a positive integer. Let $\Poly_{n+3}$ be the regular polygon with $n+3$ vertices. Fomin and Zelevinsky show in \cite{Y-systems,CAII} that clusters of a cluster algebra of type $A_n$ are in bijection with triangulations of 
$\Poly_{n+3}$, i.e., maximal collections of non-crossing diagonals, and cluster variables correspond to diagonals. Moreover, mutations correspond to flips, so two triangulations are joined by an edge in the exchange graph if and only if they are obtained from each other by  replacing a diagonal in a quadrilateral formed by two triangles of the triangulation by the another diagonal of the same quadrilateral. Furthermore, the exchange matrix of the seed whose cluster corresponds to a triangulation $\Bar{T}=\{ \tau_1, \dots, \tau_n \}$ of $\mathbf{P}_{n+3}$ is given by the skew-symmetric $n\times n$ matrix $B(\Bar{T})=(b_{ij}(\Bar{T}))$ defined by:

\begin{equation}
\label{eq:B(T)}
b_{ij}(\Bar{T})=
\begin{cases}
1 & \!\!\text{if $\tau_i$ and $\tau_j$ are two sides of a triangle
  in~$\Bar{T}$,}\\
& \hspace{-0.5cm}\text{\quad with $\tau_i$ following $\tau_j$ in  counterclockwise order;}\\
-1 & \!\!\text{if $\tau_i$ and $\tau_j$ label two sides of a triangle
  in~$\Bar{T}$,}\\
& \hspace{-0.5cm}\text{\quad with $\tau_j$ following $\tau_i$ in  counterclockwise order;}\\
0 & \!\!\text{if $\tau_i$ and $\tau_j$ do not belong to the same triangle
  in~$\Bar{T}$.}
\end{cases}
\end{equation}

 Let $(a,b)$ denote the diagonal which connects vertices $a$ and $b$ of $\Poly_{n+3}$. We indicate by $x_{(a,b)}$ the cluster variable corresponding to $(a,b)$, with the convention that $x_{(a,b)} = 1$ if $a$ and $b$ are two consecutive vertices of $\Poly_{n+3}$. 
Hence the exchange relations in a cluster 
algebra of type $A_n$ have the form 
\begin{equation} 
\label{eq:exchange-relation-A} 
x_{(a,b)} x_{(c,d)} = p^+_{ab,cd}\, x_{(a,d)}\, x_{(b,c)} 
+ p^-_{ab,cd}\, x_{(a,c)}\, x_{(b,d)} \ , 
\end{equation} 
where $a,c,b,d$ are any four vertices of $\Poly_{n+3}$ taken in 
counter-clockwise order, and 
%$\langle x_{ac} x_{bd} : x_{ab} x_{cd} \rangle$ 
%and $\langle x_{ac} x_{bd} : x_{ad} x_{bc} \rangle$ 
$p^\pm_{ab,cd}$ 
are elements of the coefficient semifield~$\PP$. 
See Figure~\ref{fig_prin_coeff}. 
%\end{proposition} 
 
%\begin{proof} 
%Compare (\ref{eq:exch-rel-fin-type}) with \cite[Lemma~4.2]{cfz}. 
%\end{proof} 
  
\begin{definition}\label{def_ca_triang}
  Let $\Bar{T}$ be a triangulation of $\Poly_{n+3}$. Let $\bold{u}_{\Bar{T}}=\{u_1,\dots,u_n\}$ be the cluster associated to $\Bar{T}$, and $\bold{y}_{\Bar{T}}=(y_1,\dots,y_n)$ be the initial coefficient vector of generators of $\mathbb{P}=\mathrm{Trop}(y_1,\dots,y_n)$. Then $\mathcal{A}^A(\Bar{T}):=\mathcal{A}(\bold{u}_{\Bar{T}},\bold{y}_{\Bar{T}},B(\Bar{T}))$ is called the \emph{cluster algebra of type $A_n$ with principal coefficients in $\Bar{T}$}.   
\end{definition}

In this case the coefficients $p^\pm_{ab,cd}$ can be explicitly determined from $\Bar{T}$. The following definition and proposition are just a restatement of Definition 17.2 and Proposition 17.3 of \cite{FT} in the case of diagonals of a polygon. 

\begin{definition} \label{elem-laminate}
Let $\gamma=(a,b)$ be a diagonal of $\Poly_{n+3}$. The \emph{elementary lamination} associated to $\gamma$ is the segment $L_{\gamma}$ which 
 begins at a point $a'\in \Poly$ located near $a$ in
the clockwise direction, and ends at a point $b' \in \Poly$ near $b$ in
the clockwise direction. % on the left.
If $\Bar{T}=\{ \tau_1,\dots,\tau_n \}$ is a triangulation of $\Poly_{n+3}$, then we 
let $L_i$ denote $L_{\tau_i}$.
\end{definition}

 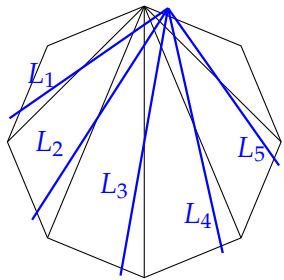
\begin{figure}[H]
        \centering
\begin{tikzpicture}[scale=0.6]
     \draw (90:3cm) -- (135:3cm) -- (180:3cm) -- (225:3cm) -- (270:3cm) -- (315:3cm) -- (360:3cm) -- (45:3cm) -- cycle;
                    \draw (90:3cm) -- node[near end, below, xshift=-1mm, yshift=-1mm] {} (180:3cm);
                    \draw (90:3cm) -- node[near end,below,xshift=1.2mm] {} (225:3cm);
                    \draw (90:3cm) -- node[near end, below,xshift=1mm] {} (270:3cm);
                    \draw (90:3cm) -- node[near end, above, xshift=-1mm] {} (315:3cm);
                    \draw (90:3cm) -- node[near end, above, xshift=-1mm] {} (360:3cm);

                    \draw[blue, line width=0.3mm] (80:3cm) -- node[near end, above, xshift=-1mm, yshift=-1mm] {$L_1$} (170:3cm);
                    \draw[blue, line width=0.3mm] (80:3cm) -- node[near end, above left,xshift=1.2mm] {$L_2$} (215:3cm);
                    \draw[blue, line width=0.3mm] (80:3cm) -- node[near end, above left,xshift=1mm] {$L_3$} (260:3cm);
                    \draw[blue, line width=0.3mm] (80:3cm) -- node[near end, below, xshift=-1.5mm] {$L_4$} (305:3cm);
                    \draw[blue, line width=0.3mm] (80:3cm) -- node[near end, below] {$L_5$} (350:3cm);

\end{tikzpicture}
        \caption{A triangulated octagon with the elementary lamination associated to each diagonal of the triangulation (in blue).}
        \label{lamination}
    \end{figure}

\begin{proposition} \label{up:skein1}
Let $\mathcal{A}^A(\Bar{T})$ be a cluster algebra of type $A_n$ with principal coefficients in a triangulation $\Bar{T}=\{ \tau_1, \dots, \tau_n \}$ of $\Poly_{n+3}$.
Let $(a,b)$ and $(c,d)$ be two diagonals
which intersect each other. Then
\begin{equation} \label{u:skein-eq1}
x_{(a,b)} x_{(c,d)} = \bold{y}^{\bold{d}_{ac,bd}} x_{(a,d)} ~x_{(b,c)}
 + \bold{y}^{\bold{d}_{ad,bc}}x_{(a,c)} ~x_{(b,d)},
\end{equation}
where
$\bold{d}_{ac,bd}$ (resp., $\bold{d}_{ad,bc}$) is the vector whose $i$-th coordinate is 1 if $L_i$ crosses both $(a,c)$ and $(b,d)$ (resp., $(a,d)$ and $(b,c)$); 0 otherwise.
\end{proposition}
\begin{example}
Let $\mathcal{A}^A(\Bar{T})$ be the cluster algebra of type $A_5$ with principal coefficients in the triangulation of the octagon in Figure \ref{fig_prin_coeff}.
    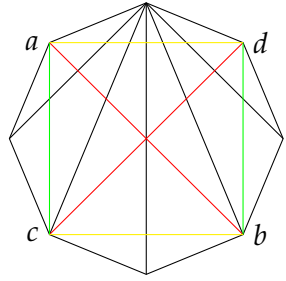
\begin{figure}[H]
        \centering
\begin{tikzpicture}[scale=0.6]
     \draw (90:3cm) -- (135:3cm) -- (180:3cm) -- (225:3cm) -- (270:3cm) -- (315:3cm) -- (360:3cm) -- (45:3cm) -- cycle;
                    \draw (90:3cm) -- node[near end, below, xshift=-1mm, yshift=-1mm] {} (180:3cm);
                    \draw (90:3cm) -- node[near end,below,xshift=1.2mm] {} (225:3cm);
                    \draw (90:3cm) -- node[near end, below,xshift=1mm] {} (270:3cm);
                    \draw (90:3cm) -- node[near end, above, xshift=-1mm] {} (315:3cm);
                    \draw (90:3cm) -- node[near end, above, xshift=-1mm] {} (360:3cm);

                   \node at (135:3cm) [left] {$a$};
                   \node at (225:3cm) [left] {$c$};
                   \node at (315:3cm) [right] {$b$};
                   \node at (45:3cm) [right] {$d$};
                   \draw[red] (135:3cm) -- (315:3cm);
                   \draw[red] (225:3cm) -- (45:3cm);
                   \draw[green] (135:3cm) -- (225:3cm);
                   \draw[green] (315:3cm) -- (45:3cm);
                   \draw[yellow] (135:3cm) -- (45:3cm);
                   \draw[yellow] (225:3cm) -- (315:3cm);

\end{tikzpicture}
\caption{An exchange relation in a triangulated octagon.}\label{fig_prin_coeff}        
    \end{figure}
By Proposition \ref{up:skein1}, $x_{(a,b)}x_{(c,d)}=x_{(a,d)}x_{(b,c)}+y_3y_4x_{(a,c)}x_{(b,d)}$.
\end{example}

\subsection{Combinatorial description of the cluster complex of type $B_n$/$C_n$}
\label{sec:ca-type-C} 
In this section, we recall the geometric model for cluster algebras of types $B$ and $C$.

Let $n$ be a positive integer. Let $\Poly_{2n+2}$ be the regular polygon with $2n+2$ vertices. Let $\theta$ denote the $180^\circ$ rotation of $\Poly_{2n+2}$. There is a natural action of $\theta$ 
on the diagonals of $\Poly_{2n+2}$. Each orbit of this action is 
either a diameter (i.e., a diagonal connecting antipodal vertices) 
or an unordered pair of centrally symmetric non-diameter diagonals 
of $\Poly_{2n+2}$.

Fomin and Zelevinsky show in \cite{Y-systems,CAII} that clusters of a cluster algebra of type $B_n$ or $C_n$ are in bijection with centrally-symmetric (that is, $\theta$-invariant) 
triangulations of $\Poly_{2n+2}$, and cluster variables correspond to $\theta$-orbits. 
Two centrally symmetric triangulations are joined by an edge in 
the exchange graph if and only if they are obtained from each other either by a flip 
involving two diameters, or by a pair of centrally symmetric flips. 
 
For a vertex $a$ of $\Poly_{2n+2}$, let $\Bar{a}$ denote the 
antipodal vertex $\theta (a)$. We indicate 
by $x_{ab}$ the cluster variable corresponding to the $\theta$-orbit $[a,b]$ of the diagonal $(a,b)$. 
Thus, we have $x_{ab} = x_{ba} = x_{\overline a \,\overline b}$, with the convention that $x_{ab} = 1$ if $a$ 
and $b$ are consecutive vertices in~$\Poly_{2n+2}$. 
 
They obtain the following concrete description of the exchange relations in 
types $B_n$ and~$C_n\,$. 
 
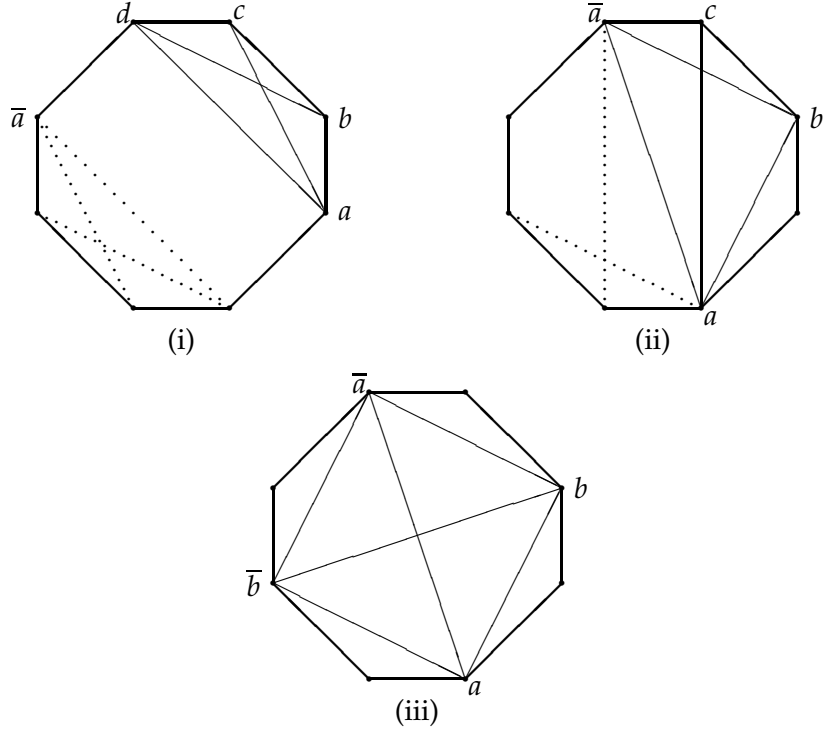
\begin{figure}[H] 
\begin{center} 
\setlength{\unitlength}{1.8pt} 
\begin{picture}(60,66)(0,-2) 
\thicklines 
  \multiput(0,20)(60,0){2}{\line(0,1){20}} 
  \multiput(20,0)(0,60){2}{\line(1,0){20}} 
  \multiput(0,40)(40,-40){2}{\line(1,1){20}} 
  \multiput(20,0)(40,40){2}{\line(-1,1){20}} 
 
  \multiput(20,0)(20,0){2}{\circle*{1}} 
  \multiput(20,60)(20,0){2}{\circle*{1}} 
  \multiput(0,20)(0,20){2}{\circle*{1}} 
  \multiput(60,20)(0,20){2}{\circle*{1}} 
 
\thinlines 
\put(20,60){\line(2,-1){40}} 
\put(20,60){\line(1,-1){40}} 
\put(40,60){\line(1,-2){20}} 
 
\put(42,62){\makebox(0,0){$c$}} 
\put(64,20){\makebox(0,0){$a$}} 
\put(64,40){\makebox(0,0){$b$}} 
\put(18,62){\makebox(0,0){$d$}} 
\put(-4,40){\makebox(0,0){$\overline a$}} 
%\put(-4,20){\makebox(0,0){$Q'$}} 

%\put(36,49){\makebox(0,0){$\alpha$}} 
%\put(47,39){\makebox(0,0){$\alpha'$}} 
%\put(30,64){\makebox(0,0){$\beta_1$}} 
%\put(64,30){\makebox(0,0){$\beta_3$}} 
%\put(37,37){\makebox(0,0){$\beta_2$}} 
%\put(53,53){\makebox(0,0){$\beta_4$}} 
 
\multiput(20,0)(-1,2){20}{\circle*{0.5}} 
\multiput(0,20)(2,-1){20}{\circle*{0.5}} 
\multiput(0,40)(2,-2){20}{\circle*{0.5}} 
 
\put(30,-7){\makebox(0,0){(i)}} 
 
\end{picture} 
\qquad\qquad\qquad 
\begin{picture}(60,66)(0,-2) 
\thicklines 
  \multiput(0,20)(60,0){2}{\line(0,1){20}} 
  \multiput(20,0)(0,60){2}{\line(1,0){20}} 
  \multiput(0,40)(40,-40){2}{\line(1,1){20}} 
  \multiput(20,0)(40,40){2}{\line(-1,1){20}} 
 
  \multiput(20,0)(20,0){2}{\circle*{1}} 
  \multiput(20,60)(20,0){2}{\circle*{1}} 
  \multiput(0,20)(0,20){2}{\circle*{1}} 
  \multiput(60,20)(0,20){2}{\circle*{1}} 
 
\thinlines 
 
\put(40,0){\line(0,1){60}} 
\put(20,60){\line(2,-1){40}} 
\put(40,0){\line(1,2){20}}

\put(20,60){\line(1,-3){20}} 
%\multiput(19.4,59.8)(3,-9){7}{\line(1,-3){2}} 

\put(42,-2){\makebox(0,0){$a$}} 
\put(42,62){\makebox(0,0){$c$}} 
%\put(64,20){\makebox(0,0){$a$}} 
\put(64,40){\makebox(0,0){$b$}} 
\put(18,62){\makebox(0,0){$\overline a$}} 
%\put(-4,40){\makebox(0,0){$\overline a$}} 
%\put(-4,20){\makebox(0,0){$Q'$}} 
 
%\put(34,49){\makebox(0,0){$\alpha$}} 
%\put(44,34){\makebox(0,0){$\alpha'$}} 
%\put(30,64){\makebox(0,0){$\beta_4$}} 
%\put(32,37){\makebox(0,0){$\beta_1$}} 
%\put(28,23){\makebox(0,0){$\beta_1'$}} 
%\put(54,20){\makebox(0,0){$\beta_2$}} 
%\put(53,53){\makebox(0,0){$\beta_3$}} 
 
\multiput(20,0)(0,2){30}{\circle*{0.5}} 
\multiput(0,20)(2,-1){20}{\circle*{0.5}} 
 
\put(30,-7){\makebox(0,0){(ii)}} 
 
\end{picture} 
\\[.4in] 
\begin{picture}(60,66)(0,-7) 
\thicklines 
  \multiput(0,20)(60,0){2}{\line(0,1){20}} 
  \multiput(20,0)(0,60){2}{\line(1,0){20}} 
  \multiput(0,40)(40,-40){2}{\line(1,1){20}} 
  \multiput(20,0)(40,40){2}{\line(-1,1){20}} 
 
  \multiput(20,0)(20,0){2}{\circle*{1}} 
  \multiput(20,60)(20,0){2}{\circle*{1}} 
  \multiput(0,20)(0,20){2}{\circle*{1}} 
  \multiput(60,20)(0,20){2}{\circle*{1}} 
 
\thinlines 
\put(20,60){\line(2,-1){40}} 
\put(0,20){\line(2,-1){40}} 
\put(20,60){\line(1,-3){20}} 
\put(40,0){\line(1,2){20}} 
\put(0,20){\line(1,2){20}} 
\put(0,20){\line(3,1){60}} 
%\multiput(0,20)(9,3){7}{\line(3,1){6}} 
 
\put(42,-2){\makebox(0,0){$a$}} 
\put(-4,20){\makebox(0,0){$\overline b$}} 
%\put(64,20){\makebox(0,0){$a$}} 
\put(64,40){\makebox(0,0){$b$}} 
\put(18,62){\makebox(0,0){$\overline a$}} 
 
%\put(30,39){\makebox(0,0){$\alpha$}} 
%\put(42,32){\makebox(0,0){$\alpha'$}} 
%\put(54,20){\makebox(0,0){$\beta_2$}} 
%\put(38,47){\makebox(0,0){$\beta_1$}} 
 
\put(30,-7){\makebox(0,0){(iii)}} 
 
\end{picture} 
 
\end{center} 
\caption{%Illustrating 
Exchanges in types $B_n$ and $C_n$} 
\label{fig:quadrilateral-D} 
\end{figure}

\begin{proposition}(\cite{CAII}, Proposition 12.9) 
\label{pr:exchanges-C} 
The exchange relations in a cluster algebra of type $B_n$ ($r=1$) or $C_n$ ($r=2$)
have the following form: 
\begin{equation} 
\label{eq:exchange-relation-C-1} 
x_{ac} x_{bd} = 
%\langle x_{ac} x_{bd} : x_{ab} x_{cd} \rangle 
p^+_{ac,bd}\, 
x_{ab}\, x_{cd} 
+ 
%\langle x_{ac} x_{bd} : x_{ad} x_{bc} \rangle 
p^-_{ac,bd}\, 
x_{ad}\, x_{bc} \ , 
\end{equation} 
for some coefficients $p^+_{ac,bd}$ and $p^-_{ac,bd}$, whenever $a,b,c,d,\overline a$ are %five vertices of $\Poly_{2n+2}$ taken 
in counter-clockwise order; 
\begin{equation} 
\label{eq:exchange-relation-C-2} 
x_{ac} x_{a \overline b} = 
%\langle x_{ac} x_{bd} : x_{ab} x_{cd} \rangle 
p^+_{ac,a \overline b}\, 
x_{ab}\, x_{a \overline c} 
+ 
%\langle x_{ac} x_{bd} : x_{ad} x_{bc} \rangle 
p^-_{ac,a\overline b}\, 
x_{a \overline a}^{2/r} \, x_{bc} \ , 
\end{equation} 
for some coefficients $p^+_{ac,a \overline b}$ and $p^-_{ac,a\overline b}$, whenever $a,b,c,\overline a$ are %four vertices of $\Poly_{2n+2}$ %taken 
in counter-clockwise order; 
\begin{equation} 
\label{eq:exchange-relation-C-3} 
x_{a\overline a} x_{b \overline b} = 
%\langle x_{ac} x_{bd} : x_{ab} x_{cd} \rangle 
p^+_{a\overline a,b \overline b}\, 
x_{ab}^r 
+ 
%\langle x_{ac} x_{bd} : x_{ad} x_{bc} \rangle 
p^-_{a\overline a,b\overline b}\, 
x_{a \overline b}^r \ , 
\end{equation} 
for some coefficients $p^+_{a\overline a,b \overline b}$ and $p^-_{a\overline a,b\overline b}$, whenever $a,b,\overline a$ are %four vertices of $\Poly_{2n+2}$ %taken 
in counter-clockwise order. 
%$\langle x_{ac} x_{bd} : x_{ab} x_{cd} \rangle$ 
%and $\langle x_{ac} x_{bd} : x_{ad} x_{bc} \rangle$ 
%$p^\pm_{ac,bd}$ 
%are elements of the coefficient semifield~$\PP$. 
See Figure~\ref{fig:quadrilateral-D}. 
\end{proposition} 
 
\section{Cluster algebras of type $B$ and $C$ with principal coefficients}\label{s_formulas}
Let $T=\{ \tau_1, \dots, \tau_{2n-1} \}$ be a $\theta$-invariant triangulation of $\Poly_{2n+2}$. It follows that $T$ has $n-1$ pairs of centrally symmetric diagonals and exactly one diameter $d$. Assuming that $d$ is oriented, in this section we associate to $T$ a cluster algebra of type $B_n$ and $C_n$ with principal coefficients in the initial seed corresponding to $T$. 

\begin{definition}\label{def_restriction}
    Let $\mathcal{D}$ be a set of diagonals of $\Poly_{2n+2}$. We define the $restriction$ $of$ $\mathcal{D}$, and we denote it by $\text{Res}(\mathcal{D})$, as the set of diagonals of $\Poly_{n+3}$ obtained from those of $\mathcal{D}$ identifying all the vertices which lie on the right of $d$. 
\end{definition}
We use the label $\ast$ for the vertex of $\Poly_{n+3}$ which is obtained by identifying the vertices of $\Poly_{2n+2}$ which lie on the right of $d$.
\begin{definition}
    Let $v \in \mathbb{Z}_{\geq 0}^{2n-1}$. We define the $restriction$ $of$ $v$, and we denote it by $\text{Res}(v)$,  as the vector of the first $n$ coordinates of $v$.
\end{definition}

Let $\Bar{T}=\text{Res}(T)=\{\tau_1,\dots,\tau_{n-1},d\}$ be the triangulation of $\Poly_{n+3}$ which is obtained from $T$ identifying all the vertices of $\Poly_{2n+2}$ which lie on the right of $d$.
Let $B(\Bar{T})=(b_{ij})$ be the skew-symmetric $n \times n$ matrix associated to $\Bar{T}$ (cf. \ref{eq:B(T)}). So $b_{ij}=1$ if and only if $\tau_i$ and $\tau_j$ are sides of a triangle of $T$, and $\tau_i$ is counterclockwise from $\tau_j$. See Figure \ref{matrix_to_triang}. Let $D=\text{diag}(1,\dots,1,2)$ be the $n\times n$ diagonal matrix with diagonal entries $(1,\dots,1,2)$. Since the symmetrizer is constant in the mutation class of a matrix (\cite{CAI}, Proposition 4.5), then $DB(\Bar{T})$ is skew-symmetrizable of type $B$ and $B(\Bar{T})D$ is skew-symmetrizable of type $C$, according to the convention of \cite{CAII}. 

\begin{example}
Figure \ref{matrix_to_triang} illustrates how to compute the  $3\times 3$ skew-symmetric  matrix $B(\Bar{T})$ associated to the $\theta$-invariant triangulation $T$ of the octagon $\Poly_8$.    

       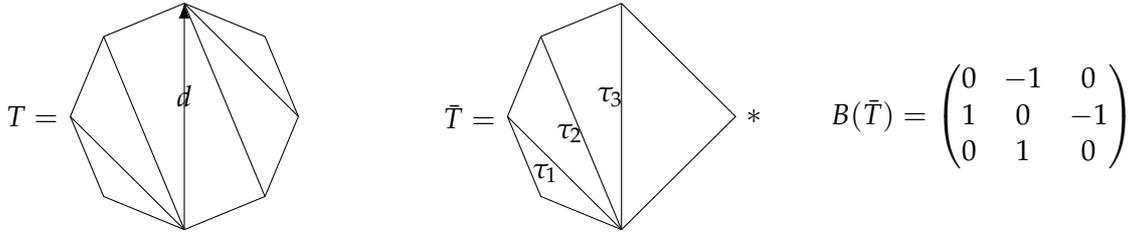
\begin{figure}[H]
        \centering
\begin{tikzpicture}[scale=0.5]
\node at (-4,0) {$T=$};                
     \draw (90:3cm) -- (135:3cm) -- (180:3cm) -- (225:3cm) -- (270:3cm) -- (315:3cm) -- (360:3cm) -- (45:3cm) -- cycle;
                    \draw (90:3cm) --(315:3cm);
                    \draw (90:3cm) -- (360:3cm);
                    \draw (270:3cm) -- (180:3cm);
                    \draw (270:3cm) -- (135:3cm);
                     \draw[-{Latex[length=2mm]}] (270:3cm) -- node[midway, above] {$d$}(90:3cm);
  \node at (7.5,0) {$\Bar{T}=$};                
\begin{scope}[xshift=11.5cm]
    \draw (90:3cm) -- (135:3cm) -- (180:3cm) -- (225:3cm) -- (270:3cm) -- (360:3cm) -- cycle;
                    %\draw (270:3cm) --(45:3cm);
                    %\draw (270:3cm) -- (360:3cm);
                    \draw (270:3cm) -- node[midway, left,xshift=0.5mm] {$\tau_1$}(180:3cm);
                    \draw (270:3cm) -- node[midway, left,xshift=1.5mm] {$\tau_2$}(135:3cm);
                     \draw (270:3cm) -- node[midway, above, xshift=-1.5mm] {$\tau_3$}(90:3cm);
                     \node[right] at (360:3cm) {$\ast$};  
\end{scope}

  \node at (21,0) {$B(\Bar{T})=\begin{pmatrix}
      0 & -1 & 0 \\ 1 & 0 & -1 \\ 0 & 1 & 0
  \end{pmatrix}$};

\end{tikzpicture}
        \caption{The matrix $B(\Bar{T})$ associated with a $\theta$-invariant triangulation of an octagon.}
        \label{matrix_to_triang}
    \end{figure}

Let $D=\text{diag}(1,1,2)$. Then the Cartan counterpart of $DB(\Bar{T})=\begin{pmatrix}
      0 & -1 & 0 \\ 1 & 0 & -1 \\ 0 & 2 & 0
  \end{pmatrix}$ is the Cartan matrix of type $B_3$, while the one of $B(\Bar{T})D=\begin{pmatrix}
      0 & -1 & 0 \\ 1 & 0 & -2 \\ 0 & 1 & 0
  \end{pmatrix}$ is the Cartan matrix of type $C_3$.
  \end{example}

  \begin{definition}\label{ca_def}
Let $T$ be a $\theta$-invariant triangulation of $\Poly_{2n+2}$. Let $\bold{u}_T=\{u_1,\dots, u_n\}$ be the cluster associated to $T$, and $\bold{y}_T=(y_1,\dots, y_n)$ be the initial coefficient vector of generators of $\mathbb{P}=\text{Trop}(y_1,\dots, y_n)$. Then $\mathcal{A}^B (T):=\mathcal{A}(\bold{u}_T,\bold{y}_T,DB(\Bar{T}))$ (resp. $\mathcal{A}^C(T):=\mathcal{A}(\bold{u}_T,\bold{y}_T,B(\Bar{T})D)$) is the $cluster$ $algebra$ $of$ $type$ $B_n$ (resp. $C_n$) $with$ $principal$ $coefficients$ $in$ $T$.
\end{definition}

\begin{remark}
    $\mathcal{A}^B (T)$ (resp. $\mathcal{A}^C(T)$), up to a change of coefficients, does not depend on $T$, but it depends only on $n$, since any two $\theta$-invariant triangulations of $\Poly_{2n+2}$ can be obtained from each other by a sequence of flips of diameters and pairs of centrally symmetric flips.
\end{remark}

\subsection{Cluster expansion formula for cluster algebras of type $B$ and $C$ }

Let $n$ be a positive integer. Let $\Poly_{2n+2}$ be the regular polygon with $2n+2$ vertices. Let $T=\{ \tau_1, \dots, \tau_n=d, \dots, \tau_{2n-1} \}$ be a $\theta$-invariant triangulation of $\Poly_{2n+2}$ with oriented diameter $d$, and let $\Bar{T}=\text{Res}(T)=\{ \tau_1, \dots, \tau_n=d \}$. Let $\mathcal{A}^B(T)$ (resp. $\mathcal{A}^C(T)$) be the cluster algebras of type $B_n$ (resp. $C_n$) with principal coefficients in $T$ (cf. Definition \ref{ca_def}), and let $\mathcal{A}^A(\Bar{T})$ be the cluster algebra of type $A_n$ with principal coefficients in $\Bar{T}$ (cf. Section \ref{sec:ca-type-A}). For a diagonal $\gamma$ of $\Poly_{n+3}$, let $F_\gamma$ and $\bold{g}_\gamma$ denote the $F$-polynomial and the $\bold{g}$-vector respectively of the cluster variable $x_\gamma \in \mathcal{A}^A(\Bar{T})$. They have an explicit description, for example in terms of perfect matchings of the snake graph associated with $\gamma$. See \cite{MS,CSI} for details.

In this section, we present a formula which expresses each cluster variable of $\mathcal{A}^B(T)$ and $\mathcal{A}^C(T)$ in terms of cluster variables of $\mathcal{A}^A(\Bar{T})$ (cf. Theorem \ref{theorem1} and Theorem \ref{theorem 2}).

\subsubsection{Type B}

\begin{definition}\label{def_type_B}
Let $[a,b] \not \subset T$ be an orbit of the action of $\theta$ on the diagonals of $\Poly_{2n+2}$. If $\text{Res}([a,b])$ contains only one diagonal $\gamma$ (as in Figure \ref{type B 1}) we define
  \begin{equation}
      F_{ab}^B=F_\gamma,
       \end{equation}
\begin{equation}
 \hspace{3cm}\bold{g}_{ab}^B=\begin{cases}
          \text{$Dg_\gamma$ \hspace{1cm}if $\gamma$ does not cross $d=\tau_n$;}\\
          \text{$Dg_\gamma+\bold{e}_n$ \hspace{0.1cm} if $\gamma$ crosses $d=\tau_n$}.
      \end{cases}
      \end{equation}
Otherwise $(a,b)$ crosses $d$, and $\text{Res}([a,b])=\{\gamma_1,\gamma_2\}$ (as in Figure \ref{type B}). We define
\begin{equation}\label{def_F^1}
    F_{ab}^B=F_{\gamma_1}F_{\gamma_2}- \bold{y}^{\bold{d}_{\gamma_1,\gamma_2}}F_{(a,\Bar{b})},
    \end{equation}
\begin{equation}\label{def_g^1}
    \hspace{-0.1cm}\bold{g}_{ab}^B=D(\bold{g}_{\gamma_1}+\bold{g}_{\gamma_2}+\bold{e}_n),
\end{equation}
with the notation of Proposition \ref{up:skein1}.

The definition is extended to any $\theta$-orbit by letting $F_{ab}^B=1$ and $\bold{g}_{ab}^B=\bold{e}_i$ if $[a,b]=\{ \tau_i,\tau_{2n-i} \} \in T$, and $F_{ab}^B=1$ and $\bold{g}_{ab}^B=\bold{0}$ if $(a,b)$ is a boundary edge of $\Poly_{2n+2}$.
   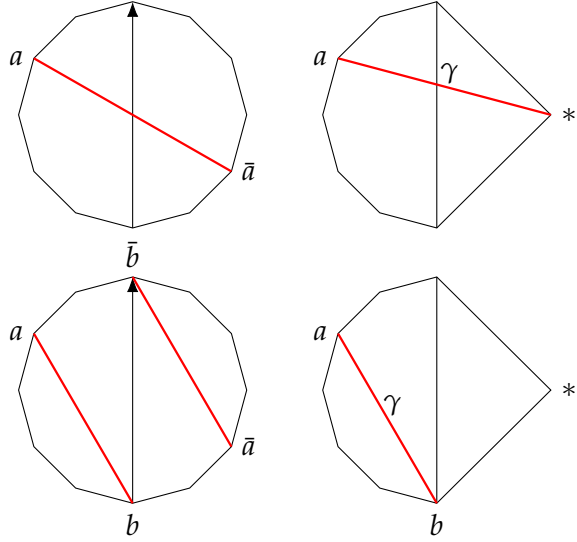
\begin{figure}[H]
                               \centering
                \begin{tikzpicture}[scale=0.5]
                    \draw (90:3cm) -- (120:3cm) -- (150:3cm) -- (180:3cm) -- (210:3cm) -- (240:3cm) -- (270:3cm) -- (300:3cm) -- (330:3cm) -- (360:3cm) -- (30:3cm) -- (60:3cm) --  cycle;
                  
                    \draw[-{Latex[length=2mm]}] (270:3cm) -- node[midway, above left,xshift=1mm] {} (90:3cm);
                  
                     \draw[red, line width=0.3mm] (150:3cm) -- (330:3cm);

                    %\node at (120:3cm) [left] {$a$};
                    %\node at (210:3cm) [left] {$\Bar{b}$};
                    \node at (330:3cm) [right] {$\Bar{a}$};
                    %\node at (60:3cm) [right] {$\hat{a}$};

\node at (150:3cm) [left] {$a$};
                    %\node at (240:3cm) [left] {$\hat{d}$};
                    %\node at (300:3cm) [right] {$d$};
                   % \node at (30:3cm) [right] {$b$};

                    \begin{scope}[xshift=8cm]
                      \draw (90:3cm) -- (120:3cm) -- (150:3cm) -- (180:3cm) -- (210:3cm) -- (240:3cm) -- (270:3cm) -- (360:3cm)  --  cycle;
                    %\draw (90:3cm) -- node[midway, above, xshift=-1mm, yshift=-1mm] {} (180:3cm);
                    %\draw (90:3cm) -- node[midway, above left,xshift=1.2mm] {} (225:3cm);
                    \draw (90:3cm) -- node[midway, above left,xshift=1mm] {} (270:3cm);
                   % \draw (90:3cm) -- node[midway, above right, xshift=-1mm] {} (315:3cm);
                    %\draw (90:3cm) -- node[midway, above right, xshift=-1mm] {} (360:3cm);
                    \draw[red, line width=0.3mm] (150:3cm) -- node[midway, above right, xshift=-2mm, yshift=-1mm] {$\textcolor{black}{\gamma}$} (360:3cm); 
                    %\draw[black, line width=0.3mm] (152:3cm) --  (358:3cm); 
                    % \draw[cyan, line width=0.3mm] (150:3cm) -- (210:3cm);
%\draw[red, line width=0.3mm] (210:3cm) -- node[midway, above right, xshift=-1mm] {$\textcolor{black}{\gamma_2}$} (360:3cm);
%\draw[black, line width=0.3mm] (212:3cm) -- (358:3cm);
                  
                    %\node at (210:3cm) [left] {$\Bar{b}$};
                    \node at (360:3cm) [right] {$\ast$};

\node at (150:3cm) [left] {$a$};
                      
                    \end{scope}
                    \begin{scope}[yshift=-7.3cm]
                         \draw (90:3cm) -- (120:3cm) -- (150:3cm) -- (180:3cm) -- (210:3cm) -- (240:3cm) -- (270:3cm) -- (300:3cm) -- (330:3cm) -- (360:3cm) -- (30:3cm) -- (60:3cm) --  cycle;
                  
                    \draw[-{Latex[length=2mm]}] (270:3cm) -- node[midway, above left,xshift=1mm] {} (90:3cm);
                  
                     \draw[red, line width=0.3mm] (150:3cm) -- (270:3cm); 
                   
 \draw[red, line width=0.3mm] (90:3cm) -- (330:3cm);

                    %\node at (120:3cm) [left] {$a$};
                    \node at (90:3cm) [above] {$\Bar{b}$};
                    \node at (330:3cm) [right] {$\Bar{a}$};
                    %\node at (60:3cm) [right] {$\hat{a}$};

\node at (150:3cm) [left] {$a$};
                    %\node at (240:3cm) [left] {$\hat{d}$};
                    %\node at (300:3cm) [right] {$d$};
                    \node at (270:3cm) [below] {$b$};

                    \begin{scope}[xshift=8cm]
                      \draw (90:3cm) -- (120:3cm) -- (150:3cm) -- (180:3cm) -- (210:3cm) -- (240:3cm) -- (270:3cm) -- (360:3cm)  --  cycle;
                    %\draw (90:3cm) -- node[midway, above, xshift=-1mm, yshift=-1mm] {} (180:3cm);
                    %\draw (90:3cm) -- node[midway, above left,xshift=1.2mm] {} (225:3cm);
                    \draw (90:3cm) -- node[midway, above left,xshift=1mm] {} (270:3cm);
                   % \draw (90:3cm) -- node[midway, above right, xshift=-1mm] {} (315:3cm);
                    %\draw (90:3cm) -- node[midway, above right, xshift=-1mm] {} (360:3cm);
                    \draw[red, line width=0.3mm] (150:3cm) -- node[midway, above right, xshift=-2mm, yshift=-1mm] {$\textcolor{black}{\gamma}$} (270:3cm); 
                    %\draw[black, line width=0.3mm] (152:3cm) --  (358:3cm); 
                    % \draw[cyan, line width=0.3mm] (150:3cm) -- (210:3cm);
%\draw[red, line width=0.3mm] (210:3cm) -- node[midway, above right, xshift=-1mm] {$\textcolor{black}{\gamma_2}$} (360:3cm);
%\draw[black, line width=0.3mm] (212:3cm) -- (358:3cm);
                  
                    %\node at (210:3cm) [left] {$\Bar{b}$};
                    \node at (360:3cm) [right] {$\ast$};
 \node at (270:3cm) [below] {$b$};
\node at (150:3cm) [left] {$a$};
                      
                    \end{scope}
                    \end{scope}
                \end{tikzpicture}
                               \caption{On the left, two $\theta$-orbits $[a,\Bar{a}]$ and $[a,b]$. On the right, their restrictions.}
                               \label{type B 1}
                           \end{figure} 
 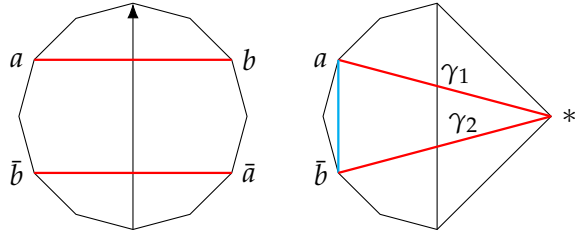
\begin{figure}[H]
                               \centering
                \begin{tikzpicture}[scale=0.5]
                    \draw (90:3cm) -- (120:3cm) -- (150:3cm) -- (180:3cm) -- (210:3cm) -- (240:3cm) -- (270:3cm) -- (300:3cm) -- (330:3cm) -- (360:3cm) -- (30:3cm) -- (60:3cm) --  cycle;
                    %\draw (90:3cm) -- node[midway, above, xshift=-1mm, yshift=-1mm] {} (180:3cm);
                    %\draw (90:3cm) -- node[midway, above left,xshift=1.2mm] {} (225:3cm);
                    \draw[-{Latex[length=2mm]}] (270:3cm) -- node[midway, above left,xshift=1mm] {} (90:3cm);
                   % \draw (90:3cm) -- node[midway, above right, xshift=-1mm] {} (315:3cm);
                    %\draw (90:3cm) -- node[midway, above right, xshift=-1mm] {} (360:3cm);
                     
                    %\draw[cyan, line width=0.3mm] (330:3cm) -- (30:3cm);
                     \draw[red, line width=0.3mm] (150:3cm) -- (30:3cm); 
                    % \draw[cyan, line width=0.3mm] (150:3cm) -- (210:3cm);
                    \draw[red, line width=0.3mm] (210:3cm) -- (330:3cm);
                    %\draw[black, line width=0.3mm] (330:3cm) -- (150:3cm);
                    %\draw[black, line width=0.3mm] (210:3cm) -- (30:3cm);

                    %\node at (120:3cm) [left] {$a$};
                    \node at (210:3cm) [left] {$\Bar{b}$};
                    \node at (330:3cm) [right] {$\Bar{a}$};
                    %\node at (60:3cm) [right] {$\hat{a}$};

\node at (150:3cm) [left] {$a$};
                    %\node at (240:3cm) [left] {$\hat{d}$};
                    %\node at (300:3cm) [right] {$d$};
                    \node at (30:3cm) [right] {$b$};

                    \begin{scope}[xshift=8cm]
                      \draw (90:3cm) -- (120:3cm) -- (150:3cm) -- (180:3cm) -- (210:3cm) -- (240:3cm) -- (270:3cm) -- (360:3cm)  --  cycle;
                    %\draw (90:3cm) -- node[midway, above, xshift=-1mm, yshift=-1mm] {} (180:3cm);
                    %\draw (90:3cm) -- node[midway, above left,xshift=1.2mm] {} (225:3cm);
                    \draw (90:3cm) -- node[midway, above left,xshift=1mm] {} (270:3cm);
                   % \draw (90:3cm) -- node[midway, above right, xshift=-1mm] {} (315:3cm);
                    %\draw (90:3cm) -- node[midway, above right, xshift=-1mm] {} (360:3cm);
                    \draw[red, line width=0.3mm] (150:3cm) -- node[midway, above right, xshift=-2mm, yshift=-1mm] {$\textcolor{black}{\gamma_1}$} (360:3cm); 
                    %\draw[black, line width=0.3mm] (152:3cm) --  (358:3cm); 
                     \draw[cyan, line width=0.3mm] (150:3cm) -- (210:3cm);
\draw[red, line width=0.3mm] (210:3cm) -- node[midway, above right, xshift=-1mm] {$\textcolor{black}{\gamma_2}$} (360:3cm);
%\draw[black, line width=0.3mm] (212:3cm) -- (358:3cm);
                  
                    \node at (210:3cm) [left] {$\Bar{b}$};
                    \node at (360:3cm) [right] {$\ast$};

\node at (150:3cm) [left] {$a$};
                      
                    \end{scope}
                \end{tikzpicture}
                               \caption{On the left, a $\theta$-orbit $[a,b]$. On the right, its restriction in red and the diagonal $(a,\Bar{b})$ in blue.}
                               \label{type B}
                           \end{figure}
\end{definition}

\begin{theorem}\label{theorem1}
    Let $T$ be a $\theta$-invariant triangulation of $\Poly_{2n+2}$ with oriented diameter $d$, and let $\mathcal{A}=\mathcal{A}^B (T)$ be the cluster algebra of type $B_n$ with principal coefficients in $T$. 
    Let $[a,b]$ be an orbit of the action of $\theta$ on the diagonals of the polygon, and $x_{ab}$ the cluster variable of $\mathcal{A}$ which corresponds to $[a,b]$. Let $F_{ab}$ and $\bold{g}_{ab}$ denote the $F$-polynomial and the $\bold{g}$-vector of $x_{ab}$, respectively. 
    Then $F_{ab}=F_{ab}^B$ and $\bold{g}_{ab}=\bold{g}_{ab}^B$.
\end{theorem}

\begin{remark}
    Since, for a diagonal $\gamma$ of $\Poly_{n+3}$, $F_\gamma$ and $\bold{g}_\gamma$ have an explicit description, for example in terms of perfect matchings of the snake graph associated with $\gamma$ \cite{MS,CSI}, Theorem \ref{theorem1} also allows us to get the expansion of cluster variables of type $B_n$ in terms of the cluster variables of the initial seed. 
\end{remark}

\begin{example}\label{example_F_poly}
By Theorem \ref{theorem1}, the $F$-polynomial of the cluster variable of type $B_3$ which corresponds to the $\theta$-orbit $[a,b]$ of $\Poly_8$ in Figure \ref{ex_f_poly} is   
\begin{align*}
 F_{ab}&=F_{\gamma_1}F_{\gamma_2}-y_1y_2y_3= (y_3+1)(y_1y_2y_3+y_1y_3+y_1+y_3+1)-y_1y_2y_3\\
&=y_1y_2y_3^2+y_1y_3^2+2y_1y_3+y_3^2+y_1+2y_3+1,    
\end{align*}

and the $\bold{g}$-vector is 
\begin{center}
$\displaystyle\bold{g}_{ab}=D(\bold{g}_{\gamma_1}+\bold{g}_{\gamma_2}+\bold{e}_3)=D(\displaystyle\begin{pmatrix}
    0 \\ 1 \\ -1
\end{pmatrix}+\displaystyle\begin{pmatrix}
    -1 \\ 1 \\ -1
\end{pmatrix}+\displaystyle\begin{pmatrix}
    0 \\ 0 \\ 1
\end{pmatrix})=D(\displaystyle\begin{pmatrix}
    -1 \\ 2 \\ -1
\end{pmatrix})=\displaystyle\begin{pmatrix}
    -1 \\ 2 \\ -2
\end{pmatrix}$.   
\end{center}
      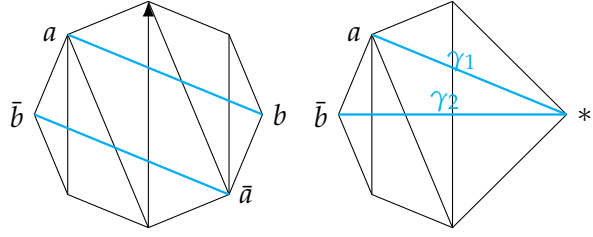
\begin{figure}[H]
                               \centering
                \begin{tikzpicture}[scale=0.5]
                    \draw (90:3cm) -- (135:3cm) -- (180:3cm) -- (225:3cm) -- (270:3cm) -- (315:3cm) -- (360:3cm) -- (45:3cm) -- cycle;
                    \draw (135:3cm) -- (225:3cm);;
                    \draw (45:3cm) -- (315:3cm);
                    \draw[-{Latex[length=2mm]}] (270:3cm) -- node[midway, above left,xshift=1mm] {} (90:3cm);
                    \draw (135:3cm) -- node[midway, above left,xshift=1mm] {} (270:3cm);
                    \draw (315:3cm) -- node[midway, above left,xshift=1mm] {} (90:3cm);
                    %\draw (90:3cm) -- node[midway, above right, xshift=-1mm] {} (315:3cm);
                    %\draw (90:3cm) -- node[midway, above right, xshift=-1mm] {} (360:3cm);
                    \draw[cyan, line width=0.3mm] (135:3cm) -- (360:3cm); 
                    \draw[cyan, line width=0.3mm] (180:3cm) -- (315:3cm);

                    \node at (135:3cm) [left] {$a$};
                    \node at (180:3cm) [left] {$\Bar{b}$};
                    \node at (360:3cm) [right] {$b$};
                    \node at (315:3cm) [right] {$\Bar{a}$};

                    \begin{scope}[xshift=8cm]
                      \draw (90:3cm) -- (135:3cm) -- (180:3cm) -- (225:3cm) -- (270:3cm) -- (360:3cm) -- cycle;
 
                    \draw (135:3cm) -- (225:3cm);;
                    %\draw (45:3cm) -- (315:3cm);
                    \draw (90:3cm) -- node[midway, above left,xshift=1mm] {} (270:3cm);
                    \draw (135:3cm) -- node[midway, above left,xshift=1mm] {} (270:3cm);
                    %\draw (45:3cm) -- node[midway, above left,xshift=1mm] {} (270:3cm);
                    %\draw (90:3cm) -- node[midway, above right, xshift=-1mm] {} (315:3cm);
                    %\draw (90:3cm) -- node[midway, above right, xshift=-1mm] {} (360:3cm);

\draw[cyan, line width=0.3mm] (180:3cm) -- node[midway, above, xshift=-1mm, yshift=-1mm] {$\gamma_2$} (360:3cm); 
\draw[cyan, line width=0.3mm] (135:3cm) -- node[midway, above, xshift=-1mm, yshift=-1mm] {$\gamma_1$} (360:3cm);

\node at (135:3cm) [left] {$a$};
\node at (180:3cm) [left] {$\bar{b}$};
\node at (360:3cm) [right] {$\ast$}; 
                    \end{scope}
                \end{tikzpicture}
                               \caption{A $\theta$-orbit $[a,b]$ in a triangulated octagon and its restriction.}
                               \label{ex_f_poly}
                           \end{figure}
               
\end{example}
In order to present the proof of Theorem \ref{theorem1}, we first need some lemmas.
\begin{lemma}\label{lemma1}
If each diagonal of $[a,b]$ crosses only one diagonal of $T$, then $F_{ab}=F_{ab}^B$ and $\bold{g}_{ab}=\bold{g}_{ab}^B$. 
\end{lemma}

\begin{proof}
    Let $T=\{\tau_1,\dots, \tau_{2n-1}\}$. Since each diagonal of $[a,b]$ crosses only one diagonal of $T$, either $[a,b]$ is a pair of diagonals which do not cross $d$ or $[a,b]=\{a,\bar{a}\}$ is the diagonal which crosses only $d$. Therefore, $\text{Res}([a,b])=\{\gamma_j\}$, where $\gamma_j$ is the diagonal of $\Poly_{n+3}$ which crosses only $\tau_j$. Let $DB(\Bar{T})=(b_{ij})$ and $B(\Bar{T})=(\Bar{b}_{ij})$. 
    We have 
    \begin{equation}
        x_{ab}u_j=y_j\displaystyle\prod_{b_{ij}>0}u_i^{b_{ij}}+\displaystyle\prod_{b_{ij}<0}u_i^{-b_{ij}},
    \end{equation}
    and
     \begin{equation}
        x_{\gamma_j}u_j=y_j\displaystyle\prod_{\Bar{b}_{ij}>0}u_i^{\Bar{b}_{ij}}+\displaystyle\prod_{\Bar{b}_{ij}<0}u_i^{-\Bar{b}_{ij}}.
    \end{equation}
    So
    \begin{equation}
        F_{ab}=y_j+1=F_{\gamma_j}=F_{ab}^B.
    \end{equation}
    If $k\neq n$,
    \begin{equation}
        (\bold{g}_{ab})_k=\displaystyle\bigg(\textbf{deg}\displaystyle\bigg(\displaystyle\frac{\displaystyle\prod_{b_{ij}<0}u_i^{-b_{ij}}}{u_j}\bigg)\bigg)_k=\displaystyle\bigg(\textbf{deg}\displaystyle\bigg(\displaystyle\frac{\displaystyle\prod_{\Bar{b}_{ij}<0}u_i^{-\Bar{b}_{ij}}}{u_j}\bigg)\bigg)_k=(\bold{g}_{\gamma_j})_k=
        (\bold{g}_{ab}^B)_k.
    \end{equation}
If $k=n$ and $j \neq n$,
    \begin{equation}
        (\bold{g}_{ab})_n=\displaystyle\bigg(\textbf{deg}\displaystyle\bigg(\displaystyle\frac{\displaystyle\prod_{b_{ij}<0}u_i^{-b_{ij}}}{u_j}\bigg)\bigg)_n=2\displaystyle\bigg(\textbf{deg}\displaystyle\bigg(\displaystyle\frac{\displaystyle\prod_{\Bar{b}_{ij}<0}u_i^{-\Bar{b}_{ij}}}{u_j}\bigg)\bigg)_n=2(\bold{g}_{\gamma_j})_n=
        (\bold{g}_{ab}^B)_n.
    \end{equation}
Finally, if $k=n$ and $j=n$,
\begin{equation}
    (\bold{g}_{ab})_n=\displaystyle\bigg(\textbf{deg}\displaystyle\bigg(\displaystyle\frac{1}{u_n}\bigg)\bigg)_n=-1=(\bold{g}_{\gamma_n})_n=2(\bold{g}_{\gamma_n})_n+1=(\bold{g}_{ab}^B)_n.
\end{equation}
    
\end{proof}
\begin{lemma}\label{lemma_c_vect_typeB}
Let $B$ be a skew-symmetric $n\times n$ matrix, and let $I$ be the $n\times n$ identity matrix. Let
$D=\text{diag}(1,\dots,1,2)$ be $n\times n$ diagonal matrix with diagonal entries $(1,\dots,1,2)$. Let 
$\mu_{i_1} \cdots \mu_{i_k}(\left[\begin{matrix}
   $B$\\ I
\end{matrix} \right])=\left[\begin{matrix}
    B' \\ C
\end{matrix} \right]$, and let $\mu_{i_1} \cdots \mu_{i_k}(\left[\begin{matrix}
    DB \\ I
\end{matrix} \right])=\left[\begin{matrix}
    DB' \\ C'
\end{matrix} \right]$, for any 
$1\leq i_1 < \cdots < i_k \leq n$. Then, 
\begin{itemize}
    \item [i)] if $i_j \neq n$ for every $j=1, \dots, k$, $C=C'$;
    \item [ii)] if $i_k=n$, the columns $(C')^1,\dots,(C')^{i_1-1}$ of $C'$ are equal to $DC^1,\dots, DC^{i_1-1}$.
\end{itemize}
\end{lemma}

\begin{proof}
    $i)$ holds since $\mu_{i_j}$ does not consider the $n$-th row for every $j=1,\dots,k$.

    We prove $ii)$ by induction on $k$. If $k=1$, $(C')^n=-\bold{e}_n$, and for $j\neq n$
    \begin{center}
        $(C')^j=\begin{cases}
        \bold{e}_j \hspace{1.2cm}\text{if $b_{nj}\leq 0$}\\
        \bold{e}_j+2\bold{e}_n \hspace{0.2cm}\text{otherwise}
    \end{cases}=DC^j$.
         \end{center}
Assume $k>1$. By inductive hypothesis, the columns $(C')^1,\dots,(C')^{i_1-1}$ of $C'$ are equal to $DC^1,\dots,$ $ DC^{i_1-1}$. Then, we mutate at $i_1-1$. If $\mu_{i_1-1}\mu_{i_1}\cdots \mu_{i_k}(\left[\begin{matrix}
   $B$\\ I
\end{matrix} \right])=\left[\begin{matrix}
    B'' \\ C''
\end{matrix} \right]$, and $\mu_{i_1-1}\mu_{i_1}\cdots \mu_{i_k}(\left[\begin{matrix}
    DB \\ I
\end{matrix} \right])=\left[\begin{matrix}
    DB'' \\ C'''
\end{matrix} \right]$, we have that $(C''')^j=D(C'')^j$ for every $j=1,\dots,i_1-2$.
\end{proof}

\begin{lemma}[\cite{S_CAUS}, Lemma 4.3]\label{lemma_schiffler}
    Let $\Bar{T}=\{ \tau_1, \dots, \tau_n \}$ be a triangulation of $\Poly_{n+3}$. Let $\gamma \not\subset \Bar{T}$ be a diagonal on which we fixed an orientation such that $\gamma$ is going from $s$ to $t$. Let $s = p_0,p_1,...,p_d,p_{d+1} = t$ be the intersection points of $\gamma$ and $\Bar{T}$ in order of occurrence on $\gamma$, and let $i_1,i_2,\dots,i_d$ be such that $p_k$ lies on $\tau_{i_k}$, for $k=1,\dots,d$. Then $\gamma \in \mu_{i_1}\dots \mu_{i_d}(\Bar{T})$, i.e. $x_\gamma \in \mu_{i_1}\dots \mu_{i_d}(\bold{u}_{\Bar{T}})$. 
\end{lemma}

\begin{proof}[Proof of Theorem~\ref{theorem1}]
    We prove the theorem by induction on the number $k$ of intersections between each diagonal of $[a,b]$ and 
    $T=\{\tau_1, \dots, \tau_n=d,\dots,\tau_{2n-1}\}$.

    If $k=0$, the theorem holds by Definition \ref{def_type_B}. If $k=1$, the theorem holds by Lemma \ref{lemma1}.
    Assume $k>1$. Let $\Bar{T}=\text{Res}(T)=\{ \tau_1, \dots, \tau_n=d \}$, and let $\bold{u}_{\Bar{T}}=\{ u_{\tau_1}, \dots, u_{\tau_n} \}=\{u_1, \dots, u_n\}$. There are three cases to consider.
    \begin{itemize}
\item [1)] Let $[a,b]=\{(a,b),(\Bar{b},\Bar{a})\}$ be such that $\text{Res}([a,b])=\{(a,b)\}$. Let $a=p_0,p_1,\dots,p_k,$
$p_{k+1} = b$ be the intersection points of $(a,b)$ and $\Bar{T}$ in order of occurrence on $(a,b)$, and let $i_1,i_2,\dots,i_k$ be such that $p_j$ lies on the diagonal $\tau_{i_j}\in \Bar{T}$, for $j=1,\dots,k$. Let $[c,d]=\{\tau_{i_1},\tau_{i_{2n-i_1}}\}$.
  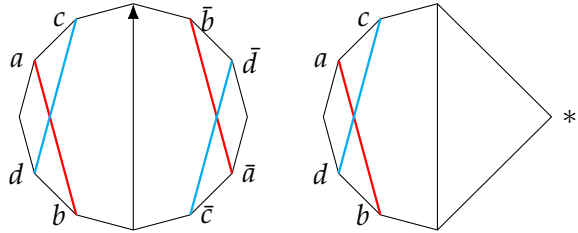
\begin{figure}[H]
                               \centering
                \begin{tikzpicture}[scale=0.5]
                    \draw (90:3cm) -- (120:3cm) -- (150:3cm) -- (180:3cm) -- (210:3cm) -- (240:3cm) -- (270:3cm) -- (300:3cm) -- (330:3cm) -- (360:3cm) -- (30:3cm) -- (60:3cm) --  cycle;
                    %\draw (90:3cm) -- node[midway, above, xshift=-1mm, yshift=-1mm] {} (180:3cm);
                    %\draw (90:3cm) -- node[midway, above left,xshift=1.2mm] {} (225:3cm);
                   \draw[-{Latex[length=2mm]}] (270:3cm) -- node[midway, above left,xshift=1mm] {} (90:3cm);
                   % \draw (90:3cm) -- node[midway, above right, xshift=-1mm] {} (315:3cm);
                    %\draw (90:3cm) -- node[midway, above right, xshift=-1mm] {} (360:3cm);
                     
                    \draw[red, line width=0.3mm] (330:3cm) -- (60:3cm);
                     \draw[red, line width=0.3mm] (150:3cm) -- (240:3cm); 
                     \draw[cyan, line width=0.3mm] (120:3cm) -- (210:3cm);
                    \draw[cyan, line width=0.3mm] (300:3cm) -- (30:3cm);

                    \node at (120:3cm) [left] {$c$};
                    \node at (210:3cm) [left] {$d$};
                    \node at (330:3cm) [right] {$\Bar{a}$};
                    \node at (60:3cm) [right] {$\Bar{b}$};

\node at (150:3cm) [left] {$a$};
                    \node at (240:3cm) [left] {$b$};
                    \node at (300:3cm) [right] {$\Bar{c}$};
                    \node at (30:3cm) [right] {$\Bar{d}$};

                    \begin{scope}[xshift=8cm]
                      \draw (90:3cm) -- (120:3cm) -- (150:3cm) -- (180:3cm) -- (210:3cm) -- (240:3cm) -- (270:3cm) -- (360:3cm)  --  cycle;
                    %\draw (90:3cm) -- node[midway, above, xshift=-1mm, yshift=-1mm] {} (180:3cm);
                    %\draw (90:3cm) -- node[midway, above left,xshift=1.2mm] {} (225:3cm);
                    \draw (90:3cm) -- node[midway, above left,xshift=1mm] {} (270:3cm);
                   % \draw (90:3cm) -- node[midway, above right, xshift=-1mm] {} (315:3cm);
                    %\draw (90:3cm) -- node[midway, above right, xshift=-1mm] {} (360:3cm);
                    \draw[red, line width=0.3mm] (150:3cm) -- (240:3cm); 
                     \draw[cyan, line width=0.3mm] (120:3cm) -- (210:3cm);

                    \node at (120:3cm) [left] {$c$};
                    \node at (210:3cm) [left] {$d$};
                    \node at (360:3cm) [right] {$\ast$};

\node at (150:3cm) [left] {$a$};
                    \node at (240:3cm) [left] {$b$};
                      
                    \end{scope}
                \end{tikzpicture}
                               \caption{On the left, the two $\theta$-orbits $[a,b]$ and $[c,d]$. On the right, their restrictions.}
                               \label{fig:enter-label}
                           \end{figure}
Then, by Lemma \ref{lemma_schiffler}, $(a,b) \in \mu_{i_1}\cdots \mu_{i_k}(\Bar{T})$. Therefore, the $\bold{c}$-vector corresponding to the exchange between $[a,b]$ and $[c,d]$ is the bottom part of the $i_1$-th column of $\mu_{i_2}\cdots \mu_{i_k}(\left[\begin{matrix}
    DB(\Bar{T}) \\ I
\end{matrix} \right])$. Since $i_j \neq n$ for each $j$, by Lemma \ref{lemma_c_vect_typeB} i), this is equal to the bottom part of the $i_1$-th column of $\mu_{i_2}\cdots \mu_{i_k}(\left[\begin{matrix}
    B(\Bar{T}) \\ I
\end{matrix} \right])$, which is given by Proposition \ref{up:skein1}. Therefore, we have the following exchange relation
\begin{equation}\label{thm1_eqn1}
u_{i_1}x_{ab}=\bold{y}^{\bold{d}_{ac,bd}}x_{ad}x_{bc}+\bold{y}^{\bold{d}_{ad,bc}}x_{ac}x_{bd}.
\end{equation}
Since $(c,d)$ is the first diagonal of $T$ that is crossed by $(a,b)$, $(a,c)$ and $(a,d)$ must be either boundary edges or diagonals of $\Bar{T}$. It follows from \ref{thm1_eqn1} that
\begin{equation}
    F_{ab}=\bold{y}^{\bold{d}_{ac,bd}}F_{bc}+\bold{y}^{\bold{d}_{ad,bc}}F_{bd}.
\end{equation}
Since each diagonal of $T$ which crosses $(b,c)$ (resp. $(b,d)$) also crosses $(a,b)$, the number of intersections between $(b,c)$ (resp. $(b,d)$) and $T$ is strictly lower than the number of crossings between $(a,b)$ and $T$. By inductive hypothesis and Proposition \ref{up:skein1},
\begin{equation}
    F_{ab}=\bold{y}^{\bold{d}_{ac,bd}}F_{(b,c)}+\bold{y}^{\bold{d}_{ad,bc}}F_{(b,d)}=F_{(a,b)}=F_{ab}^B.
\end{equation}
\begin{comment}
   It also follows from \ref{thm1_eqn1} that
\begin{equation}
    \bold{e}_{i_1}+\bold{g}_{ab}=\begin{cases}
        \bold{g}_{ad}+\bold{g}_{bc}\hspace{1cm} \text{if $\bold{y}^{\bold{d}_{ac,bd}}=1$}\\
        \bold{g}_{ac}+\bold{g}_{bd}\hspace{1cm} \text{otherwise}.
    \end{cases}
\end{equation}
By inductive hypothesis and Proposition \ref{up:skein1},

\begin{equation*}
    \bold{g}_{ab}=\begin{cases}
        -\bold{e}_{i_1}+ D\bold{g}_{(a,d)}+D\bold{g}_{(b,c)}=D(-\bold{e}_{i_1}+ \bold{g}_{(a,d)}+\bold{g}_{(b,c)}) \hspace{0.2cm}\text{if $\bold{y}^{\bold{d}_{ac,bd}}=1$}\\
        -\bold{e}_{i_1} + D\bold{g}_{(a,c)}+D\bold{g}_{(b,d)}=D(-\bold{e}_{i_1} + \bold{g}_{(a,c)}+\bold{g}_{(b,d)})\hspace{0.2cm} \text{otherwise}.
        \end{cases}\hspace{-0.3cm}=D\bold{g}_{(a,b)}=\bold{g}_{ab}^B.
\end{equation*} 
\end{comment}

\item [2)] Let $[a,b]=[a,\Bar{a}]$ be a diameter. So $\text{Res}([a,\Bar{a}])=\{(a,\ast)\}$. Let $\ast=p_0,p_1,\dots,p_s,p_{s+1} = a$ be the intersection points of $(a,\ast)$ and $\Bar{T}$ in order of occurrence on $(\ast,a)$, $s\leq k$, and let $i_1,i_2,\dots,i_s$ be such that $p_j$ lies on the diagonal $\tau_{i_j}\in \Bar{T}$, for $j=1,\dots,s$. Thus $i_1 = n$. Let $[b,\Bar{b}]=\{\tau_n\}=\{d\}$.
 \begin{figure}[H]
                               \centering
                \begin{tikzpicture}[scale=0.5]
                    \draw (90:3cm) -- (120:3cm) -- (150:3cm) -- (180:3cm) -- (210:3cm) -- (240:3cm) -- (270:3cm) -- (300:3cm) -- (330:3cm) -- (360:3cm) -- (30:3cm) -- (60:3cm) --  cycle;
                    %\draw (90:3cm) -- node[midway, above, xshift=-1mm, yshift=-1mm] {} (180:3cm);
                    %\draw (90:3cm) -- node[midway, above left,xshift=1.2mm] {} (225:3cm);
                    \draw (90:3cm) -- node[midway, above left,xshift=1mm] {} (270:3cm);
                   % \draw (90:3cm) -- node[midway, above right, xshift=-1mm] {} (315:3cm);
                    %\draw (90:3cm) -- node[midway, above right, xshift=-1mm] {} (360:3cm);
                     
                    \draw[-{Latex[length=2mm]},cyan, line width=0.3mm] (270:3cm) -- (90:3cm);
                     \draw[red, line width=0.3mm] (150:3cm) -- (330:3cm);

                    \node at (90:3cm) [above] {$b$};
                    \node at (270:3cm) [below] {$\Bar{b}$};

\node at (150:3cm) [left] {$a$};
                    \node at (330:3cm) [right] {$\Bar{a}$};

                    \begin{scope}[xshift=8cm]
                      \draw (90:3cm) -- (120:3cm) -- (150:3cm) -- (180:3cm) -- (210:3cm) -- (240:3cm) -- (270:3cm) -- (360:3cm)  --  cycle;
                    %\draw (90:3cm) -- node[midway, above, xshift=-1mm, yshift=-1mm] {} (180:3cm);
                    %\draw (90:3cm) -- node[midway, above left,xshift=1.2mm] {} (225:3cm);
                    \draw (90:3cm) -- node[midway, above left,xshift=1mm] {} (270:3cm);
                   % \draw (90:3cm) -- node[midway, above right, xshift=-1mm] {} (315:3cm);
                    %\draw (90:3cm) -- node[midway, above right, xshift=-1mm] {} (360:3cm);
                   
                    \draw[cyan, line width=0.3mm] (90:3cm) -- (270:3cm);
                     \draw[red, line width=0.3mm] (150:3cm) -- (360:3cm);

                    \node at (90:3cm) [above] {$b$};
                    \node at (270:3cm) [below] {$\Bar{b}$};

\node at (150:3cm) [left] {$a$};
                    \node at (360:3cm) [right] {$\ast$};
                      
                    \end{scope}
                \end{tikzpicture}
                               \caption{On the left, the two $\rho$-orbits $[a,\Bar{a}]$ and $[b,\Bar{b}]$. On the right, their restrictions.}
                               \label{fig:enter-label}
                           \end{figure}
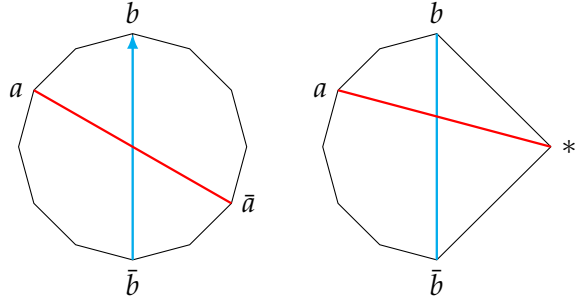    
 Then, by Lemma \ref{lemma_schiffler}, $(a,\ast) \in \mu_{i_1}\cdots \mu_{i_s}(\Bar{T})$. Therefore, the $\bold{c}$-vector corresponding to the exchange between $[a,\Bar{a}]$ and $[b,\Bar{b}]$ is the bottom part of the $i_1$-th column of $\mu_{i_2}\cdots \mu_{i_s}(\left[\begin{matrix}
    DB(\Bar{T}) \\ I
\end{matrix} \right])$. Since $i_j \neq n$ for each $j$, by Lemma \ref{lemma_c_vect_typeB} i), this is equal to the bottom part of the $i_1$-th column of $\mu_{i_2}\cdots \mu_{i_s}(\left[\begin{matrix}
    B(\Bar{T}) \\ I
\end{matrix} \right])$, which is given by Proposition \ref{up:skein1}. Therefore, we have the following exchange relation
\begin{equation}\label{thm1_eqn1_}
u_{n}x_{a\Bar{a}}=\bold{y}^{\bold{d}_{ab,\Bar{b}\ast}}x_{a\Bar{b}}+\bold{y}^{\bold{d}_{b\ast,a\Bar{b}}}x_{ab}.
\end{equation}
It follows from \ref{thm1_eqn1_} that
\begin{equation}
    F_{a\Bar{a}}=\bold{y}^{\bold{d}_{ab,\Bar{b}\ast}}F_{a\Bar{b}}+\bold{y}^{\bold{d}_{b\ast,a\Bar{b}}}F_{ab}.
\end{equation}
By inductive hypothesis and Proposition \ref{up:skein1},
\begin{equation}
    F_{a\Bar{a}}=\bold{y}^{\bold{d}_{ab,\Bar{b}\ast}}F_{(a,\Bar{b})}+\bold{y}^{\bold{d}_{b\ast,a\Bar{b}}}F_{(a,b)}=F_{(a,\ast)}=F_{a\Bar{a}}^B.
\end{equation}
\begin{comment}
    It also follows from \ref{thm1_eqn1} that
\begin{equation}
    \bold{e}_{n}+\bold{g}_{a\Bar{a}}=\begin{cases}
        \bold{g}_{a\Bar{b}}\hspace{1cm} \text{if $\bold{y}^{\bold{d}_{ab,\Bar{b}\ast}}=1$}\\
        \bold{g}_{ab}\hspace{1cm} \text{otherwise}.
    \end{cases}
\end{equation}
By inductive hypothesis and Proposition \ref{up:skein1},

\begin{equation*}
    \bold{g}_{a\Bar{a}}=\begin{cases}
        -\bold{e}_{n}+ D\bold{g}_{(a,\Bar{b})}=D(-\bold{e}_n+\bold{g}_{(a,\Bar{b})})+\bold{e}_n \hspace{0.8cm}\text{if $\bold{y}^{\bold{d}_{ab,\Bar{b}\ast}}=1$}\\
        -\bold{e}_{n} + D\bold{g}_{(a,b)}=D(-\bold{e}_n+\bold{g}_{(a,b)})+\bold{e}_n\hspace{0.8cm} \text{otherwise}.
        \end{cases}=D\bold{g}_{(a,\ast)}+\bold{e}_n=\bold{g}_{a\Bar{a}}^B.
\end{equation*} 

\end{comment}

        \item [3)] Let $[a,b]=\{(a,b),(\Bar{a},\Bar{b})\}$ be  a pair of diagonals which cross $d$, so $\text{Res}([a,b])=\{(a,\ast),(\Bar{b},\ast)\}$.
Let $a=p_0,p_1,\dots,p_s,p_{s+1} = \ast$ be the intersection points of $(a,\ast)$ and $\Bar{T}$ in order of occurrence on $(a,\ast)$, $s\leq k$, and let $i_1,i_2,\dots,i_s$ be such that $p_j$ lies on the diagonal $\tau_{i_j}\in \Bar{T}$, for $j=1,\dots,s$. So $i_s=n$. Let $[c,d]=\{\tau_{i_1},\tau_{i_{2n-i_1}}\}$. Assume that $(c,d)=\tau_{i_1}$ intersects $(a,\ast)$ (otherwise we consider $(\Bar{b},\ast)$ instead of $(a,\ast)$).

Then, by Lemma \ref{lemma_schiffler}, $(a,\ast) \in \mu_{i_1}\cdots \mu_{i_s}(\Bar{T})$. Therefore, the $\bold{c}$-vector corresponding to the exchange between $[a,b]$ and $[c,d]$ is the bottom part of the $i_1$-th column of $\mu_{i_2}\cdots \mu_{i_s}(\left[\begin{matrix}
    DB(\Bar{T}) \\ I
\end{matrix} \right])$. Since $i_s=n$, by Lemma \ref{lemma_c_vect_typeB} ii), this is equal to $DC^{i_1}$, where $C^{i_1}$ is the bottom part of the $i_1$-th column of $\mu_{i_2}\cdots \mu_{i_s}(\left[\begin{matrix}
    B(\Bar{T}) \\ I
\end{matrix} \right])$, which is given by Proposition \ref{up:skein1}. 

Now, we have two cases to consider:
\begin{itemize}
    \item [a)] $c$ is not an endpoint of $\tau_n$;
     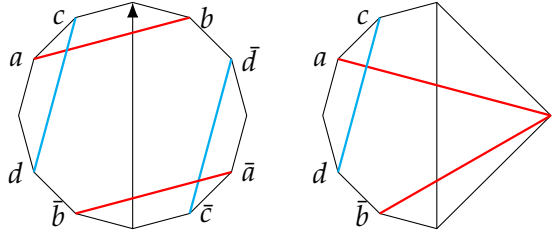
\begin{figure}[H]
                               \centering
                \begin{tikzpicture}[scale=0.5]
                    \draw (90:3cm) -- (120:3cm) -- (150:3cm) -- (180:3cm) -- (210:3cm) -- (240:3cm) -- (270:3cm) -- (300:3cm) -- (330:3cm) -- (360:3cm) -- (30:3cm) -- (60:3cm) --  cycle;
                    %\draw (90:3cm) -- node[midway, above, xshift=-1mm, yshift=-1mm] {} (180:3cm);
                    %\draw (90:3cm) -- node[midway, above left,xshift=1.2mm] {} (225:3cm);
                    \draw[-{Latex[length=2mm]}] (270:3cm) -- node[midway, above left,xshift=1mm] {} (90:3cm);
                   % \draw (90:3cm) -- node[midway, above right, xshift=-1mm] {} (315:3cm);
                    %\draw (90:3cm) -- node[midway, above right, xshift=-1mm] {} (360:3cm);
                     
                    \draw[cyan, line width=0.3mm] (300:3cm) -- (30:3cm);
                     \draw[red, line width=0.3mm] (150:3cm) -- (60:3cm); 
                     \draw[cyan, line width=0.3mm] (120:3cm) -- (210:3cm);
                    \draw[red, line width=0.3mm] (240:3cm) -- (330:3cm);

                    \node at (120:3cm) [left] {$c$};
                    \node at (210:3cm) [left] {$d$};
                    \node at (330:3cm) [right] {$\Bar{a}$};
                    \node at (60:3cm) [right] {$b$};

\node at (150:3cm) [left] {$a$};
                    \node at (240:3cm) [left] {$\Bar{b}$};
                    \node at (300:3cm) [right] {$\Bar{c}$};
                    \node at (30:3cm) [right] {$\Bar{d}$};

                    \begin{scope}[xshift=8cm]
                      \draw (90:3cm) -- (120:3cm) -- (150:3cm) -- (180:3cm) -- (210:3cm) -- (240:3cm) -- (270:3cm) -- (360:3cm)  --  cycle;
                    %\draw (90:3cm) -- node[midway, above, xshift=-1mm, yshift=-1mm] {} (180:3cm);
                    %\draw (90:3cm) -- node[midway, above left,xshift=1.2mm] {} (225:3cm);
                    \draw (90:3cm) -- node[midway, above left,xshift=1mm] {} (270:3cm);
                   % \draw (90:3cm) -- node[midway, above right, xshift=-1mm] {} (315:3cm);
                    %\draw (90:3cm) -- node[midway, above right, xshift=-1mm] {} (360:3cm);
                    \draw[red, line width=0.3mm] (150:3cm) -- (360:3cm); 
                     \draw[cyan, line width=0.3mm] (120:3cm) -- (210:3cm);
\draw[red, line width=0.3mm] (240:3cm) -- (360:3cm);

  \node at (120:3cm) [left] {$c$};
                    \node at (210:3cm) [left] {$d$};
\node at (150:3cm) [left] {$a$};
                    \node at (240:3cm) [left] {$\Bar{b}$};
                      
                    \end{scope}
                \end{tikzpicture}
                               \caption{On the left, the two $\theta$-orbits $[a,b]$ and $[c,d]$. On the right, their restrictions.}
                               \label{fig:enter-label}
                           \end{figure}
                           \item [b)] $c$ is an endpoint of $\tau_n$.
                            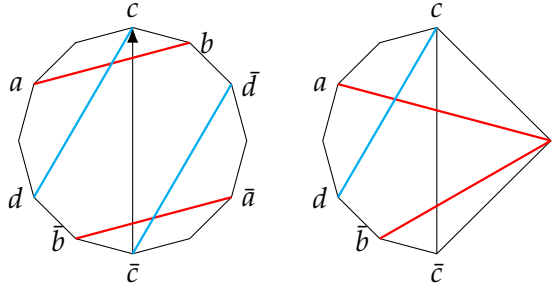
\begin{figure}[H]
                               \centering
                \begin{tikzpicture}[scale=0.5]
                    \draw (90:3cm) -- (120:3cm) -- (150:3cm) -- (180:3cm) -- (210:3cm) -- (240:3cm) -- (270:3cm) -- (300:3cm) -- (330:3cm) -- (360:3cm) -- (30:3cm) -- (60:3cm) --  cycle;
                    %\draw (90:3cm) -- node[midway, above, xshift=-1mm, yshift=-1mm] {} (180:3cm);
                    %\draw (90:3cm) -- node[midway, above left,xshift=1.2mm] {} (225:3cm);
                     \draw[-{Latex[length=2mm]}] (270:3cm) -- node[midway, above left,xshift=1mm] {} (90:3cm);
                   % \draw (90:3cm) -- node[midway, above right, xshift=-1mm] {} (315:3cm);
                    %\draw (90:3cm) -- node[midway, above right, xshift=-1mm] {} (360:3cm);
                     
                     \draw[red, line width=0.3mm] (240:3cm) -- (330:3cm); 
                     \draw[red, line width=0.3mm] (150:3cm) -- (60:3cm); 
                     \draw[cyan, line width=0.3mm] (90:3cm) -- (210:3cm);
                    \draw[cyan, line width=0.3mm] (270:3cm) -- (30:3cm);

                    \node at (90:3cm) [above] {$c$};
                    \node at (210:3cm) [left] {$d$};
                    \node at (270:3cm) [below] {$\Bar{c}$};
                    \node at (240:3cm) [left] {$\Bar{b}$};

\node at (150:3cm) [left] {$a$};
                    \node at (60:3cm) [right] {$b$};
                    \node at (330:3cm) [right] {$\Bar{a}$};
                    \node at (30:3cm) [right] {$\Bar{d}$};

                    \begin{scope}[xshift=8cm]
                      \draw (90:3cm) -- (120:3cm) -- (150:3cm) -- (180:3cm) -- (210:3cm) -- (240:3cm) -- (270:3cm) -- (360:3cm)  --  cycle;
                    %\draw (90:3cm) -- node[midway, above, xshift=-1mm, yshift=-1mm] {} (180:3cm);
                    %\draw (90:3cm) -- node[midway, above left,xshift=1.2mm] {} (225:3cm);
                    \draw (90:3cm) -- node[midway, above left,xshift=1mm] {} (270:3cm);
                   % \draw (90:3cm) -- node[midway, above right, xshift=-1mm] {} (315:3cm);
                    %\draw (90:3cm) -- node[midway, above right, xshift=-1mm] {} (360:3cm);
                    \draw[red, line width=0.3mm] (150:3cm) -- (360:3cm); 
                     \draw[cyan, line width=0.3mm] (90:3cm) -- (210:3cm);
\draw[red, line width=0.3mm] (240:3cm) -- (360:3cm);

                    \node at (90:3cm) [above] {$c$};
                    \node at (210:3cm) [left] {$d$};
                    \node at (270:3cm) [below] {$\Bar{c}$};
                    \node at (240:3cm) [left] {$\Bar{b}$};

\node at (150:3cm) [left] {$a$};
                    \end{scope}
                \end{tikzpicture}
                               \caption{On the left, the two $\rho$-orbits $[a,b]$ and $[c,d]$. On the right, their restrictions.}
                               \label{fig:enter-label}
                           \end{figure}
\end{itemize}

In case a), we have the following exchange relation:
\begin{equation}\label{thm1_eqn3}
u_{i_1}x_{ab}=\bold{y}^{D\bold{d}_{ac,d\ast}}x_{ad}x_{bc}+\bold{y}^{D\bold{d}_{ad,c\ast}}x_{ac}x_{bd}=\bold{y}^{\bold{d}_{ac,d\ast}}x_{ad}x_{bc}+\bold{y}^{\bold{d}_{ad,c\ast}}x_{ac}x_{bd},
\end{equation}
where the last equality is due to the fact that the $n$-th coordinate of $\bold{d}_{ab,c\ast}$ and $\bold{d}_{a\ast,bc}$ must be 0, since $L_n$ cannot cross both $(a,c)$ and $(d,\ast)$, nor both $(a,d)$ and $(c,\ast)$. 
It follows from \ref{thm1_eqn3} that
\begin{equation}
F_{ab}=\bold{y}^{\bold{d}_{ac,d\ast}}F_{bc}+\bold{y}^{\bold{d}_{ad,c\ast}}F_{bd},
\end{equation}
where we have used that $F_{ad}=F_{ac}=1$, since $[a,d]$ and $[a,c]$ must be either boundary edges or pairs of diagonals of $T$.

By inductive hypothesis and Proposition \ref{up:skein1},
\begin{align*}
    F_{ab}&=\bold{y}^{\bold{d}_{ac,d\ast}}(F_{(\Bar{b},\ast)}F_{(c,\ast)}-\bold{y}^{\bold{d}_{\Bar{b}\ast,c\ast}}F_{(c,\Bar{b})})+\bold{y}^{\bold{d}_{ad,c\ast}}(F_{(\Bar{b},\ast)}F_{(d,\ast)}-\bold{y}^{\bold{d}_{\Bar{b}\ast,d\ast}}F_{(d,\Bar{b})})\\
    &=F_{(\Bar{b},\ast)}(\bold{y}^{\bold{d}_{ac,d\ast}}F_{(c,\ast)}+\bold{y}^{\bold{d}_{ad,c\ast}}F_{(d,\ast)})-\bold{y}^{\bold{d}_{a\ast,\Bar{b}\ast}}(\bold{y}^{\bold{d}_{ac,d\Bar{b}}}F_{(c,\Bar{b})}+\bold{y}^{\bold{d}_{ad,c\Bar{b}}}F_{(d,\Bar{b})})\\
    &=F_{(\Bar{b},\ast)}F_{(a,\ast)}-\bold{y}^{\bold{d}_{a\ast,\Bar{b}\ast}}F_{(a,\Bar{b})}=F_{ab}^B.
\end{align*}
\begin{comment}
   It also follows from \ref{thm1_eqn3} that
\begin{equation}
    \bold{e}_{i_1}+\bold{g}_{ab}=\begin{cases}
        \bold{g}_{ad}+\bold{g}_{bc}\hspace{1cm} \text{if $\bold{y}^{\bold{d}_{ac,d\ast}}=1$}\\
        \bold{g}_{ac}+\bold{g}_{bd}\hspace{1cm} \text{otherwise}.
    \end{cases}
\end{equation}

By inductive hypothesis and Proposition \ref{up:skein1},

\begin{align*}
    \bold{g}_{ab}&=\begin{cases}
        -\bold{e}_{i_1} + \bold{g}_{(a,d)} + D(\bold{g}_{(\Bar{b},\ast)} + \bold{g}_{(c,\ast)} + \bold{e}_n)  \hspace{0.8cm}\text{if $\bold{y}^{\bold{d}_{ac,d\ast}}=1$}\\
        -\bold{e}_{i_1}+ \bold{g}_{(a,c)} + D(\bold{g}_{(\Bar{b},\ast)}+\bold{g}_{(d,\ast)}+ \bold{e}_n) \hspace{0.8cm} \text{otherwise};
        \end{cases}\\
        &=\begin{cases}
      D(-\bold{e}_{i_1}+\bold{g}_{(a,d)}+\bold{g}_{(\Bar{b},\ast)}+\bold{g}_{(c,\ast)}+\bold{e}_n) \hspace{0.8cm}\text{if $\bold{y}^{\bold{d}_{ac,d\ast}}=1$}\\
      D(-\bold{e}_{i_1}+\bold{g}_{(a,c)}+\bold{g}_{(\Bar{b},\ast)}+\bold{g}_{(d,\ast)}+\bold{e}_n)\hspace{0.8cm} \text{otherwise}.
    \end{cases}\\
    &=D(\bold{g}_{(a,\ast)}+\bold{g}_{(\Bar{b},\ast)}+\bold{e}_n)=\bold{g}_{ab}^B.
\end{align*} 
\end{comment}

On the other hand, in case b), we have the following exchange relation:
\begin{equation}\label{thm1_eqn4}
u_{i_1}x_{ab}=\bold{y}^{D\bold{d}_{ac,d\ast}}x_{ad}x_{bc}+\bold{y}^{D\bold{d}_{ad,c\ast}}x_{ac}x_{bd}=\bold{y}^{D\bold{d}_{ac,d\ast}}x_{ad}x_{bc}+\bold{y}^{\bold{d}_{ad,c\ast}}x_{ac}x_{bd},
\end{equation}
where the last equality is due to the fact that the $n$-th coordinate of $\bold{d}_{ad,c\ast}$ must be 0, since $L_n$ cannot cross both $(a,d)$ and $(c,\ast)$.
 
It follows from \ref{thm1_eqn4} that
\begin{equation}
F_{ab}=\bold{y}^{D\bold{d}_{ac,d\ast}}F_{bc}+\bold{y}^{\bold{d}_{ad,c\ast}}F_{bd},
\end{equation}
where we have used that $F_{ad}=F_{ac}=1$, since $[a,d]$ and $[a,c]$ must be either boundary edges or pairs of diagonals of $T$.

By inductive hypothesis and repeated applications of Proposition \ref{up:skein1},
\begin{align*}
   F_{ab}&=\bold{y}^{D\bold{d}_{ac,d\ast}}F_{(\Bar{b},\Bar{c})}+\bold{y}^{\bold{d}_{ad,c\ast}}(F_{(\Bar{b},\ast)}F_{(d,\ast)}-\bold{y}^{\bold{d}_{\Bar{b}\ast,d\ast}}F_{(d,\Bar{b})})\\
   &=F_{(a,\ast)}F_{(\Bar{b},\ast)}-\bold{y}^{\bold{d}_{a\ast,\Bar{b}\ast}}F_{(a,\Bar{b})}.
\end{align*}  

\begin{comment}
    
It also follows from \ref{thm1_eqn4} that
\begin{equation}
    \bold{e}_{i_1}+\bold{g}_{ab}=\begin{cases}
        \bold{g}_{ad}+\bold{g}_{bc}\hspace{1cm} \text{if $\bold{y}^{\bold{d}_{ac,d\ast}}=1$}\\
        \bold{g}_{ac}+\bold{g}_{bd}\hspace{1cm} \text{otherwise}.
    \end{cases}
\end{equation}

By inductive hypothesis and Proposition \ref{up:skein1},
\begin{align*}
    \bold{g}_{ab}&=\begin{cases}
        -\bold{e}_{i_1} + \bold{g}_{(a,d)}+ D\bold{g}_{(\Bar{b},\Bar{c})} \hspace{3.2cm}\text{if $\bold{y}^{\bold{d}_{ac,d\ast}}=1$}\\
        -\bold{e}_{i_1} + \bold{g}_{(a,c)} + D(\bold{g}_{(\Bar{b},\ast)}+\bold{g}_{(d,\ast)}+ \bold{e}_n) \hspace{0.8cm} \text{otherwise};
        \end{cases}\\
&=\begin{cases}
      D(-\bold{e}_{i_1} + \bold{g}_{(a,d)}+ \bold{g}_{(\Bar{b},\Bar{c})}) \hspace{2.8cm}\text{if $\bold{y}^{\bold{d}_{ac,d\ast}}=1$}\\
      D(-\bold{e}_{i_1} + \bold{g}_{(a,c)} + \bold{g}_{(\Bar{b},\ast)}+\bold{g}_{(d,\ast)}+ \bold{e}_n)\hspace{0.8cm} \text{otherwise}.
    \end{cases}\\
&=D(\bold{g}_{(a,\ast)}+\bold{g}_{(\Bar{b},\ast)}+\bold{e}_n)=\bold{g}_{ab}^B.
\end{align*}
\end{comment}

    \end{itemize}
    Similarly one can prove that $\bold{g}_{ab}=\bold{g}_{ab}^B$.
\end{proof}
\begin{comment}

\begin{proof}
    We prove the theorem by induction on the distance $d\geq 1$ of the cluster $(\tilde{x'}, \tilde{B'})$ from the initial cluster $(\tilde{\bold{x}},\tilde{B})$ in the exchange graph of $\mathcal{A}$. 

    If $d=1$, $(\tilde{x'}, \tilde{B'})=(\tilde{x}\setminus \{x_k\} \cup \{x_{k'}\}, \tilde{B'})$, where $x_{k'}$ corresponds to a $\rho$-orbit $[a,b]$, which is obtained from the initial triangulation $T$ by applying one flip. Therefore, each diagonal of $[a,b]$ crosses only one diagonal of $T$. The theorem holds by Lemma \ref{lemma1}.

    Assume now the theorem true for every $d \le n$, i.e., the $F$-polynomial (resp. the $\bold{g}$-vector) of each cluster variable in each cluster, at a distance not greater than $n$ from the initial one, is given by $F^B$ (resp. $g^B$). Let $(\tilde{\bold{x'}},\tilde{B'})$ be a cluster such that $d=n$. We mutate it once. If $d$ decreases, we have nothing to prove. Otherwise, we conclude by Lemma \ref{lemma2}.
\end{proof}
    
\end{comment}
\subsubsection{Type C}
\begin{definition}
    Let $[a,b]$ be an orbit of the action of $\theta$ on the diagonals of $\Poly_{2n+2}$. We define the $rotated$ $restriction$ $of$ $[a,b]$, and we denote it by $\tilde{\text{Res}}([a,b])$, as follows. 
    
    \begin{itemize}
        \item If $[a,b]=[a,\Bar{a}]$ is a diameter, so $\text{Res}([a,\Bar{a}])=\{\gamma\}$, then $\tilde{\text{Res}}([a,\Bar{a}]):=\{\tilde\gamma_1, \tilde\gamma_2\}$, where $\tilde\gamma_1=\gamma$ and $\tilde\gamma_2$, if it exists, is the diagonal of $\Poly_{n+3}$ which intersects the same diagonals of $T$ as $\gamma$ but $d$. If there is no such diagonal, $\tilde{\text{Res}}([a,\Bar{a}]):=\{\tilde\gamma_1\}$. A possible situation is represented in Figure \ref{res_tilde_1}.

         \begin{figure}[H]
                               \centering
                \begin{tikzpicture}[scale=0.5]
                    \draw (90:3cm) -- (120:3cm) -- (150:3cm) -- (180:3cm) -- (210:3cm) -- (240:3cm) -- (270:3cm) -- (300:3cm) -- (330:3cm) -- (360:3cm) -- (30:3cm) -- (60:3cm) --  cycle;
                    %\draw (90:3cm) -- node[midway, above, xshift=-1mm, yshift=-1mm] {} (180:3cm);
                    %\draw (90:3cm) -- node[midway, above left,xshift=1.2mm] {} (225:3cm);
                    \draw[-{Latex[length=2mm]}] (270:3cm) -- node[midway, above left,xshift=1mm] {} (90:3cm);
                   % \draw (90:3cm) -- node[midway, above right, xshift=-1mm] {} (315:3cm);
                    %\draw (90:3cm) -- node[midway, above right, xshift=-1mm] {} (360:3cm);

                    \draw[black, line width=0.3mm] (150:3cm) -- (330:3cm);
                    %\draw[green, line width=0.3mm] (150:3cm) -- (90:3cm);
                    %\draw[green, line width=0.3mm] (90:3cm) -- (30:3cm);
                    %\draw[yellow, line width=0.3mm] (150:3cm) -- (240:3cm);
                    %\draw[yellow, line width=0.3mm] (30:3cm) -- (300:3cm);

                    %\node at (120:3cm) [left] {$a$};
                
                    %\node at (180:3cm) [left] {$\hat{c}$};
                    %\node at (360:3cm) [right] {$c$};
                    %\node at (60:3cm) [right] {$\hat{a}$};

\node at (150:3cm) [left] {$a$};
                    %\node at (210:3cm) [left] {$\hat{b}$};
                    %\node at (330:3cm) [right] {$b$};
                 
                    \node at (330:3cm) [right] {$b=\Bar{a}$};

                    \begin{scope}[xshift=8cm]
                      \draw (90:3cm) -- (120:3cm) -- (150:3cm) -- (180:3cm) -- (210:3cm) -- (240:3cm) -- (270:3cm) -- (360:3cm)  --  cycle;
                    %\draw (90:3cm) -- node[midway, above, xshift=-1mm, yshift=-1mm] {} (180:3cm);
                    %\draw (90:3cm) -- node[midway, above left,xshift=1.2mm] {} (225:3cm);
                    \draw (90:3cm) -- node[midway, above left,xshift=1mm] {} (270:3cm);
                   % \draw (90:3cm) -- node[midway, above right, xshift=-1mm] {} (315:3cm);
                    %\draw (90:3cm) -- node[midway, above right, xshift=-1mm] {} (360:3cm);
                    %\draw[red, line width=0.3mm] (150:3cm) -- (360:3cm); 
                     %\draw[cyan, line width=0.3mm] (120:3cm) -- (360:3cm);
                     %\draw[cyan, line width=0.3mm] (270:3cm) -- (180:3cm);
                    
\draw[black, line width=0.3mm] (150:3cm) -- node[midway, above right, xshift=-2mm, yshift=-1mm] {$\textcolor{black}{\tilde\gamma_2}$}(270:3cm);
\draw[black, line width=0.3mm] (150:3cm) -- node[midway, above right, xshift=-2mm, yshift=-1mm] {$\textcolor{black}{\tilde\gamma_1}$}(360:3cm);
%\draw[black, line width=0.3mm] (152:3cm) -- (358:3cm);
%\draw[black, line width=0.3mm] (152:3cm) -- (268:3cm);

                    \node at (360:3cm) [right] {$\ast$};

\node at (150:3cm) [left] {$a$};

                    \end{scope}
                \end{tikzpicture}
                               \caption{On the left, a diameter $[a,\Bar{a}]$. On the right, its rotated restriction.}
                               \label{res_tilde_1}
                           \end{figure}
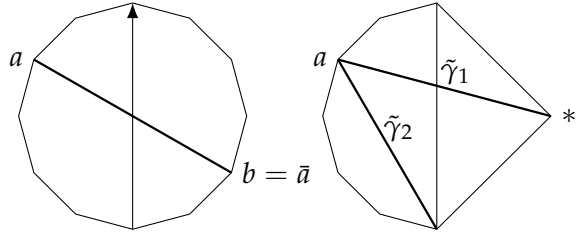

        \item If $[a,b]$ is a pair of diagonals which do not cross $d$, then $\tilde{\text{Res}}([a,b]):=\text{Res}([a,b])$. A possible situation is represented in Figure \ref{res_tilde_2}.

         \begin{figure}[H]
                               \centering
                \begin{tikzpicture}[scale=0.5]
                    \draw (90:3cm) -- (120:3cm) -- (150:3cm) -- (180:3cm) -- (210:3cm) -- (240:3cm) -- (270:3cm) -- (300:3cm) -- (330:3cm) -- (360:3cm) -- (30:3cm) -- (60:3cm) --  cycle;
                    %\draw (90:3cm) -- node[midway, above, xshift=-1mm, yshift=-1mm] {} (180:3cm);
                    %\draw (90:3cm) -- node[midway, above left,xshift=1.2mm] {} (225:3cm);
                    \draw[-{Latex[length=2mm]}] (270:3cm) -- node[midway, above left,xshift=1mm] {} (90:3cm);
                   % \draw (90:3cm) -- node[midway, above right, xshift=-1mm] {} (315:3cm);
                    %\draw (90:3cm) -- node[midway, above right, xshift=-1mm] {} (360:3cm);

                    \draw[black, line width=0.3mm] (150:3cm) -- (270:3cm);
                    \draw[black, line width=0.3mm] (90:3cm) -- (330:3cm);
                    %\draw[green, line width=0.3mm] (150:3cm) -- (90:3cm);
                    %\draw[green, line width=0.3mm] (90:3cm) -- (30:3cm);
                    %\draw[yellow, line width=0.3mm] (150:3cm) -- (240:3cm);
                    %\draw[yellow, line width=0.3mm] (30:3cm) -- (300:3cm);

                    %\node at (120:3cm) [left] {$a$};
                
                    %\node at (180:3cm) [left] {$\hat{c}$};
                    %\node at (360:3cm) [right] {$c$};
                    %\node at (60:3cm) [right] {$\hat{a}$};

\node at (150:3cm) [left] {$a$};
                    %\node at (210:3cm) [left] {$\hat{b}$};
                    %\node at (330:3cm) [right] {$b$};
                 
                    \node at (330:3cm) [right] {$\Bar{a}$};
                    \node at (270:3cm) [below] {$b$};
\node at (90:3cm) [above] {$\Bar{b}$};
                    \begin{scope}[xshift=8cm]
                      \draw (90:3cm) -- (120:3cm) -- (150:3cm) -- (180:3cm) -- (210:3cm) -- (240:3cm) -- (270:3cm) -- (360:3cm)  --  cycle;
                    %\draw (90:3cm) -- node[midway, above, xshift=-1mm, yshift=-1mm] {} (180:3cm);
                    %\draw (90:3cm) -- node[midway, above left,xshift=1.2mm] {} (225:3cm);
                    \draw (90:3cm) -- node[midway, above left,xshift=1mm] {} (270:3cm);
                   % \draw (90:3cm) -- node[midway, above right, xshift=-1mm] {} (315:3cm);
                    %\draw (90:3cm) -- node[midway, above right, xshift=-1mm] {} (360:3cm);
                    %\draw[red, line width=0.3mm] (150:3cm) -- (360:3cm); 
                     %\draw[cyan, line width=0.3mm] (120:3cm) -- (360:3cm);
                     %\draw[cyan, line width=0.3mm] (270:3cm) -- (180:3cm);
                    
%\draw[black, line width=0.3mm] (150:3cm) -- node[midway, above right, xshift=-2mm, yshift=-1mm] {$\textcolor{black}{\tilde\gamma_2}$}(270:3cm);
\draw[black, line width=0.3mm] (150:3cm) -- node[midway, above right, xshift=-2mm, yshift=-1mm] {}(270:3cm);
%\draw[black, line width=0.3mm] (152:3cm) -- (358:3cm);
%\draw[black, line width=0.3mm] (152:3cm) -- (268:3cm);

                    \node at (360:3cm) [right] {$\ast$};
\node at (270:3cm) [below] {$b$};
\node at (150:3cm) [left] {$a$};

                    \end{scope}
                \end{tikzpicture}
                               \caption{On the left, a $\theta$-orbit $[a,b]$. On the right, its rotated restriction.}
                               \label{res_tilde_2}
                           \end{figure}
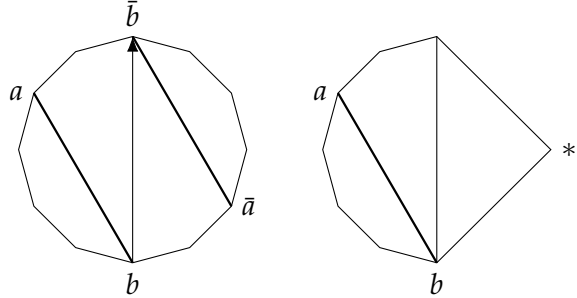

        \item If $[a,b]$ is a pair of diagonals which cross $d$, then $\text{Res}([a,b])=\{\gamma_1,\gamma_2\}$, where $\gamma_1$ and $\gamma_2$ are two diagonals of $\Poly_{n+3}$ that share the right endpoint, and such that $\gamma_2$ is obtained from $\gamma_1$ by rotating counterclockwise (resp. clockwise) its left endpoint if $\tau_{n-1}$ is counterclockwise (resp. clockwise) from $\tau_n$. We define $\tilde{\text{Res}}([a,b]):=\{\tilde\gamma_1, \tilde\gamma_2\}$, where  $\tilde\gamma_1=\gamma_1$ and $\tilde\gamma_2$, if it exists, is the diagonal of $\Poly_{n+3}$ which intersects the same diagonals of $T$ as $\gamma_2$ but the diameter. If there is no such diagonal, $\tilde{\text{Res}}([a,b]):=\{\tilde\gamma_1\}$. A possible situation is represented in Figure \ref{res_tilde_3}.

         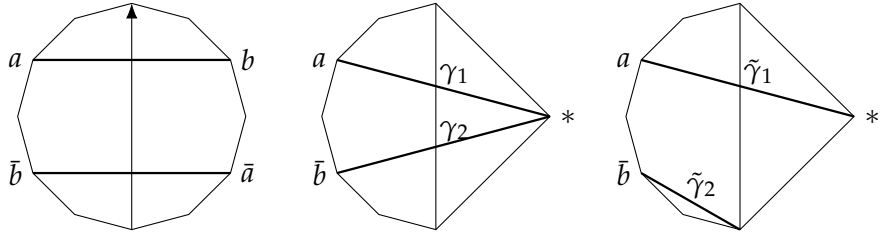
\begin{figure}[H]
                               \centering
                \begin{tikzpicture}[scale=0.5]
                    \draw (90:3cm) -- (120:3cm) -- (150:3cm) -- (180:3cm) -- (210:3cm) -- (240:3cm) -- (270:3cm) -- (300:3cm) -- (330:3cm) -- (360:3cm) -- (30:3cm) -- (60:3cm) --  cycle;
                    %\draw (90:3cm) -- node[midway, above, xshift=-1mm, yshift=-1mm] {} (180:3cm);
                    %\draw (90:3cm) -- node[midway, above left,xshift=1.2mm] {} (225:3cm);
                    \draw[-{Latex[length=2mm]}] (270:3cm) -- node[midway, above left,xshift=1mm] {} (90:3cm);
                   % \draw (90:3cm) -- node[midway, above right, xshift=-1mm] {} (315:3cm);
                    %\draw (90:3cm) -- node[midway, above right, xshift=-1mm] {} (360:3cm);

                    \draw[black, line width=0.3mm] (150:3cm) -- (30:3cm);
                    \draw[black, line width=0.3mm] (210:3cm) -- (330:3cm);
                    %\draw[green, line width=0.3mm] (150:3cm) -- (90:3cm);
                    %\draw[green, line width=0.3mm] (90:3cm) -- (30:3cm);
                    %\draw[yellow, line width=0.3mm] (150:3cm) -- (240:3cm);
                    %\draw[yellow, line width=0.3mm] (30:3cm) -- (300:3cm);

                    %\node at (120:3cm) [left] {$a$};
                
                    %\node at (180:3cm) [left] {$\hat{c}$};
                    %\node at (360:3cm) [right] {$c$};
                    %\node at (60:3cm) [right] {$\hat{a}$};

\node at (150:3cm) [left] {$a$};
                    %\node at (210:3cm) [left] {$\hat{b}$};
                    %\node at (330:3cm) [right] {$b$};
                 
                    \node at (330:3cm) [right] {$\Bar{a}$};
                    \node at (30:3cm) [right] {$b$};
                    \node at (210:3cm) [left] {$\Bar{b}$};

                    \begin{scope}[xshift=8cm]
                      \draw (90:3cm) -- (120:3cm) -- (150:3cm) -- (180:3cm) -- (210:3cm) -- (240:3cm) -- (270:3cm) -- (360:3cm)  --  cycle;
                    %\draw (90:3cm) -- node[midway, above, xshift=-1mm, yshift=-1mm] {} (180:3cm);
                    %\draw (90:3cm) -- node[midway, above left,xshift=1.2mm] {} (225:3cm);
                    \draw (90:3cm) -- node[midway, above left,xshift=1mm] {} (270:3cm);
                   % \draw (90:3cm) -- node[midway, above right, xshift=-1mm] {} (315:3cm);
                    %\draw (90:3cm) -- node[midway, above right, xshift=-1mm] {} (360:3cm);
                    %\draw[red, line width=0.3mm] (150:3cm) -- (360:3cm); 
                     %\draw[cyan, line width=0.3mm] (120:3cm) -- (360:3cm);
                     %\draw[cyan, line width=0.3mm] (270:3cm) -- (180:3cm);
                    
\draw[black, line width=0.3mm] (210:3cm) -- node[midway, above right, xshift=-2mm, yshift=-1mm] {$\textcolor{black}{\gamma_2}$}(360:3cm);
\draw[black, line width=0.3mm] (150:3cm) -- node[midway, above right, xshift=-2mm, yshift=-1mm] {$\textcolor{black}{\gamma_1}$}(360:3cm);
%\draw[black, line width=0.3mm] (152:3cm) -- (358:3cm);
%\draw[black, line width=0.3mm] (152:3cm) -- (268:3cm);

                    \node at (360:3cm) [right] {$\ast$};

\node at (150:3cm) [left] {$a$};
\node at (210:3cm) [left] {$\Bar{b}$};

                    \end{scope}
                    
                    \begin{scope}[xshift=16cm]
                      \draw (90:3cm) -- (120:3cm) -- (150:3cm) -- (180:3cm) -- (210:3cm) -- (240:3cm) -- (270:3cm) -- (360:3cm)  --  cycle;
                    %\draw (90:3cm) -- node[midway, above, xshift=-1mm, yshift=-1mm] {} (180:3cm);
                    %\draw (90:3cm) -- node[midway, above left,xshift=1.2mm] {} (225:3cm);
                    \draw (90:3cm) -- node[midway, above left,xshift=1mm] {} (270:3cm);
                   % \draw (90:3cm) -- node[midway, above right, xshift=-1mm] {} (315:3cm);
                    %\draw (90:3cm) -- node[midway, above right, xshift=-1mm] {} (360:3cm);
                    %\draw[red, line width=0.3mm] (150:3cm) -- (360:3cm); 
                     %\draw[cyan, line width=0.3mm] (120:3cm) -- (360:3cm);
                     %\draw[cyan, line width=0.3mm] (270:3cm) -- (180:3cm);
                    
\draw[black, line width=0.3mm] (210:3cm) -- node[midway, above right, xshift=-2mm, yshift=-1mm] {$\textcolor{black}{\tilde\gamma_2}$}(270:3cm);
\draw[black, line width=0.3mm] (150:3cm) -- node[midway, above right, xshift=-2mm, yshift=-1mm] {$\textcolor{black}{\tilde\gamma_1}$}(360:3cm);
%\draw[black, line width=0.3mm] (152:3cm) -- (358:3cm);
%\draw[black, line width=0.3mm] (152:3cm) -- (268:3cm);

                    \node at (360:3cm) [right] {$\ast$};

\node at (150:3cm) [left] {$a$};
\node at (210:3cm) [left] {$\Bar{b}$};

                    \end{scope}
                \end{tikzpicture}
                               \caption{From left to right, a $\theta$-orbit $[a,b]$, its restriction and its rotated restriction.}
                               \label{res_tilde_3}
                           \end{figure}

    \end{itemize}
 
\end{definition}

\begin{definition}
    Let $v \in \mathbb{Z}_{\geq 0}^{2n-1}$. We define the $rotated$ $restriction$ $of$ $v$, and we denote it by $\tilde{\text{Res}}(v)$,  as the vector of the first $n$ coordinates of $v$, with the $n$-th one divided by 2.
\end{definition}

\begin{definition}\label{def_type_C}
  Let $[a,b] \not \subset T$ be an orbit of the action of $\theta$ on the diagonals of $\Poly_{2n+2}$. If $\tilde{\text{Res}}([a,b])=\{\tilde\gamma\}$ contains only one diagonal $\tilde\gamma$ of $\Poly_{n+3}$,  we define
  \begin{equation}
      F_{ab}^C=F_{\Tilde{\gamma}},
\end{equation}
\begin{equation}
      \bold{g}_{ab}^C=\begin{cases}
          \text{$\bold{g}_{\Tilde{\gamma}}+\bold{e}_i$ \hspace{0.5cm} if $\tau_i$ and $\tau_n$ are two different sides of a triangle of $T$},\\ \hspace{1.9cm}\text{$\tau_i$ is clockwise from $\tau_n$, and $\tilde\gamma$ crosses $\tau_n$;}\\
          \text{$\bold{g}_{\tilde\gamma}$ \hspace{1.4cm}otherwise.}
      \end{cases}
  \end{equation} 
  
  Otherwise there are two cases to consider:
  \begin{itemize}
      \item $(a,b)=(a,\Bar{a})$ is a diameter. Then $\tilde{\text{Res}}([a,\Bar{a}])=\{\tilde\gamma_1, \tilde\gamma_2\}$, and there are uniquely determined two $\theta$-orbits $[a,\Bar{c}]$ and $[a,\Bar{b}]$, such that $\tilde{\text{Res}}([a,\Bar{c}])=\{\tilde\gamma_1\}$ and $\tilde{\text{Res}}([a,\Bar{b}])=\{\tilde\gamma_2\}$. A possible situation is represented in Figure \ref{type$C$-1}.
 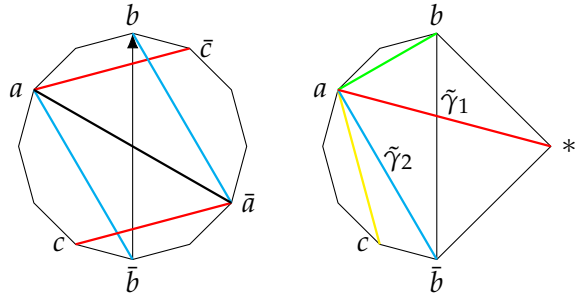
\begin{figure}[H]
                               \centering
                \begin{tikzpicture}[scale=0.5]
                    \draw (90:3cm) -- (120:3cm) -- (150:3cm) -- (180:3cm) -- (210:3cm) -- (240:3cm) -- (270:3cm) -- (300:3cm) -- (330:3cm) -- (360:3cm) -- (30:3cm) -- (60:3cm) --  cycle;
                    %\draw (90:3cm) -- node[midway, above, xshift=-1mm, yshift=-1mm] {} (180:3cm);
                    %\draw (90:3cm) -- node[midway, above left,xshift=1.2mm] {} (225:3cm);
                    \draw[-{Latex[length=2mm]}] (270:3cm) -- node[midway, above left,xshift=1mm] {} (90:3cm);
                   % \draw (90:3cm) -- node[midway, above right, xshift=-1mm] {} (315:3cm);
                    %\draw (90:3cm) -- node[midway, above right, xshift=-1mm] {} (360:3cm);
                     
                    \draw[cyan, line width=0.3mm] (150:3cm) -- (270:3cm);
                     \draw[red, line width=0.3mm] (150:3cm) -- (60:3cm); 
                     \draw[red, line width=0.3mm] (330:3cm) -- (240:3cm);
                     \draw[cyan, line width=0.3mm] (330:3cm) -- (90:3cm);
                    \draw[black, line width=0.3mm] (150:3cm) -- (330:3cm);
                    %\draw[green, line width=0.3mm] (150:3cm) -- (90:3cm);
                    %\draw[green, line width=0.3mm] (90:3cm) -- (30:3cm);
                    %\draw[yellow, line width=0.3mm] (150:3cm) -- (240:3cm);
                    %\draw[yellow, line width=0.3mm] (30:3cm) -- (300:3cm);

                    %\node at (120:3cm) [left] {$a$};
                    \node at (90:3cm) [above,yshift=-1] {$b$};
                    \node at (270:3cm) [below,yshift=1] {$\Bar{b}$};
                    %\node at (180:3cm) [left] {$\hat{c}$};
                    %\node at (360:3cm) [right] {$c$};
                    %\node at (60:3cm) [right] {$\hat{a}$};

\node at (150:3cm) [left] {$a$};
                    %\node at (210:3cm) [left] {$\hat{b}$};
                    %\node at (330:3cm) [right] {$b$};
                    \node at (60:3cm) [right] {$\Bar{c}$};
                    \node at (240:3cm) [left] {$c$};
                    \node at (330:3cm) [right] {$\Bar{a}$};

                    \begin{scope}[xshift=8cm]
                      \draw (90:3cm) -- (120:3cm) -- (150:3cm) -- (180:3cm) -- (210:3cm) -- (240:3cm) -- (270:3cm) -- (360:3cm)  --  cycle;
                    %\draw (90:3cm) -- node[midway, above, xshift=-1mm, yshift=-1mm] {} (180:3cm);
                    %\draw (90:3cm) -- node[midway, above left,xshift=1.2mm] {} (225:3cm);
                    \draw (90:3cm) -- node[midway, above left,xshift=1mm] {} (270:3cm);
                   % \draw (90:3cm) -- node[midway, above right, xshift=-1mm] {} (315:3cm);
                    %\draw (90:3cm) -- node[midway, above right, xshift=-1mm] {} (360:3cm);
                    %\draw[red, line width=0.3mm] (150:3cm) -- (360:3cm); 
                     %\draw[cyan, line width=0.3mm] (120:3cm) -- (360:3cm);
                     %\draw[cyan, line width=0.3mm] (270:3cm) -- (180:3cm);
                     \draw[green, line width=0.3mm] (150:3cm) -- (90:3cm);
%\draw[red, line width=0.3mm] (210:3cm) -- (270:3cm);
                   \draw[yellow, line width=0.3mm] (150:3cm) -- (240:3cm);
\draw[cyan, line width=0.3mm] (150:3cm) -- node[midway, above right, xshift=-2mm, yshift=-1mm] {$\textcolor{black}{\tilde\gamma_2}$}(270:3cm);
\draw[red, line width=0.3mm] (150:3cm) -- node[midway, above right, xshift=-2mm, yshift=-1mm] {$\textcolor{black}{\tilde\gamma_1}$}(360:3cm);
%\draw[black, line width=0.3mm] (152:3cm) -- (358:3cm);
%\draw[black, line width=0.3mm] (152:3cm) -- (268:3cm);

                    \node at (360:3cm) [right] {$\ast$};
 \node at (240:3cm) [left] {$c$};
\node at (150:3cm) [left] {$a$};
\node at (90:3cm) [above,yshift=-1] {$b$};
                    \node at (270:3cm) [below,yshift=1] {$\Bar{b}$};

                    \end{scope}
                \end{tikzpicture}
                               \caption{On the left, the $\theta$-orbits $[a,\Bar{a}]$, $[a,\Bar{c}]$, $[a,\Bar{b}]$. On the right, their rotated restrictions, and the diagonals $(a,b)$ and $(a,c)$.}
                               \label{type$C$-1}
                           \end{figure}     
      
      We define
      \begin{equation}
        F_{a\Bar{a}}^C=F_{\tilde\gamma_1}F_{\tilde\gamma_2}- \bold{y}^{\Tilde{\text{Res}}(\bold{d}_{a\ast,c\Bar{b}}+\bold{d}_{a\Bar{b},b\ast})}F_{(a,b)}F_{(a,c)},
\end{equation}
\begin{equation}\label{eq_g_vect_1}
      \bold{g}_{a\Bar{a}}^C=\begin{cases}
          \text{$\bold{g}_{\tilde\gamma_1}+\bold{g}_{\tilde\gamma_2}+\bold{e}_i-\bold{g}_{(\Bar{b},c)}$  \hspace{0.6cm}if $\tau_i$ and $\tau_n$ are two different sides of a triangle of $T$},\\ \hspace{4.2cm}\text{and $\tau_i$ is clockwise from $\tau_n$;}\\
          \text{$\bold{g}_{\tilde\gamma_1}+\bold{g}_{\tilde\gamma_2}$ \hspace{2.6cm}otherwise.}
      \end{cases}
  \end{equation}
     
      \item $[a,b]$ is a pair of diagonals which cross $d$, and $\tilde{\text{Res}}([a,b])=\{\tilde\gamma_1,\tilde\gamma_2\}$, where $\tilde\gamma_1$ and $\tilde\gamma_2$ are two diagonals of $\Poly_{n+3}$. There are uniquely determined two $\theta$-orbits $[a,d]$ and $[b,c]$, such that $\tilde{\text{Res}}([a,d])=\{\tilde\gamma_1\}$ and $\tilde{\text{Res}}([b,c])=\{\tilde\gamma_2\}$. A possible situation is represented in Figure \ref{type$C$-2}.
     
         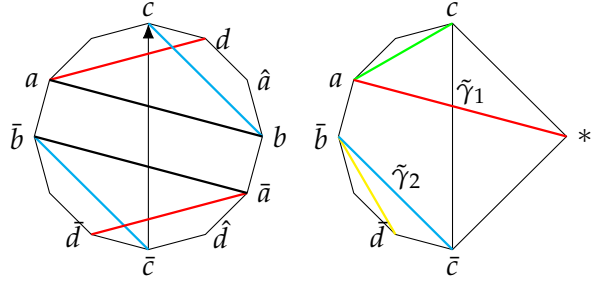
\begin{figure}[H]
                               \centering
                \begin{tikzpicture}[scale=0.5]
                    \draw (90:3cm) -- (120:3cm) -- (150:3cm) -- (180:3cm) -- (210:3cm) -- (240:3cm) -- (270:3cm) -- (300:3cm) -- (330:3cm) -- (360:3cm) -- (30:3cm) -- (60:3cm) --  cycle;
                    %\draw (90:3cm) -- node[midway, above, xshift=-1mm, yshift=-1mm] {} (180:3cm);
                    %\draw (90:3cm) -- node[midway, above left,xshift=1.2mm] {} (225:3cm);
                    \draw[-{Latex[length=2mm]}] (270:3cm) -- node[midway, above left,xshift=1mm] {} (90:3cm);
                   % \draw (90:3cm) -- node[midway, above right, xshift=-1mm] {} (315:3cm);
                    %\draw (90:3cm) -- node[midway, above right, xshift=-1mm] {} (360:3cm);
                     
                    \draw[cyan, line width=0.3mm] (180:3cm) -- (270:3cm);
                     \draw[red, line width=0.3mm] (150:3cm) -- (60:3cm); 
                     \draw[red, line width=0.3mm] (330:3cm) -- (240:3cm);
                     \draw[cyan, line width=0.3mm] (360:3cm) -- (90:3cm);
                    \draw[black, line width=0.3mm] (150:3cm) -- (360:3cm);
                    \draw[black, line width=0.3mm] (180:3cm) -- (330:3cm);
                    %\draw[green, line width=0.3mm] (150:3cm) -- (90:3cm);
                    %\draw[green, line width=0.3mm] (90:3cm) -- (30:3cm);
                    %\draw[yellow, line width=0.3mm] (210:3cm) -- (240:3cm);
                    %\draw[yellow, line width=0.3mm] (330:3cm) -- (300:3cm);

                    %\node at (120:3cm) [left] {$a$};
                    \node at (90:3cm) [above,yshift=-1] {$c$};
                    \node at (270:3cm) [below,yshift=1] {$\Bar{c}$};
                    \node at (180:3cm) [left] {$\Bar{b}$};
                    \node at (360:3cm) [right] {$b$};
                    \node at (330:3cm) [right] {$\Bar{a}$};
                    %\node at (60:3cm) [right] {$\hat{a}$};

\node at (150:3cm) [left] {$a$};
                    %\node at (210:3cm) [left] {$\hat{b}$};
                    %\node at (330:3cm) [right] {$b$};
                    \node at (300:3cm) [right,xshift=-1] {$\hat{d}$};
                    \node at (60:3cm) [right] {$d$};
                    \node at (240:3cm) [left] {$\Bar{d}$};
                    \node at (30:3cm) [right] {$\hat{a}$};

                    \begin{scope}[xshift=8cm]
                      \draw (90:3cm) -- (120:3cm) -- (150:3cm) -- (180:3cm) -- (210:3cm) -- (240:3cm) -- (270:3cm) -- (360:3cm)  --  cycle;
                    %\draw (90:3cm) -- node[midway, above, xshift=-1mm, yshift=-1mm] {} (180:3cm);
                    %\draw (90:3cm) -- node[midway, above left,xshift=1.2mm] {} (225:3cm);
                    \draw (90:3cm) -- node[midway, above left,xshift=1mm] {} (270:3cm);
                   % \draw (90:3cm) -- node[midway, above right, xshift=-1mm] {} (315:3cm);
                    %\draw (90:3cm) -- node[midway, above right, xshift=-1mm] {} (360:3cm);
                    %\draw[red, line width=0.3mm] (150:3cm) -- (360:3cm); 
                     %\draw[cyan, line width=0.3mm] (120:3cm) -- (360:3cm);
                     %\draw[cyan, line width=0.3mm] (270:3cm) -- (180:3cm);
                     \draw[green, line width=0.3mm] (150:3cm) -- (90:3cm);
%\draw[red, line width=0.3mm] (210:3cm) -- (270:3cm);
                    \draw[yellow, line width=0.3mm] (180:3cm) -- (240:3cm);
\draw[cyan, line width=0.3mm] (180:3cm) -- node[midway, above right, xshift=-2mm, yshift=-1mm] {$\textcolor{black}{\tilde\gamma_2}$}(270:3cm);
\draw[red, line width=0.3mm] (150:3cm) -- node[midway, above right, xshift=-2mm, yshift=-1mm] {$\textcolor{black}{\tilde\gamma_1}$}(360:3cm);
%\draw[black, line width=0.3mm] (152:3cm) -- node[midway, above right, xshift=-1mm, yshift=-1.5mm] {$\Tilde{\gamma_1}$} (358:3cm);
%\draw[black, line width=0.3mm] (212:3cm) -- node[midway, above right, xshift=-1mm, yshift= -2mm] {$\Tilde{\gamma_2}$} (268:3cm);

                    \node at (360:3cm) [right] {$\ast$};
 \node at (240:3cm) [left] {$\Bar{d}$};
\node at (150:3cm) [left] {$a$};
\node at (90:3cm) [above,yshift=-1] {$c$};
                    \node at (270:3cm) [below,yshift=1] {$\Bar{c}$};
                     \node at (180:3cm) [left] {$\Bar{b}$};
                      
                    \end{scope}
                \end{tikzpicture}
                               \caption{On the left, the $\theta$-orbits $[a,b]$, $[a,d]$, $[b,c]$. On the right, their rotated restrictions, and the diagonals $(a,c)$ and $(\Bar{b},\Bar{d})$.}
                               \label{type$C$-2}
                           \end{figure} 
 We define
      \begin{equation}
        F_{ab}^C=F_{\tilde\gamma_1}F_{\tilde\gamma_2}- \bold{y}^{\Tilde{\text{Res}}(\bold{d}_{\Bar{b}\ast,\Bar{d}\Bar{c}}+\bold{d}_{a\Bar{c},c\ast})}F_{(a,c)}F_{(\Bar{b},\Bar{d})},
\end{equation}
\begin{equation}\label{eq_g_vect_2}
      \bold{g}_{ab}^C=\begin{cases}
          \text{$\bold{g}_{\tilde\gamma_1}+\bold{g}_{\tilde\gamma_2}+\bold{e}_i-\bold{g}_{(\Bar{c},\Bar{d})}$  \hspace{0.6cm}if $\tau_i$ and $\tau_n$ are two different sides of a triangle of $T$},\\ \hspace{4.2cm}\text{and $\tau_i$ is clockwise from $\tau_n$;}\\
          \text{$\bold{g}_{\tilde\gamma_1}+\bold{g}_{\tilde\gamma_2}$ \hspace{2.6cm}otherwise.}
      \end{cases}
  \end{equation}
  \end{itemize}
The definition is extended to any $\theta$-orbit by letting $F_{ab}^C=1$ and $\bold{g}_{ab}^C=\bold{e}_i$ if $[a,b]=\{ \tau_i,\tau_{2n-i} \} \in T$, and $F_{ab}^C=1$ and $\bold{g}_{ab}^C=\bold{0}$ if $(a,b)$ is a boundary edge of $\Poly_{2n+2}$. 
\end{definition}
\begin{remark}
    $(\Bar{b},c)$ in \ref{eq_g_vect_1} and $(\Bar{c},\Bar{d})$) in \ref{eq_g_vect_2} are either diagonals of $\Bar{T}$ or boundary edges, since $\tilde{\text{Res}}([a,\Bar{c}])=\{\tilde\gamma_1\}$ and $\tilde{\text{Res}}([a,d])=\{\tilde\gamma_1\}$ respectively. Remember that by convention $x_{(a,b)}=1$ if $(a,b)$ is a boundary edge, and so in that case $\bold{g}_{(a,b)}=\bold{0}$. 
\end{remark}
\begin{remark}
    We note that $F_{a\Bar{a}}^C$ (resp. $F_{ab}^C$ for $[a,b]$ pair of diagonals which cross $d$) are well-defined polynomial in $y_1,\dots,y_n$, since if $L_n$ crosses $(a,\ast)$ and $(c,\Bar{b})$ (resp. $(\Bar{b},\ast)$ and $(\Bar{d},\Bar{c})$), then it also crosses $(a,\Bar{b})$ and $(b,\ast)$ (resp. $(a,\Bar{c})$ and $(c,\ast)$).
\end{remark}

\begin{theorem}\label{theorem 2}
 Let $T$ be a $\theta$-invariant triangulation of $\Poly_{2n+2}$ with oriented diameter $d$, and let $\mathcal{A}=\mathcal{A}^C(T)$ be the cluster algebra of type $C_n$ with principal coefficients in $T$. 
    Let $[a,b]$ be an orbit of the action of $\theta$ on the diagonals of the polygon, and $x_{ab}$ the cluster variable of $\mathcal{A}$ which corresponds to $[a,b]$. Let $F_{ab}$ and $\bold{g}_{ab}$ denote the $F$-polynomial and the $\bold{g}$-vector of $x_{ab}$, respectively. 
    Then $F_{ab}=F_{ab}^C$ and $\bold{g}_{ab}=\bold{g}_{ab}^C$.
\end{theorem}
\begin{remark}
As observed for Theorem \ref{theorem1}, since for a diagonal $\gamma$ of $\Poly_{n+3}$, $F_\gamma$ and $\bold{g}_\gamma$ have an explicit description, for example in terms of perfect matchings of the snake graph associated with $\gamma$, Theorem \ref{theorem 2} also allows us to get the expansion of cluster variables of type $C_n$ in terms of the cluster variables of the initial seed. 
\end{remark}
\begin{example}\label{example_F_poly_type_C}
By Theorem \ref{theorem 2}, the $F$-polynomial of the cluster variable of type $C_3$ which corresponds to the $\theta$-orbit $[a,b]$ of $\Poly_8$ in Figure \ref{ex_f_poly_C} is   
\begin{center}
 $F_{ab}=F_{\tilde\gamma_1}F_{\tilde\gamma_2}-y_3F_{(a,c)}= (y_3y_2+y_3+1)(y_1+1)-y_3(y_2+1)=y_1y_2y_3+y_1y_3+y_1+1$,    
\end{center}

and the $\bold{g}$-vector is 
\begin{center}
$\displaystyle\bold{g}_{ab}=\bold{g}_{\tilde\gamma_1}+\bold{g}_{\tilde\gamma_2}+\bold{e}_2-\bold{e}_2=\bold{g}_{\tilde\gamma_1}+\bold{g}_{\tilde\gamma_2}=\displaystyle\begin{pmatrix}
    0 \\ 0 \\ -1
\end{pmatrix}+\displaystyle\begin{pmatrix}
    -1 \\ 0 \\ 1
\end{pmatrix}=\displaystyle\begin{pmatrix}
    -1 \\ 0 \\ 0
\end{pmatrix}$.  
\end{center}

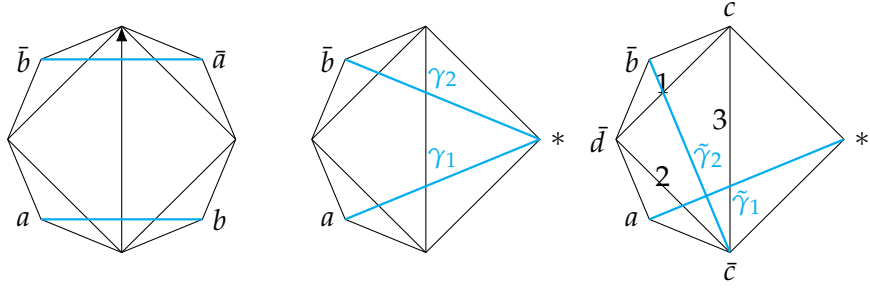
\begin{figure}[H]
                               \centering
                \begin{tikzpicture}[scale=0.5]
                    \draw (90:3cm) -- (135:3cm) -- (180:3cm) -- (225:3cm) -- (270:3cm) -- (315:3cm) -- (360:3cm) -- (45:3cm) -- cycle;
                    \draw (270:3cm) -- (180:3cm);;
                    \draw (90:3cm) -- (360:3cm);
                    \draw[-{Latex[length=2mm]}] (270:3cm) -- node[midway, above left,xshift=1mm] {} (90:3cm);
                    \draw (180:3cm) -- node[midway, above left,xshift=1mm] {} (90:3cm);
                    \draw (360:3cm) -- node[midway, above left,xshift=1mm] {} (270:3cm);
                    %\draw (90:3cm) -- node[midway, above right, xshift=-1mm] {} (315:3cm);
                    %\draw (90:3cm) -- node[midway, above right, xshift=-1mm] {} (360:3cm);
                    \draw[cyan, line width=0.3mm] (135:3cm) -- (45:3cm); 
                    \draw[cyan, line width=0.3mm] (225:3cm) -- (315:3cm);

                    \node at (135:3cm) [left] {$\Bar{b}$};
                    \node at (225:3cm) [left] {$a$};
                    \node at (45:3cm) [right] {$\Bar{a}$};
                    \node at (315:3cm) [right] {$b$};
                    
                    \begin{scope}[xshift=8cm]
                      \draw (90:3cm) -- (135:3cm) -- (180:3cm) -- (225:3cm) -- (270:3cm) -- (360:3cm) -- cycle;
 
                    \draw (180:3cm) -- (270:3cm);;
                    %\draw (45:3cm) -- (315:3cm);
                    \draw (90:3cm) -- node[midway, above left,xshift=1mm] {} (270:3cm);
                    \draw (180:3cm) -- node[midway, above left,xshift=1mm] {} (90:3cm);
                    %\draw (45:3cm) -- node[midway, above left,xshift=1mm] {} (270:3cm);
                    %\draw (90:3cm) -- node[midway, above right, xshift=-1mm] {} (315:3cm);
                    %\draw (90:3cm) -- node[midway, above right, xshift=-1mm] {} (360:3cm);

\draw[cyan, line width=0.3mm] (225:3cm) -- node[midway,above] {$\gamma_1$} (360:3cm); 
\draw[cyan, line width=0.3mm] (135:3cm) -- node[midway, above] {$\gamma_2$} (360:3cm);

\node at (135:3cm) [left] {$\Bar{b}$};
\node at (225:3cm) [left] {$a$};
\node at (360:3cm) [right] {$\ast$}; 
                    \end{scope}

                    \begin{scope}[xshift=16cm]
                      \draw (90:3cm) -- (135:3cm) -- (180:3cm) -- (225:3cm) -- (270:3cm) -- (360:3cm) -- cycle;
 
                    \draw (180:3cm) -- node[midway, above left,xshift=1mm] {2}(270:3cm);;
                    %\draw (45:3cm) -- (315:3cm);
                    \draw (90:3cm) -- node[midway, above left,xshift=1mm] {3} (270:3cm);
                    \draw (180:3cm) -- node[midway, left,xshift=1mm] {1} (90:3cm);
                    %\draw (45:3cm) -- node[midway, above left,xshift=1mm] {} (270:3cm);
                    %\draw (90:3cm) -- node[midway, above right, xshift=-1mm] {} (315:3cm);
                    %\draw (90:3cm) -- node[midway, above right, xshift=-1mm] {} (360:3cm);

\draw[cyan, line width=0.3mm] (225:3cm) -- node[midway,below] {$\tilde\gamma_1$} (360:3cm); 
\draw[cyan, line width=0.3mm] (135:3cm) -- node[midway, right,xshift=-1mm] {$\tilde\gamma_2$} (270:3cm);

\node at (135:3cm) [left] {$\Bar{b}$};
\node at (90:3cm) [above] {$c$};
\node at (270:3cm) [below] {$\Bar{c}$};
\node at (180:3cm) [left] {$\Bar{d}$};
\node at (225:3cm) [left] {$a$};
\node at (360:3cm) [right] {$\ast$}; 
                    \end{scope}
                \end{tikzpicture}
                               \caption{A $\theta$-orbit $[a,b]$ in a triangulated octagon, its restriction and its rotated restriction.}
                               \label{ex_f_poly_C}
                           \end{figure}

\end{example}

The proof of Theorem \ref{theorem 2} is similar to the one of Theorem \ref{theorem1}. For completeness we report it in Section \ref{sec_proof}.

\section{The categorification}\label{s_categ}

\subsection{Symmetric quivers and their representations}\label{symmetric_quivers}
In this section, first we report basic definitions of quiver, quiver algebra and their representations, in order to fix the notation. Standard references for these notions are for instance \cite{libro_blu,AR_book}. Then we recall some facts about symmetric quivers and their representations from \cite{DW} and \cite{boos2021degenerations}.

\vspace{0.25cm}
Let $k=\mathbb{C}$ be the field of complex numbers.

A $quiver$ is a finite oriented graph given by a quadruple $Q=(Q_0,Q_1,s,t)$, where $Q_0$ denotes the finite set of vertices of $Q$, $Q_1$ denotes the finite set of edges and $s,t:Q_1 \to Q_0$ are two functions that provide the orientation $\alpha:s(\alpha)\to t(\alpha)$ of arrows. The $path$ $algebra$ $kQ$ of $Q$ is defined to be the $k$-vector space with a basis given by the set of all paths in $Q$. The multiplication of two paths is defined by concatenation of paths. Let $R$ be the two-sided ideal generated by the arrows of $Q$. An ideal $I \subseteq kQ$ is said to be $admissible$ if there is an integer $m\geq 2$ such that $R^m \subseteq I \subseteq R^2$. Let $I$ be an admissible ideal. Then $(Q,I)$ is called a $bound$ $quiver$ and the quotient algebra $\mathcal{A}=kQ/I$ is called a $quiver$ $algebra$. 

A $representation$ of $Q$ (or $Q$-representation) is a pair $(V,f)$, where $V$ is a $Q_0$-graded vector space, and $f$ is a collection of maps $f_\alpha$, $\alpha \in Q_1$, such that $f_\alpha:V_{s(\alpha)}\to V_{t(\alpha)}$ is a linear map. A representation of $(Q,I)$ is a $Q$-representation satisfying the relations from $I$.

\begin{definition}
    A $symmetric$ $quiver$ is a pair $(Q,\sigma)$, where $Q$ is a finite quiver and $\sigma$ is an involution of $Q_0$ and of $Q_1$ which reverses the orientation of arrows. 
\end{definition}

\begin{example}
    Let $Q = 1 \xrightarrow[]{\alpha} 2 \xrightarrow[]{\beta} 3$ and $Q'= 1 \xrightarrow[]{\alpha} 2 \xleftarrow[]{\beta} 3$ be two quivers of type $A_3$. Then $Q$ is symmetric, with the involution $\sigma$ given by $\sigma(1)=3$, $\sigma(2)=2$ and $\sigma(\alpha)=\beta$, while $Q'$ is not symmetric, i.e., it cannot be endowed with the structure of a symmetric quiver.
\end{example}

\begin{definition}
    Let $(Q,\sigma)$ be a symmetric quiver. Let $I \subset kQ$ be an admissible ideal such that $\sigma(I)=I$. $(Q,I,\sigma)$ is called a $bound$ $symmetric$ $quiver$ and the pair $(\mathcal{A}=kQ/I,\sigma)$ is called a $symmetric$ $quiver$ $algebra$.
\end{definition}

\begin{definition}
    A $symmetric$ $representation$ of a bound symmetric quiver $(Q,I,\sigma)$ is a triple 
    $(V,f,$ $ \langle \cdot , \cdot \rangle)$, where $(V,f)$ is a representation of $(Q,I)$, $\langle \cdot , \cdot \rangle$ is a nondegenerate symmetric or skew-symmetric scalar product on $V$ such that its restriction to $V_i \times V_j$ is 0 if $j \neq \sigma(i)$, and $\langle f_\alpha(v) , w \rangle + \langle v , f_{\sigma(\alpha)}(w) \rangle=0$, for every $\alpha : i \to j \in Q_1$, $v \in V_i$, $w \in V_{\sigma(j)}$. If $\langle \cdot , \cdot \rangle$ is symmetric (resp. skew-symmetric), $(V,f, \langle \cdot , \cdot \rangle)$ is called $orthogonal$ (resp $symplectic$).
\end{definition}
\begin{remark}
If $\bold{d}=(\mathrm{dim}(V_i))$ is the dimension vector of a symmetric representation $(V,f, \langle \cdot , \cdot \rangle)$ of a bound symmetric quiver $(Q,I,\sigma)$, then $d_i=d_{\sigma(i)}$. If the dimension vector $\bold{d}$ of a $(Q,I)$-representation has this property, we say that it is \emph{symmetric}.   
\end{remark}

\begin{definition}
    If $(V,f, \langle \cdot , \cdot \rangle)$ and $(V',f', \langle \cdot , \cdot \rangle')$ are symmetric representations of a bound symmetric quiver $Q$, then their direct sum is given by $(V\oplus V',f \oplus f', \langle \cdot , \cdot \rangle + \langle \cdot , \cdot \rangle')$. A symmetric representation is called $indecomposable$ if it is nontrivial and it is not isomorphic to the direct sum of two nontrivial symmetric representations.
\end{definition}

\begin{definition}
    Let $L=(V,f)$ be a representation of a bound symmetric quiver $Q$. The $twisted$ $dual$ of $L$ is the $\mathcal{A}$-representation $\nabla L = (\nabla V, \nabla f)$, where $(\nabla V)_i=V_{\sigma(i)}^\ast$ and $(\nabla f)_\alpha = - f_{\sigma (\alpha)}^\ast$ ($\ast$ denotes the linear dual).
\end{definition}

\begin{remark}
   If $L$ is symmetric, the scalar product $\langle \cdot , \cdot \rangle$ induces an isomorphism from $V$ to $\nabla V$. 
\end{remark}

\begin{lemma}[Lemma 2.10, \cite{boos2021degenerations}]\label{ind_symm}
 Let $M$ be an indecomposable symmetric representation of a bound symmetric quiver $Q$. Then, one and only one of the following three cases can occur: 
 \begin{itemize}
     \item [(I)] $M$ is indecomposable as a $Q$-representation; in this case, $M$ is called of type (I), for “indecomposable”;
     \item [(S)] there exists an indecomposable $Q$-representation $L$ such that $M=L\oplus \nabla L$ and $L\ncong \nabla L$; in this case, $M$ is called of type (S), for “split”;
     \item [(R)] there exists an indecomposable $Q$-representation $L$ such that $M=L\oplus \nabla L$ and $L\cong \nabla L$; in this case, $M$ is called of type (R) for “ramified”.
 \end{itemize}
\end{lemma}

\subsection{$\rho$-orbits as orthogonal and symplectic representations}\label{quiver_to_triang}

Let $\Bar{T}$ be a triangulation of $\Poly_{n+3}$, and let $Q(\Bar{T})$ be the quiver associated to $\Bar{T}$ as in \cite{FST,LF}, so that there is an arrow from the vertex $j$ to the vertex $i$ if and only if $\tau_i$ and $\tau_j$ are sides of a triangle of $\Bar{T}$, and $\tau_i$ is counterclockwise from $\tau_j$, and the relations are given by all paths $i\to j \to k$ such that there exists an arrow $k\to i$.  Then $Q(\Bar{T})$ is a cluster-tilted bound quiver of type $A_{n}$ (see \cite{S_Quiver_rep}, 3.4.1). Since $\Bar{T}$ is a triangulation of the polygon, any other diagonal $\gamma$ which is not already in $\Bar{T}$ will cut through a certain number of diagonals in $\Bar{T}$; in fact, any such diagonal $\gamma$ is uniquely determined by the set of diagonals in $\Bar{T}$ that $\gamma$ crosses. To such a diagonal $\gamma$, it is associated a representation $L=(V,f)$ of $Q(\Bar{T})$ defined as follows:
\begin{center}
    $V_i=\begin{cases}
        \text{$k$ if $\gamma$ crosses the diagonal $i;$}\\
        \text{0 otherwise;}
    \end{cases}$
\end{center}
and $f_\alpha = 1$ whenever $V_{s(\alpha)}=V_{t(\alpha)}=k$, and $f_\alpha=0$ otherwise.

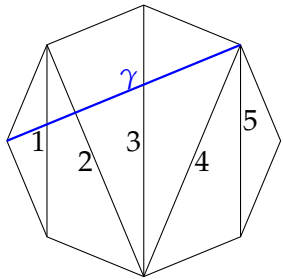
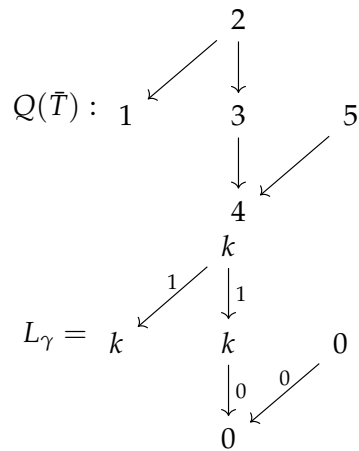
\begin{example}\label{ex_quiver}
\end{example}

     \begin{figure}[H]
        \centering
\begin{tikzpicture}[scale=0.6]
     \draw (90:3cm) -- (135:3cm) -- (180:3cm) -- (225:3cm) -- (270:3cm) -- (315:3cm) -- (360:3cm) -- (45:3cm) -- cycle;
                    \draw (90:3cm) -- node[midway, left,xshift=1mm] {3} (270:3cm);
                    \draw (135:3cm) -- node[midway, left,xshift=1.2mm] {1} (225:3cm);
                    \draw (135:3cm) -- node[midway, left,xshift=1mm] {2} (270:3cm);
                    \draw (45:3cm) -- node[midway, right, xshift=-1mm] {4} (270:3cm);
                    \draw (45:3cm) -- node[midway, above right, xshift=-1mm] {5} (315:3cm);
                    \draw[blue, line width=0.3mm] (180:3cm) -- node[midway, above right, xshift=-2mm, yshift=-1mm] {$\gamma$} (45:3cm);
                 
%\node at (10,2) {$Q(T): 1 \leftarrow 2 \rightarrow 3 \rightarrow 4 \leftarrow 5$};
\node at (10,2) {$Q(\Bar{T}): \begin{tikzcd}
& 2 \arrow[dl] \arrow[d] \\
1  & 3 \arrow[d] & 5 \arrow[dl]  \\
&  4 
\end{tikzcd}$};
\node at (10,-3) {$L_\gamma=\begin{tikzcd}
& k \arrow[dl, "1" above] \arrow[d,"1"] \\
k  & k \arrow[d,"0"] & 0 \arrow[dl,"0" above]  \\
&  0 
\end{tikzcd}$};

\end{tikzpicture}
        %\caption{An example of cluster for a cluster algebra of type $B_3$.}
        \label{fig:enter-label}
    \end{figure}
Sometimes we will use indices of vertices with a nonzero dimensional vector space to indicate representations. For instance, for $L_\gamma$ of the previous example the shorthand is $\tabbedCenterstack{2\\13}$.

The map $\gamma \mapsto L_\gamma$ is a bijection from the set of diagonals that are not in $\Bar{T}$ and the set of isoclasses of indecomposable representations of $Q(\Bar{T})$. 

\begin{remark}\label{orbits_rmk}
 Let $d$ be a diameter of $\Poly_{2n+2}$. Let $\rho$ denote the reflection of the polygon along $d$. It induces an action on the diagonals of the polygon.
If $T'$ is a $\rho$-invariant triangulation of $\Poly_{2n+2}$, then $(Q(T'),\sigma_\rho)$ is a cluster-tilted bound symmetric quiver of type $A_{2n-1}$, with involution $\sigma_\rho$ induced by $\rho$.

\begin{example}
Let $\rho$ be the reflection of the octagon along the diameter $d$ in Figure \ref{ex_quiver_to_rho_inv_triang}. Let $\sigma_\rho$ be the involution of $Q(T')$ defined by $\sigma_\rho(1)=\rho(1)=5$, $\sigma_\rho(2)=\rho(2)=4$, $\sigma_\rho(3)=\rho(3)=3$, and $\sigma_\rho(\alpha)=\delta$, $\sigma_\rho(\beta)=\gamma$. Then $(Q(T'),\sigma_\rho)$ is a symmetric quiver of type $A_5$.
    \begin{figure}[H]
        \centering
        \begin{tikzpicture}[scale=0.6]
            \node at (-4,0) {$T'=$};
     \draw (90:3cm) -- (135:3cm) -- (180:3cm) -- (225:3cm) -- (270:3cm) -- (315:3cm) -- (360:3cm) -- (45:3cm) -- cycle;
                     \draw (90:3cm) -- node[midway, above, xshift=-1mm, yshift=-1mm] {$1$} (180:3cm);
                    \draw (90:3cm) -- node[midway, above left,xshift=1.2mm] {$2$} (225:3cm);
                    \draw (90:3cm) -- node[midway] {$3=d$} (270:3cm);
                    \draw (90:3cm) -- node[midway, above right, xshift=-2mm] {$4$} (315:3cm);
                    \draw (90:3cm) -- node[midway, above right, xshift=-2mm] {$5$} (360:3cm);
   \node at (9,0) {$Q(T'): 1 \xleftarrow[]{\alpha} 2 \xleftarrow[]{\beta} 3 \xleftarrow[]{\gamma} 4 \xleftarrow[]{\delta} 5$};         
        \end{tikzpicture}
        \caption{A $\rho$-invariant triangulation of $\Poly_8$ and the associated quiver.}
        \label{ex_quiver_to_rho_inv_triang}
    \end{figure}
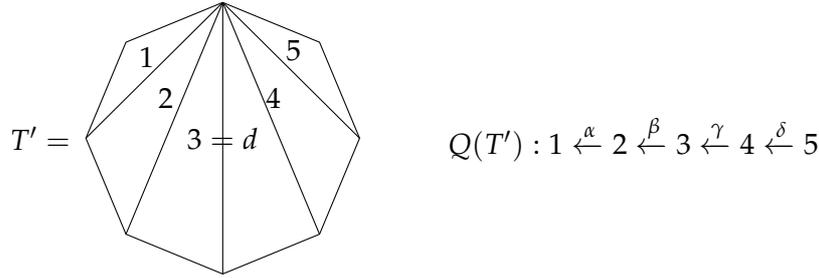

\end{example}

Moreover, if $[a,b]^\rho=\{\alpha_1, \alpha_2\}$ is a $\rho$-orbit and $\alpha_1$ corresponds to the indecomposable representation of $Q(T')$ $L_{\alpha_1}$, then $\alpha_2$ corresponds to $L_{\alpha_2}=\nabla L_{\alpha_1}$. In fact, if we denote by $\bold{d}_{\alpha_i}$ the vector of indices of diagonals of $T'$ crossed by $\alpha_i$, i.e. the dimension vector of $L_{\alpha_i}$, we have that both $\bold{d}_{\alpha_1}$ and $\bold{d}_{\alpha_2}$ are not symmetric, while $\bold{d}_{\alpha_1}+\bold{d}_{\alpha_2}$ is. It follows from Lemma \ref{ind_symm} that $L_{\alpha_1} \oplus L_{\alpha_2}$ is symmetric indecomposable of type S, so $L_{\alpha_2}=\nabla L_{\alpha_1}$.  

On the other hand, if $[a,b]^\rho=\{\alpha\}$, then $\alpha$ corresponds to the $\nabla$-invariant indecomposable representation of $Q(T')$ $L_\alpha$, since $\bold{d}_\alpha$ is symmetric.

\vspace{0.3cm}

Let $T'=\{ \tau_1, \dots, \tau_{2n-1} \}$ be a $\rho$-invariant triangulation of $\Poly_{2n+2}$. Then it has $n-1$ $\rho$-invariant pairs of diagonals not orthogonal to $d$ and exactly one $\rho$-invariant diagonal $\tau_n$. We have two cases to consider.
\begin{itemize}
    \item [$\tau_n=d$] In this case $Q(T')$ has a fixed vertex $n$ and no fixed arrows. Therefore, every $\rho$-invariant diagonal $\alpha$ which is not in $T'$ crosses $\tau_n$. So $L_\alpha$ is orthogonal indecomposable of type I, while $L_\alpha \oplus L_\alpha$ is symplectic indecomposable of type R, since in the latter case the nonzero vector space at vertex $n$ of the quiver must be a symplectic space, so it must have dimension 2. 
    \begin{example}
    \end{example}
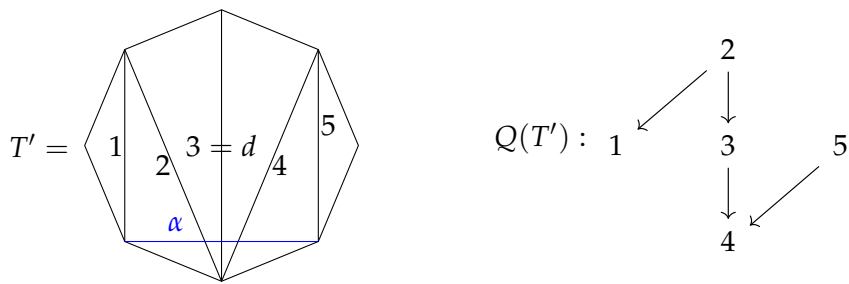
    \begin{figure}[H]
        \centering
\begin{tikzpicture}[scale=0.6]
     \draw (90:3cm) -- (135:3cm) -- (180:3cm) -- (225:3cm) -- (270:3cm) -- (315:3cm) -- (360:3cm) -- (45:3cm) -- cycle;
                    \draw (90:3cm) -- node[midway] {$3=d$} (270:3cm);
                    \draw (135:3cm) -- node[midway, left,xshift=1.2mm] {1} (225:3cm);
                    \draw (135:3cm) -- node[midway, left,xshift=1mm] {2} (270:3cm);
                    \draw (45:3cm) -- node[midway, right, xshift=-1mm] {4} (270:3cm);
                    \draw (45:3cm) -- node[midway, above right, xshift=-1mm] {5} (315:3cm);
                    \draw[blue] (225:3cm) -- node[left, above, xshift=-6mm] {$\alpha$} (315:3cm);
                  
 \node at (-4,0) {$T'=$};                 
%\node at (10,2) {$Q(T): 1 \leftarrow 2 \rightarrow 3 \rightarrow 4 \leftarrow 5$};
\node at (10,0) {$Q(T'): \begin{tikzcd}
& 2 \arrow[dl] \arrow[d] \\
1  & 3 \arrow[d] & 5 \arrow[dl]  \\
&  4 
\end{tikzcd}$};
%\node at (10,-3) {$L_\gamma=\begin{tikzcd}
%& k \arrow[dl, "1" above] \arrow[d,"1"] \\
%k  & k \arrow[d,"0"] & 0 \arrow[dl,"0" above]  \\
%&  0 
%\end{tikzcd}$};

\end{tikzpicture}
    \end{figure}
    \item [$\tau_n \neq d$] In this case $Q(T')$ has a fixed vertex $n$ and a fixed arrow $\beta: i \to j$. Therefore, every $\rho$-invariant diagonal $\alpha$ which is not in $T'$ crosses $i$ and $j$, while it cannot cross $\tau_n$. Let $\{v\}$ be a basis of the 1-dimensional vector space of $L_\alpha$ at vertex $i$ and let $\{w\}$ be a basis of the 1-dimensional vector space of $L_\alpha$ at vertex $j$. If $(L_\alpha,\langle \cdot , \cdot \rangle)$ is a symmetric representation of $Q(T')$, then by definition
    \begin{equation}
        \langle w , v \rangle = \langle f_\beta(v) , v \rangle = - \langle v , f_{\sigma_\rho(\beta)}(v) \rangle = -\langle v , f_\beta(v) \rangle = -\langle v , w \rangle.
    \end{equation}
    Since $\langle \cdot , \cdot \rangle$ is a non-degenerate scalar product, it must be skew-symmetric. It follows from Lemma \ref{ind_symm} that $L_\alpha$ is symplectic indecomposable of type I, while $L_\alpha \oplus L_\alpha$ is orthogonal indecomposable of type R.

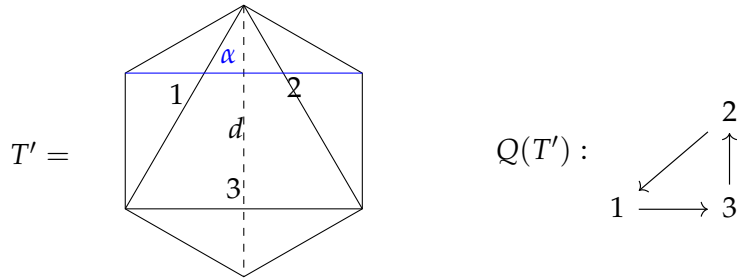
    \begin{example}
      \end{example}
      \begin{figure}[H]
    \centering
   \begin{tikzpicture}[scale=0.6,rotate=30]
    \draw (0:3cm) -- (60:3cm) -- (120:3cm) -- (180:3cm) -- (240:3cm) -- (300:3cm) -- cycle;
                    \draw (60:3cm) -- node[midway, above, xshift=-1mm, yshift=-1mm] {1} (180:3cm);
                    \draw (300:3cm) -- node[midway, above left,xshift=1.2mm] {2} (60:3cm);
                    \draw (300:3cm) -- node[midway, above left,xshift=1mm] {3} (180:3cm);
                   \draw[dashed] (60:3cm) -- node[midway, above, xshift=-1mm, yshift=-1mm] {$d$} (240:3cm);
                    \draw[blue] (120:3cm) -- node[left, above, xshift=-2mm] {$\alpha$} (0:3cm);
                    \node at (-4,2) {$T'=$};
                   \node at (7,-4.5) {$Q(T'): \begin{tikzcd}
& 2 \arrow[dl]  \\
1 \arrow[r] & 3 \arrow[u] 
\end{tikzcd}$};
  \end{tikzpicture}
  \end{figure}
  
\end{itemize}

\end{remark}

\vspace{0.8cm}

Let $T$ be a $\theta$-invariant triangulation of $\Poly_{2n+2}$ with oriented diameter $d$. Then $Q(T)$ is not symmetric. 

\begin{example}
Let $T$ be $\theta$-invariant triangulation of the octagon in Figure \ref{ex_quiver_to_theta_inv_triang}. Then the quiver $Q(T)$ is not symmetric.
    \begin{figure}[H]
        \centering
        \begin{tikzpicture}[scale=0.6]
            \node at (-4,0) {$T=$};
     \draw (90:3cm) -- (135:3cm) -- (180:3cm) -- (225:3cm) -- (270:3cm) -- (315:3cm) -- (360:3cm) -- (45:3cm) -- cycle;
                     \draw (90:3cm) -- node[midway, above, xshift=-1mm, yshift=-1mm] {1} (180:3cm);
                    \draw (90:3cm) -- node[midway, above left,xshift=1.2mm] {2} (225:3cm);
                    \draw (90:3cm) -- node[midway, above left,xshift=1mm] {3} (270:3cm);
                    \draw (270:3cm) -- node[midway, above right, xshift=-2mm] {4} (45:3cm);
                    \draw (270:3cm) -- node[midway, above right, xshift=-1mm] {5} (360:3cm);
   \node at (9,0) {$Q(T): 1 \xleftarrow[]{} 2 \xleftarrow[]{} 3 \xrightarrow[]{} 4 \xrightarrow[]{} 5$};         
        \end{tikzpicture}
        \caption{A $\theta$-invariant triangulation of $\Poly_8$ and the associated quiver.}
        \label{ex_quiver_to_theta_inv_triang}
    \end{figure}
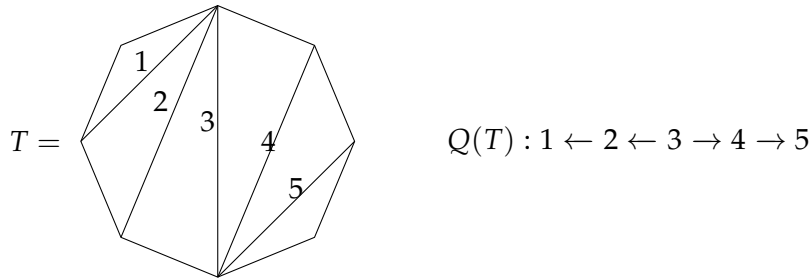

\end{example}

In order to get a symmetric quiver, we define an involution on the polygon that we call $F_d$. 
\begin{definition}
    $F_d$ is the operation on $\Poly_{2n+2}$ which consists of the following three steps in order:
    \begin{itemize}
        \item [1)] cut the polygon along $d$;
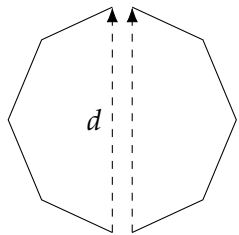
         \begin{figure}[H]
        \centering
\begin{tikzpicture}[scale=0.5]
     \draw (95:3cm) -- (135:3cm) -- (180:3cm) -- (225:3cm) -- (265:3cm) ; 
     \draw (275:3cm)-- (315:3cm) -- (360:3cm) -- (45:3cm) -- (85:3cm) ;
         \draw[dashed,-{Latex[length=2mm]}](265:3cm) -- node[midway, left] {$d$} (95:3cm); 
         \draw[dashed,-{Latex[length=2mm]}](275:3cm) -- node[midway, left, xshift=0.4cm,yshift=-0.5cm] {} (85:3cm); 

                   \end{tikzpicture}
                   \end{figure}
        \item [2)] reflect the right part with respect to the axis of symmetry of $d$;
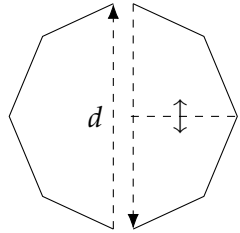
              \begin{figure}[H]
        \centering
\begin{tikzpicture}[scale=0.5]
     \draw (95:3cm) -- (135:3cm) -- (180:3cm) -- (225:3cm) -- (265:3cm) ; 
     \draw (275:3cm)-- (315:3cm) -- (360:3cm) -- (45:3cm) -- (85:3cm) ;
         \draw[dashed,-{Latex[length=2mm]}](265:3cm) -- node[midway, left] {$d$} (95:3cm); 
         \draw[dashed,-{Latex[length=2mm]}](85:3cm) -- node[midway, left, xshift=0.4cm,yshift=-0.5cm] {} (275:3cm); 
         \draw[dashed](0.2,0) -- node[midway, left, xshift=0.4cm,yshift=-0.5cm] {} (360:3cm);
         \node at (1.5,0) {$\big\updownarrow$};

                   \end{tikzpicture}
                   \end{figure}
        \item [3)] glue  again the right part along $d$.
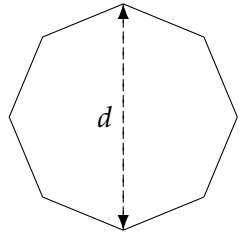
         \begin{figure}[H]
        \centering
\begin{tikzpicture}[scale=0.5]
     \draw (90:3cm) -- (135:3cm) -- (180:3cm) -- (225:3cm) -- (270:3cm) ; 
     \draw (270:3cm)-- (315:3cm) -- (360:3cm) -- (45:3cm) -- (90:3cm) ;
         \draw[dashed,-{Latex[length=2mm]}](270:3cm) -- node[midway, left] {$d$} (90:3cm); 
         \draw[dashed,-{Latex[length=2mm]}](90:3cm) -- node[midway, left, xshift=0.4cm,yshift=-0.5cm] {} (270:3cm); 

                   \end{tikzpicture}
                   \end{figure}
    \end{itemize}
   
\end{definition}
\begin{remark}
    $F_d$ induces an action on isotopy classes of diagonals of the polygon.
\end{remark}
 
\begin{lemma}\label{lemma_f_d}
Under the bijection $F_d$, $\theta$-orbits correspond to $\rho$-orbits. In particular, diameters correspond to $\rho$-invariant diagonals, while pairs of centrally symmetric diagonals correspond to $\rho$-invariant pairs of diagonals which are not orthogonal to $d$. 
\end{lemma}
\begin{proof}
Let $[a,b]$ be a $\theta$-orbit. We have three cases to consider:
\begin{itemize}
    \item [i)] $(a,b)$ is a diameter (illustrated in Figure \ref{fig_i});
    \item [ii)] $[a,b]$ is a pair of centrally symmetric diagonals which cross $d$ (illustrated in Figure \ref{fig_ii});
    \item [iii)] $[a,b]$ is a pair of centrally symmetric diagonals which do not cross $d$ (illustrated in Figure \ref{fig_iii}).
\end{itemize}

 \begin{figure}[H]
        \centering
\begin{tikzpicture}[scale=0.5]
     \draw (90:3cm) -- (135:3cm) -- (180:3cm) -- (225:3cm) -- (270:3cm) -- (315:3cm) -- (360:3cm) -- (45:3cm) -- cycle;
       \draw (135:3cm) -- (315:3cm); 
         \draw[dashed,-{Latex[length=2mm]}](270:3cm) -- node[midway, left, xshift=0.4cm,yshift=-0.5cm] {$d$} (90:3cm);
         \node at (135:3) [left] {$a$};
         \node at (315:3) [right] {$b$};
 \node at (4,0) {$\xrightarrow[]{F_d}$};                  
\begin{scope}[xshift=8cm]
    \draw (90:3cm) -- (135:3cm) -- (180:3cm) -- (225:3cm) -- (270:3cm) -- (315:3cm) -- (360:3cm) -- (45:3cm) -- cycle;
     \draw (0,0) -- (135:3cm); 
      \draw (0,0) -- (45:3cm); 
        \node at (45:3) [right] {$b$}; 
       \draw[dashed,-{Latex[length=2mm]}](270:3cm) -- node[midway, left, xshift=0.4cm,yshift=-0.5cm] {$d$} (90:3cm);
        \node at (135:3) [left] {$a$};
\end{scope}
\node at (12,0) {$\cong$};
\begin{scope}[xshift=16cm]
    \draw (90:3cm) -- (135:3cm) -- (180:3cm) -- (225:3cm) -- (270:3cm) -- (315:3cm) -- (360:3cm) -- (45:3cm) -- cycle;
     \draw (45:3cm) -- (135:3cm); 
        \node at (45:3) [right] {$b$}; 
       \draw[dashed,-{Latex[length=2mm]}](270:3cm) -- node[midway, left, xshift=0.4cm,yshift=-0.5cm] {$d$} (90:3cm) ;
        \node at (135:3) [left] {$a$};
\end{scope}
                   
\end{tikzpicture}
        \caption{The action of $F_d$ on the diameter $(a,b)$.}
        \label{fig_i}
    \end{figure}
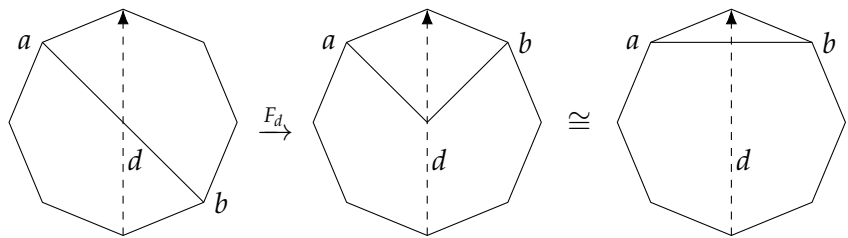  
   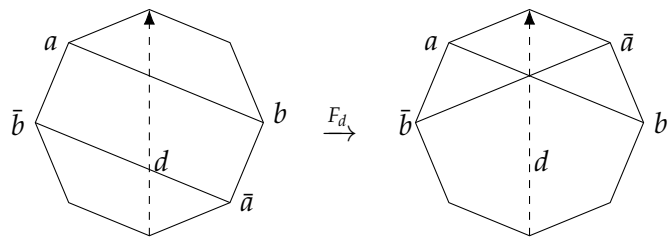
\begin{figure}[H]
        \centering
\begin{tikzpicture}[scale=0.5]
     \draw (90:3cm) -- (135:3cm) -- (180:3cm) -- (225:3cm) -- (270:3cm) -- (315:3cm) -- (360:3cm) -- (45:3cm) -- cycle;
       \draw (135:3cm) -- (360:3cm); 
       \draw (180:3cm) -- (315:3cm); 
         \draw[dashed,-{Latex[length=2mm]}](270:3cm) -- node[midway, left, xshift=0.4cm,yshift=-0.5cm] {$d$} (90:3cm);
         \node at (135:3) [left] {$a$};
         \node at (365:3) [right] {$b$};
         \node at (315:3) [right] {$\Bar{a}$};
         \node at (180:3) [left] {$\Bar{b}$};
 \node at (5,0) {$\xrightarrow[]{F_d}$};                  

\begin{scope}[xshift=10cm]
    \draw (90:3cm) -- (135:3cm) -- (180:3cm) -- (225:3cm) -- (270:3cm) -- (315:3cm) -- (360:3cm) -- (45:3cm) -- cycle;
     \draw (360:3cm) -- (135:3cm); 
     \draw (45:3cm) -- (180:3cm); 
        \node at (360:3) [right] {$b$}; 
       \draw[dashed,-{Latex[length=2mm]}](270:3cm) -- node[midway, left, xshift=0.4cm,yshift=-0.5cm] {$d$} (90:3cm) ;
        \node at (135:3) [left] {$a$};
        \node at (180:3) [left,xshift=0.1cm] {$\Bar{b}$};
         \node at (45:3) [right] {$\Bar{a}$};
\end{scope}
                   
\end{tikzpicture}
        \caption{The action of $F_d$ on the $\theta$-orbit $[a,b]$ whose diagonals cross $d$.}
        \label{fig_ii}
    \end{figure} 

     \begin{figure}[H]
        \centering
\begin{tikzpicture}[scale=0.5]
     \draw (90:3cm) -- (135:3cm) -- (180:3cm) -- (225:3cm) -- (270:3cm) -- (315:3cm) -- (360:3cm) -- (45:3cm) -- cycle;
       \draw (135:3cm) -- (270:3cm); 
       \draw (90:3cm) -- (315:3cm); 
         \draw[dashed,-{Latex[length=2mm]}](270:3cm) -- node[midway, left, xshift=0.4cm,yshift=-0.5cm] {$d$} (90:3cm);
         \node at (135:3) [left] {$a$};
         \node at (270:3) [below] {$b$};
          \node at (315:3) [right] {$\Bar{a}$};
         \node at (90:3) [above] {$\Bar{b}$};
 \node at (5,0) {$\xrightarrow[]{F_d}$};                  
\begin{scope}[xshift=10cm]
   \draw (90:3cm) -- (135:3cm) -- (180:3cm) -- (225:3cm) -- (270:3cm) -- (315:3cm) -- (360:3cm) -- (45:3cm) -- cycle;
       \draw (135:3cm) -- (270:3cm); 
       \draw (270:3cm) -- (45:3cm); 
         \draw[dashed,-{Latex[length=2mm]}](270:3cm) -- node[midway, left, xshift=0.4cm,yshift=-0.5cm] {$d$} (90:3cm);
         \node at (135:3) [left] {$a$};
         \node at (270:3) [below] {$b$};
          \node at (45:3) [right] {$\Bar{a}$};
         \node at (90:3) [above] {$\Bar{b}$};
\end{scope}

\end{tikzpicture}
        \caption{The action of $F_d$ on the $\theta$-orbit $[a,b]$ whose diagonals do not cross $d$.}
        \label{fig_iii}
    \end{figure}
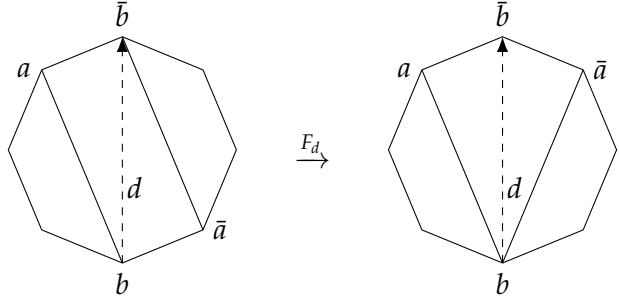 
\end{proof}

\begin{remark}
Let $T'$ be the element in the isotopy class of $F_d(T)$ which is also a triangulation. It follows from Lemma \ref{lemma_f_d} that $T'$ is a $\rho$-invariant triangulation of $\Poly_{2n+2}$ which contains the diameter $d$. Then $Q(T')$ is a cluster-tilted bound symmetric quiver of type $A_{2n-1}$ with a fixed vertex and no fixed arrows (cf. Remark \ref{orbits_rmk}).
\end{remark}

Now, let $\mathcal{A}=\mathcal{A}^B (T)$ be the cluster algebra of type $B$ with principal coefficients in $T$ defined in Section \ref{sec:ca-type-C}. Let $[a,b]$ be a $\theta$-orbit and let $x_{ab}$ be the cluster variable which corresponds to $[a,b]$. If $F_d([a,b])=\{\alpha\}$ consists of only one $\rho$-invariant diagonal, then $x_{ab}$ corresponds to the orthogonal indecomposable $Q(T')$-representation $L_\alpha$ of type I (cf. Remark \ref{orbits_rmk}). Otherwise, $F_d([a,b])=\{\alpha_1,\alpha_2\}$. In this case, $x_{ab}$ corresponds to $L_{\alpha_1} \oplus L_{\alpha_2}$ which is an orthogonal indecomposable $Q(T')$-representation of type S by Remark \ref{orbits_rmk}. 

Moreover, the restriction on $\theta$-orbits corresponds to an operation on orthogonal indecomposable $Q(T')$-representations defined in the following way: 
\begin{definition}
    Let $M=(V,f,\langle \cdot, \cdot \rangle)$ be an orthogonal indecomposable $Q(T')$- representation. Then the \emph{restriction} of $M$ is $\text{Res}(M)=(\text{Res}(V),\text{Res}(f))$, where $\text{Res}(V)_i=V_i$ if $i \leq n$, $\text{Res}(V)_i=0$ otherwise; and $\text{Res}(f)_\alpha=f_\alpha$ if $\alpha : i \to j$, with $i,j \leq n$, $\text{Res}(f)_\alpha=0$ otherwise. In other words, if $[a,b]$ is the $\theta$-orbit which corresponds to $M$, and $\text{Res}([a,b])=\{\gamma_1,\gamma_2\}$ (resp. $\text{Res}([a,b])=\{\gamma\}$), then $\text{Res}(M)=L_{\gamma_1}\oplus L_{\gamma_2}$ (resp. $\text{Res}(M)=L_{\gamma}$).
\end{definition}

\begin{remark}
    Note that $\text{Res}(M)$ is no longer orthogonal. Moreover, $\text{Res}(M)$ is a representation of the quiver associated to the triangulation of $\Poly_{n+3}$ obtained from $T'$ by identifying the vertices which lie on the right of $d$, i.e. $\Bar{T}=\text{Res}(T')=\text{Res}(T)$ (the part of $T$ on the left of $d$ is equal to the one of $T'$ on the left of $d$).
\end{remark}
On the other hand, let $\mathcal{A}=\mathcal{A}^C(T)$ be the cluster algebra of type $C$ with principal coefficients in $T$ defined in Section \ref{sec:ca-type-C}. Let $[a,b]$ be a $\theta$-orbit and let $x_{ab}$ be the cluster variable which corresponds to $[a,b]$. If $F_d([a,b])=\{\alpha\}$ consists of only one $\rho$-invariant diagonal, then $x_{ab}$ corresponds to the symplectic indecomposable $Q(T')$-representation $L_\alpha \oplus L_\alpha$ of type R (cf. Remark \ref{orbits_rmk}). Otherwise, $F_d([a,b])=\{\alpha_1,\alpha_2\}$. As before, $x_{ab}$ corresponds to the symplectic indecomposable $Q(T')$-representation $L_{\alpha_1} \oplus L_{\alpha_2}= L_{\alpha_1} \oplus \nabla L_{\alpha_1}$ of type S. 

Moreover, the rotated restriction on $\theta$-orbits corresponds to the operation on symplectic $Q(T')$-representations defined in the following way: 

\begin{definition}
Let $M$ be an indecomposable symplectic representation of $Q(T')$, and let $[a,b]$ be the $\theta$-orbit that corresponds to $M$. If $\Tilde{\text{Res}}([a,b])=\{\tilde\gamma_1,\tilde\gamma_2\}$ (resp. $\Tilde{\text{Res}}([a,b])=\{\tilde\gamma\}$), then $\Tilde{\text{Res}}(M)=L_{\tilde\gamma_1}\oplus L_{\tilde\gamma_2}$ (resp. $\Tilde{\text{Res}}(M)=L_{\tilde\gamma}$).
\end{definition}

\begin{comment}
    Moreover, the rotated restriction on $\theta$-orbits corresponds to an operation on symplectic $Q(T')$-representations defined in the following way: 
\begin{definition}
    Let $M=(V,f,\langle \cdot, \cdot \rangle)$ be a symplectic $Q(T')$-representation. Then the $rotated$ $restriction$ $of$ $M$ is $\Tilde{\text{Res}}(M)=(\Tilde{\text{Res}}(V),\Tilde{\text{Res}}(f))$, where $\Tilde{\text{Res}}(V)=V_i$ if $i < n$, $\Tilde{\text{Res}}(V)_n=K^{\text{dim}(V_n)/2}$, $\text{Res}(V)_i=0$ otherwise; and $\Tilde{\text{Res}}(f)_\alpha=f_\alpha$ if $\alpha : i \to j$, with $i,j < n$, $\Tilde{\text{Res}}(f)_\alpha=(f_\alpha)_1$, if $f_\alpha=((f_\alpha)_1,(f_\alpha)_2) : V_i \to V_n = k^{\text{dim}(V_n)/2}\oplus k^{\text{dim}(V_n)/2}$, $i<n$, $\Tilde{\text{Res}}(f)_\alpha=0$ otherwise.
\end{definition}
\end{comment}
\begin{remark}
Note that $\Tilde{\text{Res}}(M)$ is no longer symplectic. Moreover, as for $\text{Res}(M)$, $\Tilde{\text{Res}}(M)$ is a representation of the quiver associated to the triangulation $\Bar{T}=\text{Res}(T')=\text{Res}(T)$ of $\Poly_{n+3}$.
\end{remark}
%\newpage
\begin{example}
\end{example}

   \begin{figure}[H]
        \centering
\begin{tikzpicture}[scale=0.6]
\node at (-4,7) {$T=$};
\begin{scope}[yshift=7cm]
     \draw (90:3cm) -- (135:3cm) -- (180:3cm) -- (225:3cm) -- (270:3cm) -- (315:3cm) -- (360:3cm) -- (45:3cm) -- cycle;
                    \draw (90:3cm) -- node[midway, above, xshift=-1mm, yshift=-1mm] {} (180:3cm);
                    \draw (90:3cm) -- node[midway, above left,xshift=1.2mm] {} (225:3cm);
                    \draw (90:3cm) -- node[midway, above left,xshift=1mm] {} (270:3cm);
                    \draw (270:3cm) -- node[midway, above right, xshift=-1mm] {} (45:3cm);
                    \draw (270:3cm) -- node[midway, above right, xshift=-1mm] {} (360:3cm);
                  
\end{scope}
\node at (6,7) {$T'=$};
\begin{scope}[xshift=10cm, yshift=7cm]
     \draw (90:3cm) -- (135:3cm) -- (180:3cm) -- (225:3cm) -- (270:3cm) -- (315:3cm) -- (360:3cm) -- (45:3cm) -- cycle;
                    \draw (90:3cm) -- node[midway, above, xshift=-1mm, yshift=-1mm] {} (180:3cm);
                    \draw (90:3cm) -- node[midway, above left,xshift=1.2mm] {} (225:3cm);
                    \draw (90:3cm) -- node[midway, above left,xshift=1mm] {} (270:3cm);
                    \draw (90:3cm) -- node[midway, above right, xshift=-1mm] {} (315:3cm);
                    \draw (90:3cm) -- node[midway, above right, xshift=-1mm] {} (360:3cm);
                   
\end{scope}
     \draw (90:3cm) -- (135:3cm) -- (180:3cm) -- (225:3cm) -- (270:3cm) -- (315:3cm) -- (360:3cm) -- (45:3cm) -- cycle;
                    \draw (90:3cm) -- node[midway, above, xshift=-1mm, yshift=-1mm] {1} (180:3cm);
                    \draw (90:3cm) -- node[midway, above left,xshift=1.2mm] {2} (225:3cm);
                    \draw (90:3cm) -- node[midway, above left,xshift=1mm] {3} (270:3cm);
                    \draw (90:3cm) -- node[midway, above right, xshift=-1mm] {4} (315:3cm);
                    \draw (90:3cm) -- node[midway, above right, xshift=-1mm] {5} (360:3cm);
                    \draw[cyan, line width=0.3mm] (135:3cm) -- (315:3cm); 
                    \draw[cyan, line width=0.3mm] (225:3cm) -- (45:3cm);
                    \draw[green, line width=0.3mm] (135:3cm) -- (45:3cm);
                    %\draw[red, line width=0.3mm] (225:3cm) -- (315:3cm);
                    \draw[magenta, line width=0.3mm] (135:3cm) -- (360:3cm);
                    \draw[magenta, line width=0.3mm] (45:3cm) -- (180:3cm);
                 
\node at (9,2) {$Q(T'): 1 \leftarrow 2 \leftarrow 3 \leftarrow 4 \leftarrow 5$};
\node at (10,-1) {$\{x_{\tabbedCenterstack{\textcolor{green}{5} \\ \textcolor{green}{4} \\ \textcolor{green}{3} \\ \textcolor{green}{2} \\ \textcolor{green}{1}}}, x_{\tabbedCenterstack{\textcolor{magenta}{4} \\ \textcolor{magenta}{3} \\ \textcolor{magenta}{2} \\ \textcolor{magenta}{1}}\oplus \tabbedCenterstack{\textcolor{magenta}{5} \\ \textcolor{magenta}{4} \\ \textcolor{magenta}{3} \\ \textcolor{magenta}{2}}} ,x_{\tabbedCenterstack{\textcolor{cyan}{3} \\ \textcolor{cyan}{2} \\ \textcolor{cyan}{1}}\oplus \tabbedCenterstack{\textcolor{cyan}{5} \\ \textcolor{cyan}{4} \\ \textcolor{cyan}{3}}} \}$ is a cluster of $\mathcal{A}^B (T)$, while};

\node at (9.5,-4.5) {$\{x_{\tabbedCenterstack{\textcolor{green}{5} \\ \textcolor{green}{4} \\ \textcolor{green}{3} \\ \textcolor{green}{2} \\ \textcolor{green}{1}}\oplus\tabbedCenterstack{\textcolor{green}{5} \\ \textcolor{green}{4} \\ \textcolor{green}{3} \\ \textcolor{green}{2} \\ \textcolor{green}{1}}}, x_{\tabbedCenterstack{\textcolor{magenta}{4} \\ \textcolor{magenta}{3} \\ \textcolor{magenta}{2} \\ \textcolor{magenta}{1}}\oplus \tabbedCenterstack{\textcolor{magenta}{5} \\ \textcolor{magenta}{4} \\ \textcolor{magenta}{3} \\ \textcolor{magenta}{2}}} ,x_{\tabbedCenterstack{\textcolor{cyan}{3} \\ \textcolor{cyan}{2} \\ \textcolor{cyan}{1}}\oplus \tabbedCenterstack{\textcolor{cyan}{5} \\ \textcolor{cyan}{4} \\ \textcolor{cyan}{3}}} \}$ is a cluster of $\mathcal{A}^C (T)$.};
                   
\end{tikzpicture}
        \caption{An example of cluster for a cluster algebra of type $B_3$ and $C_3$.}
        \label{fig:enter-label}
    \end{figure}
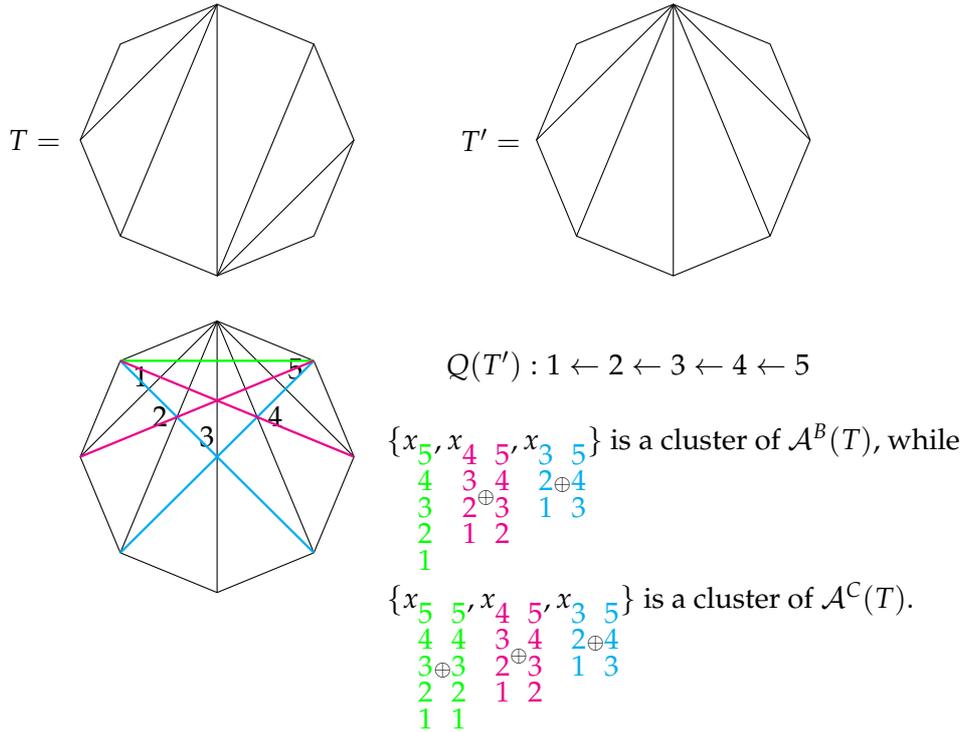

Finally, Theorem \ref{theorem1} and Theorem \ref{theorem 2} give two formulas (the former for type $B_n$ and the latter for type $C_n$) to express each cluster variable associated to a $\theta$-orbit, on the one hand in terms of the cluster variables of the initial seed, on the other hand in terms of cluster variables of type $A_n$. It follows from the above correspondence that, given a cluster-tilted bound symmetric quiver $Q$ of type $A_{2n-1}$ with no fixed arrows, they allow us to express the type $B_n$ (resp. type $C_n$) cluster variable that corresponds to an orthogonal (resp. symplectic) indecomposable representation of $Q$, on the one hand in terms of the initial cluster variables, on the other hand in terms of (ordinary) representations of $Q(\Bar{T})$, where $\Bar{T}=\text{Res}(T')$, and $T'$ is the triangulation of $\Poly_{2n+2}$ such that $Q=Q(T')$. In other words, we get a Caldero-Chapoton like map (see \cite{CC}) from the category of symmetric representations of cluster tilted bound symmetric quivers of type $A_{2n-1}$ (with no fixed arrows) to cluster algebras of type $B_n$ and $C_n$.

\begin{remark}
    The techniques presented in this section could be used to produce a categorification of other classes of skew-symmetrizable cluster algebras through the representation theory of symmetric quivers. For example, they could provide an alternative categorification of non skew-symmetric cluster algebras associated by Felikson, Shapiro and Tumarkin \cite{FeST} to surfaces with marked points and order-2 orbifold points. These algebras have been categorified in the work of Geuenich and Labardini-Fragoso \cite{GeuLF1,GeuLF2} by species with potential. 
\end{remark}

\subsection{Categorical interpretation of Theorem \ref{theorem1} in the acyclic case}\label{cat_interpretation}

In this section we assume that $Q$ is an acyclic quiver with $n$ vertices. 

First, we recall the cluster multiplication formula of \cite{CEHR}, Section 7. Then, we use it to obtain a categorical interpretation of Theorem \ref{theorem1}.

    Let X, S be $Q$-representations such that 
    $\text{dim\hspace{-0.05cm} Ext}^1(S,X)=1$. Then, by the Auslander-Reiten formulas, there are nonzero morphisms $f:X \to \tau S$ and $g: \tau^{-1}X \to S$ which are unique up to scalar, where $\tau$ is the Auslander-Reiten translation. 
    
    We use the following notation from \cite{CEHR}:
    
\begin{center}
    $X_S:=\text{ker}(f) \subset X$; \hspace{2cm}$S^X:=\text{im}(g) \subseteq S$.
\end{center}

Let $M$ be a finite
dimensional representation of $Q$. The $\bold{g}$–$vector$ \cite{DWZ} of $M$ is the integer vector $\bold{g}_M \in \mathbb{Z}_{\geq 0}^n$ given by $(\bold{g}_M)_i := - \langle S_i , M \rangle$, where $S_i$ is the simple at vertex $i$, and $\langle -, - \rangle$ is the Euler-Ringel form of $Q$. Let $B$ be the exchange matrix of $Q$. The $CC$–map is a map $M \mapsto CC(M)$
which associates to $M$ a Laurent polynomial $CC(M) \in \mathbb{Z}[y_1 ,\dots ,y_n ,x_1^{\pm 1},\dots,x_n^{\pm 1}]$, defined as follows
\begin{center}
    $CC(M) := \displaystyle\sum_{\bold{e} \in \mathbb{Z}_{\geq 0}^n} \chi(Gr_{\bold{e}}(M))\bold{y}^{\bold{e}}\bold{x}^{B\bold{e}+\bold{g}_M}$,
\end{center}
where $Gr_\bold e(M)$ is the quiver Grassmannian. Moreover, the $F$-$polynomial$ \cite{DWZ} of $M$ is $F_M:=CC(M)_{|x_1=\dots=x_n=1}$.

Let X, S be $Q$-representations such that $\text{dim \hspace{-0.05cm}Ext}^1(S,X)=1$. Then, by \cite[Lemma 31]{CEHR}, there exists an exact sequence $0 \to X/X_S \to \tau S^X \to I \to 0$,
where $I$ is either injective or zero. Let $I = I_1^{f_1} \oplus I_2^{f_2} \oplus \cdots \oplus I_n^{f_n}$  be the indecomposable decomposition of $I$, and let $\bold{f} = (f_1 ,\cdots ,f_n)$.

\begin{theorem}[\cite{CEHR}, Theorem 67]\label{thm67}
   Let X, S be $Q$-representations such that 
    
    $\text{dim \hspace{-0.05cm}Ext}^1(S,X)=1$. Let $\xi \in \mathrm{Ext}^1(S,X)$ be a non-split short exact sequence with middle term $Y$. Then
    \begin{equation}\label{eq_thm67}
        CC(X)CC(S) = CC(Y) + \bold{y}^{\textbf{dim} S^X} CC(X_S \oplus S/S^X )\bold{x}^{\bold{f}}.
    \end{equation}
  
    Moreover, if $\mathrm{Ext}^1(X,S) = 0$, and both $X$ and $S$ are rigid and indecomposable, then formula \ref{eq_thm67} is an exchange relation between the cluster variables $CC(X)$ and $CC(S)$ for the cluster algebra $\mathcal{A}(\bold{x},\bold{y},B)$ with principal coefficients at the initial seed $(\bold{x},\bold{y}, B)$.
 
\end{theorem}

\begin{remark}
    Let $Q$ be a symmetric quiver, and let $L$ be an ordinary representation of $Q$ such that 
    $\text{dim \hspace{-0.05cm}Ext}^1(\nabla L, L)=1$. By definition, $L_{\nabla L}= \text{ker}(L \to \tau \nabla L)$, and $\nabla L^L= \text{im}(\tau^{-1} L \to \nabla L)$. So we have that $\nabla (L_{\nabla L})=\text{coker}(\tau^{-1} L \to \nabla L)=\nabla L / \nabla L^L$, where we have used the fact that $\nabla \tau = \tau^{-1} \nabla$ (\cite{DW}, Proposition 3.4). Therefore, $L_{\nabla L}\oplus \nabla L / \nabla L^L$ is a symmetric representation of $Q$.
\end{remark}

Now, let $Q$ be a symmetric quiver of type $A_{2n-1}$. Observe that, in this case, if $M$ is a representation of $Q$, then 
\begin{equation}
    CC(M)= \displaystyle\sum_{\{\bold{e}=\textbf{dim}N\in \mathbb{Z}^n | N\subseteq M\}} \bold{y}^{\bold{e}}\bold{x}^{B\bold{e}+\bold{g}_M},
\end{equation}
since $Gr_\bold e(M)$ is either empty or a point.

Let $T'$ be the triangulation of $\Poly_{2n+2}$ such that $Q=Q(T')$. Since $Q$ has a fixed vertex $n$ and no fixed arrows, then $T'$ contains a diameter $d=\tau_n$, and if $\rho$ is the reflection along $d$, $T'$ is $\rho$-invariant. Let $[a,b]=\{(a,b),(\Bar{b},\Bar{a})\}$ be a $\theta$-orbit such that each diagonal of $[a,b]$ crosses $d$, so $\text{Res}([a,b])=\{(a,\ast),(\Bar{b},\ast)\}$, and let $(a,\Bar{a}), (\Bar{b},b)$ be the diameters starting in $a$ and $\Bar{b}$ respectively, so that $\text{Res}([a,\Bar{a}])=\{(a,\ast)\}$ and $\text{Res}([b,\Bar{b}])=\{(\Bar{b},\ast)\}$, see Figure \ref{fig_22} (the restriction is with respect to $d$). Therefore $[a,b]$ corresponds via $F_d$ to $L_{(a,\rho({\Bar{b})})}\oplus \nabla L_{(a,\rho({\Bar{b})})}$, with $\text{dim Ext}^1(\nabla L_{(a,\rho({\Bar{b})})}, L_{(a,\rho({\Bar{b})})})$ $=1$. Then, there exists a non-degenerate square in the Auslander-Reiten quiver of $Q$ from $L_{(a,\rho({\Bar{b})})}$ to $\nabla L_{(a,\rho({\Bar{b})})}=L_{(\Bar{b},\rho(a))}$, whose middle vertices $L_{(a,\rho(a))}$, $L_{(\Bar{b},\rho({\Bar{b})})}$ are $\nabla$-invariant. In other words, there is the non-split short exact sequence
\begin{equation}\label{ses}
    0 \to L_{(a,\rho({\Bar{b})})} \to L_{(a,\rho(a))} \oplus L_{(\Bar{b},\rho({\Bar{b})})} \to \nabla L_{(a,\rho({\Bar{b})})} \to 0.
\end{equation}

 \begin{figure}[H]
        \centering
\begin{tikzpicture}[scale=0.5]
    \draw (90:3cm) -- (120:3cm) -- (150:3cm) -- (180:3cm) -- (210:3cm) -- (240:3cm) -- (270:3cm) -- (300:3cm) -- (330:3cm) -- (360:3cm) -- (30:3cm) -- (60:3cm) --  cycle;
       \draw (120:3cm) -- (360:3cm); 
       \draw (180:3cm) -- (300:3cm); 
         \draw[dashed,-{Latex[length=2mm]}](270:3cm) -- node[midway, left, xshift=0.4cm,yshift=-0.5cm] {$d$} (90:3cm);
         \node at (120:3) [left] {$a$};
         \node at (365:3) [right,yshift=-0.1cm] {$b$};
         \node at (300:3) [right] {$\Bar{a}$};
         \node at (180:3) [left] {$\Bar{b}$};                 
\node at (6.3,0) [left] {$\xlongleftrightarrow{F_d}$};                 
\begin{scope}[xshift=10cm]
   \draw (90:3cm) -- (120:3cm) -- (150:3cm) -- (180:3cm) -- (210:3cm) -- (240:3cm) -- (270:3cm) -- (300:3cm) -- (330:3cm) -- (360:3cm) -- (30:3cm) -- (60:3cm) --  cycle;
     \draw (360:3cm) -- (120:3cm); 
     \draw (60:3cm) -- (180:3cm); 
        \node at (360:3) [right] {$\rho({\Bar{b})}$}; 
       \draw[dashed,-{Latex[length=2mm]}](270:3cm) -- node[midway, left, xshift=0.4cm,yshift=-0.5cm] {$d$} (90:3cm) ;
        \node at (120:3) [left] {$a$};
        \node at (180:3) [left,xshift=0.1cm] {$\Bar{b}$};
         \node at (60:3) [right] {$\rho(a)$};
\end{scope}
                   
\end{tikzpicture}
        \caption{The action of $F_d$ on the $\theta$-orbit $[a,b]$ whose diagonals cross $d$.}
        \label{fig_22}
    \end{figure} 

By Theorem \ref{thm67}, we have that 
\begin{align*}
    F_{L_{(a,\rho({\Bar{b})})}\oplus \nabla L_{(a,\rho({\Bar{b})})}}=F_{L_{(a,\rho(a))}\oplus L_{(\Bar{b},\rho({\Bar{b})})}}+\bold{y}^{\textbf{dim}\nabla L_{(a,\rho({\Bar{b})})}^{L_{(a,\rho({\Bar{b})})}}}F_{(L_{(a,\rho({\Bar{b})})})_{\nabla L_{(a,\rho({\Bar{b})})}}\oplus \nabla L_{(a,\rho({\Bar{b})})} / \nabla L_{(a,\rho({\Bar{b})})}^{L_{(a,\rho({\Bar{b})})}}}.
\end{align*}
On the other hand, by Proposition \ref{up:skein1},
\begin{align*}
    F_{L_{(a,\rho(\Bar{b}))}\oplus \nabla L_{(a,\rho(\Bar{b}))}}&=F_{L_{(a,\rho(a))}\oplus L_{(\Bar{b},\rho(\Bar{b}))}}+\bold{y}^{\bold{d}_{a\rho(a),\Bar{b}\rho(\Bar{b})}}F_{L_{(a,\Bar{b})}\oplus \nabla L_{(a,\Bar{b})}}.
\end{align*}
Thus
\begin{align*}
  \textbf{dim}\nabla L_{(a,\rho(\Bar{b}))}^{L_{(a,\rho(\Bar{b}))}}&=\bold{d}_{a\rho(a),\Bar{b}\rho(\Bar{b})},   
\end{align*}
and
\begin{align*}
   (L_{(a,\rho(\Bar{b}))})_{\nabla L_{(a,\rho(\Bar{b}))}}\oplus \nabla L_{(a,\rho({\Bar{b})})} / \nabla L_{(a,\rho(\Bar{b}))}^{L_{(a,\rho(\Bar{b}))}}=L_{(a,\Bar{b})}\oplus\nabla L_{(a,\Bar{b})}.
\end{align*}

Let $\mathcal{A}^B(T)$ be the cluster algebra of type $B_n$ with principal coefficients in the $\theta$-invariant triangulation $T$ of $\Poly_{2n+2}$ in the isotopy class of $F_d(T')$. Let $M$ be an orthogonal indecomposable representation of $Q(T')$. We denote by $F_M$ and $\bold{g}_M$ the $F$-polynomial and the $\bold{g}$-vector respectively of the cluster variable of $\mathcal{A}^B(T)$ that corresponds to $M$, and by $F_{\text{Res}(M)}$ and $\bold{g}_{\text{Res}(M)}$ the $F$-polynomial and the $\bold{g}$-vector respectively of the $Q(T')$-representation $\text{Res}(M)$. Then from the above discussion, it follows that Theorem \ref{theorem1} can be reformulated as:

\begin{comment}
    
Moreover, let $\Bar{T}$ be the triangulation of $\Poly_{n+3}$ obtained from $T$ by identifying the vertices which lie on the right of $d$, and let $\mathcal{A}(\Bar{T})=\mathcal{A}(\bold{u}_{\Bar{T}},\bold{y}_{\Bar{T}},B(\Bar{T}))$ be the cluster algebra of type A with principal coefficients in $\Bar{T}$.\end{comment}

\begin{theorem}\label{cat_interpr}
Let $M$ be an orthogonal indecomposable $Q(T')$-representation. If $\text{Res}(M)=(V,f)$ is indecomposable as $Q(T')$-representation, then 
    \begin{equation}\label{e_1}
        \text{$F_M=F_{\text{Res}(M)}$,}
    \end{equation}
    and
    \begin{equation}\label{e_2}
 \bold{g}_{M}=\begin{cases}
          \text{$D \bold{g}_{\text{Res}(M)}$ \hspace{1.8cm}if $\textbf{dim} V_n =0$;}\\
          \text{$D \bold{g}_{\text{Res}(M)}+\bold{e}_n$ \hspace{1cm}if $\textbf{dim} V_n \neq 0$.}
      \end{cases}
      \end{equation}
      Otherwise, $M=L\oplus \nabla L$ with $\text{dim Ext}^1(\nabla L, L)=1$, and there exists a non-split short exact sequence
    \begin{center}
        $0 \to L \to G_1 \oplus G_2 \to \nabla L \to 0 $,
    \end{center}
    where $G_1$ and $G_2$ are orthogonal indecomposable $Q(T')$-representations of type I. Then
    \begin{equation}\label{e_3}
        F_M= F_{\text{Res}(M)} - \bold{y}^{\text{Res}(\textbf{dim}\nabla L^L)}F_{\text{Res}(L_{\nabla L}\oplus \nabla L / \nabla L^L)},
    \end{equation}
    %\begin{equation}
        %=F_{\text{Res}(G_1)}F_{\text{Res}(G_2)}-\bold{y}^{\text{Res}(\text{dim}\nabla L^L)}F_{\text{Res}(L_{\nabla L}\oplus \nabla L / \nabla L^L)},
    %\end{equation}
    and
    \begin{equation}\label{e_4}
        \bold{g}_{M}=D(\bold{g}_{\text{Res}(M)}+\bold{e}_n).
    \end{equation}

\end{theorem}

\begin{remark}
    Observe that on the right hand sides of \ref{e_1}, \ref{e_2}, \ref{e_3}, \ref{e_4} we have only $F$-polynomials and $\bold{g}$-vectors of ordinary type $A$ quiver representations.  
\end{remark}

\begin{example}
Let  \[
\begin{tikzcd}
& 2 \arrow[dl] \arrow[d] \\
Q:\hspace{0.3cm}1  & 3 \arrow[d] & 5 \arrow[dl]  \\
&  4 
\end{tikzcd}
\] be the quiver of Example \ref{ex_quiver}. We compute the $F$-polynomial and the $\bold{g}$-vector of Example \ref{example_F_poly} using Theorem \ref{cat_interpr}. Let $M=L \oplus \nabla L=\begin{smallmatrix}
    35\\4
\end{smallmatrix}\oplus\begin{smallmatrix}
    2\\13
\end{smallmatrix}$ be an orthogonal indecomposable $Q$-representation. We have the short exact sequence
\begin{center}
       $0 \to \begin{smallmatrix}
           35\\4
       \end{smallmatrix} \to \begin{smallmatrix}
           2\\135\\4
       \end{smallmatrix}\oplus \begin{smallmatrix}
           3
       \end{smallmatrix} \to \begin{smallmatrix}
           2\\13
       \end{smallmatrix}\to 0$.
   \end{center}

   Since the sequence is almost split, $L_{\nabla L}=0$ and $\nabla L^{L}=\nabla L$. Therefore 

\begin{center}
    $F_M=F_{\text{Res}(  \begin{smallmatrix}
       35\\4
   \end{smallmatrix}\oplus\begin{smallmatrix}
       2\\13
   \end{smallmatrix})}-\bold{y}^{\displaystyle\text{Res}(\textbf{dim}\begin{smallmatrix}
      2\\13 
   \end{smallmatrix})}=F_{\begin{smallmatrix}
       3
   \end{smallmatrix}}F_{\begin{smallmatrix}
       2\\13
   \end{smallmatrix}}-y_1y_2y_3=y_1y_2y_3^2+y_1y_3^2+2y_1y_3+y_3^2+y_1+2y_3+1$.
\end{center}
 On the other hand, the $\bold{g}$-vector is 
\begin{center}
$\displaystyle\bold{g}_{M}=D(\bold{g}_{\bold{g}_{\text{Res}(M)}}+\bold{e}_3)=D(\bold{g}_{\begin{smallmatrix}
    3
\end{smallmatrix}\oplus \begin{smallmatrix}
    2\\13
\end{smallmatrix}}+\bold{e}_3) =D(\displaystyle\begin{pmatrix}
    -1 \\ 2 \\ -2
\end{pmatrix}+\displaystyle\begin{pmatrix}
    0 \\ 0 \\ 1
\end{pmatrix})=D(\displaystyle\begin{pmatrix}
    -1 \\ 2 \\ -1
\end{pmatrix})=\displaystyle\begin{pmatrix}
    -1 \\ 2 \\ -2
\end{pmatrix}$.   
\end{center}  
   
\end{example}

\section{Proof of Theorem \ref{theorem 2}}\label{sec_proof}
In order to present the proof of Theorem \ref{theorem 2}, we first need some lemmas.
\begin{lemma}\label{lemma_2}
    If each diagonal of $[a,b]$ crosses only one diagonal of $T$, then $F_{ab}=F_{ab}^C$ and $\bold{g}_{ab}=\bold{g}_{ab}^C$. 
\end{lemma}

\begin{proof}
   With the notation of the proof of Lemma \ref{lemma1}, $\Tilde{\text{Res}}([a,b])=\text{Res}([a,b])=\{\gamma_j\}$, where $\gamma_j$ is the diagonal of $\Poly_{n+3}$ which crosses only $\tau_j$. Let $B(\Bar{T})D=(b_{ij})$ and $B(\Bar{T})=(\Bar{b}_{ij})$.
    We have 
    \begin{equation}
        x_{ab}u_j=y_j\displaystyle\prod_{b_{ij}>0}u_i^{b_{ij}}+\displaystyle\prod_{b_{ij}<0}u_i^{-b_{ij}},
    \end{equation}
     and
     \begin{equation}
        x_{\gamma_j}u_j=y_j\displaystyle\prod_{\Bar{b}_{ij}>0}u_i^{\Bar{b}_{ij}}+\displaystyle\prod_{\Bar{b}_{ij}<0}u_i^{-\Bar{b}_{ij}}.
    \end{equation}
    So
    \begin{equation}
        F_{ab}=y_j+1=F_{\gamma_j}=F_{ab}^C.
    \end{equation}

If $j=n$ and $k$ is such that $\tau_k$ and $\tau_n$ are both sides of a triangle of $T$, and $\tau_k$ is clockwise from $\tau_n$, then $b_{kn}=-2$, while $\Bar{b}_{kn}=-1$. So
    \begin{equation}
        (\bold{g}_{ab})_k=\displaystyle\bigg(\text{deg}\displaystyle\bigg(\displaystyle\frac{\displaystyle\prod_{b_{in}<0}u_i^{-b_{in}}}{u_n}\bigg)\bigg)_k=\displaystyle\bigg(\text{deg}\displaystyle\bigg(\displaystyle\frac{\displaystyle\prod_{\Bar{b}_{in}<0}u_i^{-\Bar{b}_{in}}}{u_n}\bigg)\bigg)_k + 1= (\bold{g}_{\gamma_n})_k +1=
        (\bold{g}_{ab}^C)_k.
    \end{equation}
   Otherwise, 
    \begin{equation}
        (\bold{g}_{ab})_k=\displaystyle\bigg(\text{deg}\displaystyle\bigg(\displaystyle\frac{\displaystyle\prod_{b_{in}<0}u_i^{-b_{in}}}{u_n}\bigg)\bigg)_k=\displaystyle\bigg(\text{deg}\displaystyle\bigg(\displaystyle\frac{\displaystyle\prod_{\Bar{b}_{in}<0}u_i^{-\Bar{b}_{in}}}{u_n}\bigg)\bigg)_k=(\bold{g}_{\gamma_n})_k=
        (\bold{g}_{ab}^C)_k.
    \end{equation}
\end{proof}

\begin{lemma}\label{lemma_c_vect_typeC}
Let $B$ be a skew-symmetric $n\times n$ matrix, and let $I$ be the $n\times n$ identity matrix. Let
$D=\text{diag}(1,\dots,1,2)$ be $n\times n$ diagonal matrix with diagonal entries $(1,\dots,1,2)$. 
\begin{itemize}
    \item [i)] Let 
$\mu_{i_1} \cdots \mu_{i_k}(\left[\begin{matrix}
   $B$\\ I
\end{matrix} \right])=\left[\begin{matrix}
    B' \\ C
\end{matrix} \right]$, and let $\mu_{i_1} \cdots \mu_{i_k}(\left[\begin{matrix}
    BD \\ I
\end{matrix} \right])=\left[\begin{matrix}
    B'D \\ C'
\end{matrix} \right]$, for any 
$1\leq i_1 < \cdots < i_k \leq n$. Then, 
 $C^k=(C')^k$ for any $k\neq n$.
 \item [ii)] Let 
$\mu_{i_k} \cdots \mu_{i_1}(\left[\begin{matrix}
   $B$\\ I
\end{matrix} \right])=\left[\begin{matrix}
    B' \\ C
\end{matrix} \right]$, and let $\mu_{i_k} \cdots \mu_{i_1}(\left[\begin{matrix}
    BD \\ I
\end{matrix} \right])=\left[\begin{matrix}
    B'D \\ C'
\end{matrix} \right]$, for any 
$1 \leq i_1 < \cdots < i_k < n$. Then $((C')^n)_i=\begin{cases}
    \text{$2(C^n)_i$ if $i\neq n$,}\\
    \text{$(C^n)_n$ if $i=n$}.
\end{cases}$
 \end{itemize}
\end{lemma}

\begin{proof}
$B$ and $BD$ differ only in the $n$-th column, and the $n$-th column of $BD$ is equal to the $n$-th one of $B$ multiplied by 2.
$i)$ follows from the fact the 2 can appear in the bottom part of the matrix only in the $n$-th column, since we mutate at $n$ only eventually once at the beginning. In $ii)$, we start mutating from the left. So in the bottom part of the $n$-th column, other than the last coordinate, only the entries corresponding to $i_1,\dots,i_k$ can be nonzero. For each $j$, $\mu_{i_j}\cdots \mu_{i_1}(BD)=\mu_{i_j}\cdots \mu_{i_1}(B)D$, since the symmetrizer is constant in the mutation class of $B$ (\cite{CAII}, Proposition 4.5), i.e. $\mu_{i_j}\cdots \mu_{i_1}(BD)$ is equal to $\mu_{i_j}\cdots \mu_{i_1}(B)$ with the $n$-th column multiplied by 2. So for any $i\neq n$, $((C')^n)_i \neq 0$ if and only if $(C^n)_i \neq 0$, and $((C')^n)_i=2(C^n)_i$. Finally, $((C')^n)_n$ doesn't change after mutations, as well as $(C^n)_n$, so $((C')^n)_n=1=(C^n)_n$   
\end{proof}

\begin{proof}[Proof of Theorem \ref{theorem 2}]
    We prove the theorem by induction on the number $k$ of intersections between each diagonal of $[a,b]$ and $T=\{\tau_1, \dots, \tau_n=d,\dots,\tau_{2n-1}\}$.

    If $k=0$, the theorem holds by Definition \ref{def_type_C}. If $k=1$, the theorem holds by Lemma \ref{lemma_2}.
    Assume $k>1$.Let $\Bar{T}=\text{Res}(T)=\{ \tau_1, \dots, \tau_n=d \}$, and let $\bold{u}_{\Bar{T}}=\{ u_{\tau_1}, \dots, u_{\tau_n} \}=\{u_1, \dots, u_n\}$. There are three cases to consider.
    \begin{itemize}
\item [1)] Let $[a,b]=\{(a,b),(\Bar{b},\Bar{a})\}$ be such that $\Tilde{\text{Res}}([a,b])=\{(a,b)\}$. Let $a=p_0,p_1,\dots,p_k,$ $p_{k+1}=b$ be the intersection points of $(a,b)$ and $\Bar{T}$ in order of occurrence on $(a,b)$, and let $i_1,i_2,\dots,i_k$ be such that $p_j$ lies on the diagonal $\tau_{i_j}\in \Bar{T}$, for $j=1,\dots,k$. Let $[c,d]=\{\tau_{i_1},\tau_{i_{2n-i_1}}\}$.
  \begin{figure}[H]
                               \centering
                \begin{tikzpicture}[scale=0.5]
                    \draw (90:3cm) -- (120:3cm) -- (150:3cm) -- (180:3cm) -- (210:3cm) -- (240:3cm) -- (270:3cm) -- (300:3cm) -- (330:3cm) -- (360:3cm) -- (30:3cm) -- (60:3cm) --  cycle;
                    %\draw (90:3cm) -- node[midway, above, xshift=-1mm, yshift=-1mm] {} (180:3cm);
                    %\draw (90:3cm) -- node[midway, above left,xshift=1.2mm] {} (225:3cm);
                    \draw[-{Latex[length=2mm]}] (270:3cm) -- node[midway, above left,xshift=1mm] {} (90:3cm);
                   % \draw (90:3cm) -- node[midway, above right, xshift=-1mm] {} (315:3cm);
                    %\draw (90:3cm) -- node[midway, above right, xshift=-1mm] {} (360:3cm);
                     
                    \draw[red, line width=0.3mm] (330:3cm) -- (60:3cm);
                     \draw[red, line width=0.3mm] (150:3cm) -- (240:3cm); 
                     \draw[cyan, line width=0.3mm] (120:3cm) -- (210:3cm);
                    \draw[cyan, line width=0.3mm] (300:3cm) -- (30:3cm);

                    \node at (120:3cm) [left] {$c$};
                    \node at (210:3cm) [left] {$d$};
                    \node at (330:3cm) [right] {$\Bar{a}$};
                    \node at (60:3cm) [right] {$\Bar{b}$};

\node at (150:3cm) [left] {$a$};
                    \node at (240:3cm) [left] {$b$};
                    \node at (300:3cm) [right] {$\Bar{c}$};
                    \node at (30:3cm) [right] {$\Bar{d}$};

                    \begin{scope}[xshift=8cm]
                      \draw (90:3cm) -- (120:3cm) -- (150:3cm) -- (180:3cm) -- (210:3cm) -- (240:3cm) -- (270:3cm) -- (360:3cm)  --  cycle;
                    %\draw (90:3cm) -- node[midway, above, xshift=-1mm, yshift=-1mm] {} (180:3cm);
                    %\draw (90:3cm) -- node[midway, above left,xshift=1.2mm] {} (225:3cm);
                    \draw (90:3cm) -- node[midway, above left,xshift=1mm] {} (270:3cm);
                   % \draw (90:3cm) -- node[midway, above right, xshift=-1mm] {} (315:3cm);
                    %\draw (90:3cm) -- node[midway, above right, xshift=-1mm] {} (360:3cm);
                    \draw[red, line width=0.3mm] (150:3cm) -- (240:3cm); 
                     \draw[cyan, line width=0.3mm] (120:3cm) -- (210:3cm);

                    \node at (120:3cm) [left] {$c$};
                    \node at (210:3cm) [left] {$d$};
                    \node at (360:3cm) [right] {$\ast$};

\node at (150:3cm) [left] {$a$};
                    \node at (240:3cm) [left] {$b$};
                      
                    \end{scope}
                \end{tikzpicture}
                               \caption{On the left, the two $\theta$-orbits $[a,b]$ and $[c,d]$. On the right, their rotated restrictions.}
                               \label{fig:enter-label}
                           \end{figure}
Then, by Lemma \ref{lemma_schiffler}, $(a,b) \in \mu_{i_1}\cdots \mu_{i_k}(\Bar{T})$. Therefore, the $\bold{c}$-vector corresponding to the exchange between $[a,b]$ and $[c,d]$ is the bottom part of the $i_1$-th column of $\mu_{i_2}\cdots \mu_{i_k}(\left[\begin{matrix}
    B(\Bar{T})D \\ I
\end{matrix} \right])$. Since $i_1 \neq n$, by Lemma \ref{lemma_c_vect_typeC} i), this is equal to the bottom part of the $i_1$-th column of $\mu_{i_2}\cdots \mu_{i_k}(\left[\begin{matrix}
    B(\Bar{T}) \\ I
\end{matrix} \right])$, which is given by Proposition \ref{up:skein1}. Therefore, we have the following exchange relation
\begin{equation}\label{thm2_eqn1}
u_{i_1}x_{ab}=\bold{y}^{\bold{d}_{ac,bd}}x_{ad}x_{bc}+\bold{y}^{\bold{d}_{ad,bc}}x_{ac}x_{bd}.
\end{equation}
Since $(c,d)$ is the first diagonal of $T$ that is crossed by $(a,b)$, $(a,c)$ and $(a,d)$ must be either boundary edges or diagonals of $\Bar{T}$. It follows from \ref{thm2_eqn1} that
\begin{equation}
    F_{ab}=\bold{y}^{\bold{d}_{ac,bd}}F_{bc}+\bold{y}^{\bold{d}_{ad,bc}}F_{bd}.
\end{equation}
By inductive hypothesis and Proposition \ref{up:skein1},
\begin{equation}
    F_{ab}=\bold{y}^{\bold{d}_{ac,bd}}F_{(b,c)}+\bold{y}^{\bold{d}_{ad,bc}}F_{(b,d)}=F_{(a,b)}=F_{ab}^C.
\end{equation}
\begin{comment}
    It also follows from \ref{thm2_eqn1} that
\begin{equation}
    \bold{e}_{i_1}+\bold{g}_{ab}=\begin{cases}
        \bold{g}_{ad}+\bold{g}_{bc}\hspace{1cm} \text{if $\bold{y}^{\bold{d}_{ac,bd}}=1$}\\
        \bold{g}_{ac}+\bold{g}_{bd}\hspace{1cm} \text{otherwise}.
    \end{cases}
\end{equation}
By inductive hypothesis and Proposition \ref{up:skein1},

\begin{equation*}
    \bold{g}_{ab}=\begin{cases}
        -\bold{e}_{i_1}+ \bold{g}_{(a,d)}+\bold{g}_{(b,c)} \hspace{0.6cm}\text{if $\bold{y}^{\bold{d}_{ac,bd}}=1$}\\
        -\bold{e}_{i_1} + \bold{g}_{(a,c)} +\bold{g}_{(b,d)}\hspace{0.6cm} \text{otherwise}.
        \end{cases}=\bold{g}_{(a,b)}=\bold{g}_{ab}^C.
\end{equation*}
\end{comment}

\item [2)] Let $[a,\Bar{a}]$ be a diameter. So $\Tilde{\text{Res}}([a,\Bar{a}])=\{(a,\ast),(a,\Bar{b})\}$. Let $\ast=p_0,p_1,\dots,p_s,$ $p_{s+1} = a$ be the intersection points of $(a,\ast)$ and $\Bar{T}$ in order of occurrence on $(\ast,a)$, $s\leq k$, and let $i_1,i_2,\dots,i_s$ be such that $p_j$ lies on the diagonal $\tau_{i_j}\in \Bar{T}$, for $j=1,\dots,s$. Thus $i_1 = n$. Let $[b,\Bar{b}]=\{\tau_n\}=\{d\}$. We have two cases to consider:
\begin{itemize}
    \item [i)] there is no $i\in \{1,\dots,n\}$ such that $\tau_i$ and $\tau_n$ are both sides of a triangle of $T$, and $\tau_i$ is clockwise from $\tau_n$; \begin{figure}[H]
                          \centering
                \begin{tikzpicture}[scale=0.5]
                    \draw (90:3cm) -- (120:3cm) -- (150:3cm) -- (180:3cm) -- (210:3cm) -- (240:3cm) -- (270:3cm) -- (300:3cm) -- (330:3cm) -- (360:3cm) -- (30:3cm) -- (60:3cm) --  cycle;
                    %\draw (90:3cm) -- node[midway, above, xshift=-1mm, yshift=-1mm] {} (180:3cm);
                    %\draw (90:3cm) -- node[midway, above left,xshift=1.2mm] {} (225:3cm);
                    \draw (90:3cm) -- node[midway, above left,xshift=1mm] {} (270:3cm);
                   % \draw (90:3cm) -- node[midway, above right, xshift=-1mm] {} (315:3cm);
                    %\draw (90:3cm) -- node[midway, above right, xshift=-1mm] {} (360:3cm);
                     
                    \draw[-{Latex[length=2mm]},cyan,line width=0.3mm] (270:3cm) -- node[midway, above left,xshift=1mm] {} (90:3cm);
                     \draw[red, line width=0.3mm] (150:3cm) -- (330:3cm);

                    \node at (90:3cm) [above] {$b$};
                    \node at (270:3cm) [below] {$\Bar{b}$};

\node at (150:3cm) [left] {$a$};
                    \node at (330:3cm) [right] {$\Bar{a}$};
                    \begin{scope}[xshift=8cm]
                      \draw (90:3cm) -- (120:3cm) -- (150:3cm) -- (180:3cm) -- (210:3cm) -- (240:3cm) -- (270:3cm) -- (360:3cm)  --  cycle;
                    %\draw (90:3cm) -- node[midway, above, xshift=-1mm, yshift=-1mm] {} (180:3cm);
                    %\draw (90:3cm) -- node[midway, above left,xshift=1.2mm] {} (225:3cm);
                    \draw (90:3cm) -- node[midway, above left,xshift=1mm] {} (270:3cm);
                   % \draw (90:3cm) -- node[midway, above right, xshift=-1mm] {} (315:3cm);
                    %\draw (90:3cm) -- node[midway, above right, xshift=-1mm] {} (360:3cm);
                   
                    \draw[cyan, line width=0.3mm] (90:3cm) -- (270:3cm);
                     \draw[red, line width=0.3mm] (150:3cm) -- (360:3cm); 
                     \draw[red, line width=0.3mm] (150:3cm) -- (270:3cm);

  \node at (90:3cm) [above] {$b$};
                    \node at (270:3cm) [below] {$\Bar{b}$};
 \node at (240:3cm) [left] {$c$};
\node at (150:3cm) [left] {$a$};
                    \node at (360:3cm) [right] {$\ast$};
                      
                    \end{scope}
                \end{tikzpicture}
                               \caption{On the left, the two $\rho$-orbits $[a,\Bar{a}]$ and $[b,\Bar{b}]$. On the right, their rotated restrictions.}
                               \label{fig:enter-label}
                           \end{figure}    
\item [ii)] there exists $i \in \{1,\dots,n\}$ such that $\tau_i$ is clockwise from $\tau_n$. 
                            \begin{figure}[H]
                          \centering
                \begin{tikzpicture}[scale=0.5]
                    \draw (90:3cm) -- (120:3cm) -- (150:3cm) -- (180:3cm) -- (210:3cm) -- (240:3cm) -- (270:3cm) -- (300:3cm) -- (330:3cm) -- (360:3cm) -- (30:3cm) -- (60:3cm) --  cycle;
                    %\draw (90:3cm) -- node[midway, above, xshift=-1mm, yshift=-1mm] {} (180:3cm);
                    %\draw (90:3cm) -- node[midway, above left,xshift=1.2mm] {} (225:3cm);
                    \draw[-{Latex[length=2mm]},cyan,line width=0.3mm] (270:3cm) -- node[midway, above left,xshift=1mm] {} (90:3cm);
                   % \draw (90:3cm) -- node[midway, above right, xshift=-1mm] {} (315:3cm);
                    %\draw (90:3cm) -- node[midway, above right, xshift=-1mm] {} (360:3cm);
                     
                    \draw[cyan, line width=0.3mm] (90:3cm) -- (270:3cm);
                     \draw[red, line width=0.3mm] (150:3cm) -- (330:3cm);

                    \node at (90:3cm) [above] {$\Bar{b}$};
                    \node at (270:3cm) [below] {$b$};

\node at (150:3cm) [left] {$a$};
                    \node at (330:3cm) [right] {$\Bar{a}$};
                    \begin{scope}[xshift=8cm]
                      \draw (90:3cm) -- (120:3cm) -- (150:3cm) -- (180:3cm) -- (210:3cm) -- (240:3cm) -- (270:3cm) -- (360:3cm)  --  cycle;
                    %\draw (90:3cm) -- node[midway, above, xshift=-1mm, yshift=-1mm] {} (180:3cm);
                    %\draw (90:3cm) -- node[midway, above left,xshift=1.2mm] {} (225:3cm);
                    \draw (90:3cm) -- node[midway, above left,xshift=1mm] {} (270:3cm);
                   % \draw (90:3cm) -- node[midway, above right, xshift=-1mm] {} (315:3cm);
                    %\draw (90:3cm) -- node[midway, above right, xshift=-1mm] {} (360:3cm);
                   
                    \draw[cyan, line width=0.3mm] (90:3cm) -- (270:3cm);
                     \draw[red, line width=0.3mm] (150:3cm) -- (360:3cm); 
                     \draw[red, line width=0.3mm] (150:3cm) -- (90:3cm);

  \node at (90:3cm) [above] {$\Bar{b}$};
                    \node at (270:3cm) [below] {$b$};
 \node at (120:3cm) [above] {$c$};
\node at (150:3cm) [left] {$a$};
                    \node at (360:3cm) [right] {$\ast$};
                      
                    \end{scope}
                \end{tikzpicture}
                               \caption{On the left, the two $\rho$-orbits $[a,\Bar{a}]$ and $[b,\Bar{b}]$. On the right, their rotated restrictions.}
                               \label{fig:enter-label}
                           \end{figure}    
                           
\end{itemize}
We prove i). The proof of ii) is analogous.
 By Lemma \ref{lemma_schiffler}, $(a,\ast) \in \mu_{i_1}\cdots \mu_{i_s}(\Bar{T})$. Therefore, the $\bold{c}$-vector corresponding to the exchange between $[a,\Bar{a}]$ and $[b,\Bar{b}]$ is the bottom part of the $i_1$-th column of $\mu_{i_2}\cdots \mu_{i_s}(\left[\begin{matrix}
    B(\Bar{T})D \\ I
\end{matrix} \right])$. By Lemma \ref{lemma_c_vect_typeC} ii), this is equal to the bottom part of the $i_1$-th column of $\mu_{i_2}\cdots \mu_{i_s}(\left[\begin{matrix}
    B(\Bar{T}) \\ I
\end{matrix} \right])$, which is given by Proposition \ref{up:skein1}, with all coordinates multiplied by two except the $n$-th one. If $v \in \mathbb{Z}_{\geq 0}$, we indicate by $\Bar{v}$ the vector whose coordinates are multiplied by two but the $n$-th one. Therefore, we have the following exchange relation
\begin{equation}\label{thm2_eqn5}
u_{n}x_{a\Bar{a}}=\bold{y}^{\Bar{\bold{d}}_{ab,\Bar{b}\ast}}x_{a\Bar{b}}^2+\bold{y}^{\Bar{\bold{d}}_{b\ast,a\Bar{b}}}x_{ab}^2.
\end{equation}

We note that $\bold{y}^{\Bar{\bold{d}}_{ab,\Bar{b}\ast}}=1$, since it cannot exist $i$ such that $L_i$ intersects both $(a,b)$ and $(\Bar{b},\ast)$. 

It follows from \ref{thm2_eqn5} that
\begin{equation}
    F_{a\Bar{a}}=F_{a\Bar{b}}^2+\bold{y}^{\Bar{\bold{d}}_{b\ast,a\Bar{b}}}F_{ab}^2.
\end{equation}
By inductive hypothesis and repeated applications of Proposition \ref{up:skein1},
\begin{equation}
    F_{a\Bar{a}}=F_{(a,\Bar{b})}^2+\bold{y}^{\Bar{\bold{d}}_{b\ast,a\Bar{b}}}F_{(a,b)}^2=F_{(a,\ast)}F_{(a,\Bar{b})}-\bold{y}^{\Tilde{\text{Res}}(\bold{d}_{a\ast,c\Bar{b}}+\bold{d}_{a\Bar{b},b\ast})}F_{(a,b)}F_{(a,c)}=F_{a\Bar{a}}^C.
\end{equation}

        \item [3)] Let $[a,b]=\{(a,b),(\Bar{b},\Bar{a}\}$ be such that $\Tilde{\text{Res}([a,b])}=\{(a,\ast),(\Bar{b},\Bar{e})\}$.

        Let $a=p_0,p_1,\dots,p_s,p_{s+1} = \ast$ be the intersection points of $(a,\ast)$ and $\Bar{T}$ in order of occurrence on $(a,\ast)$, and let $i_1,i_2,\dots,i_s$ be such that $p_j$ lies on the diagonal $\tau_{i_j}\in \Bar{T}$, for $j=1,\dots,s$. So $i_s=n$. Let $[c,d]=\{\tau_{i_1},\tau_{i_{2n-i_1}}\}$. Assume that $(c,d)=\tau_{i_1}$ intersects $(a,\ast)$ (otherwise we consider $(\Bar{b},\Bar{e})$ instead of $(a,\ast)$).

        \begin{figure}[H]
                               \centering
                \begin{tikzpicture}[scale=0.5]
                    \draw (90:3cm) -- (120:3cm) -- (150:3cm) -- (180:3cm) -- (210:3cm) -- (240:3cm) -- (270:3cm) -- (300:3cm) -- (330:3cm) -- (360:3cm) -- (30:3cm) -- (60:3cm) --  cycle;
                    %\draw (90:3cm) -- node[midway, above, xshift=-1mm, yshift=-1mm] {} (180:3cm);
                    %\draw (90:3cm) -- node[midway, above left,xshift=1.2mm] {} (225:3cm);
                     \draw[-{Latex[length=2mm]}] (270:3cm) -- node[midway, above left,xshift=1mm] {} (90:3cm);
                   % \draw (90:3cm) -- node[midway, above right, xshift=-1mm] {} (315:3cm);
                    %\draw (90:3cm) -- node[midway, above right, xshift=-1mm] {} (360:3cm);
                     
                    \draw[cyan, line width=0.3mm] (300:3cm) -- (360:3cm);
                     \draw[red, line width=0.3mm] (150:3cm) -- (30:3cm); 
                     \draw[cyan, line width=0.3mm] (120:3cm) -- (180:3cm);
                    \draw[red, line width=0.3mm] (210:3cm) -- (330:3cm);

                    \node at (120:3cm) [left] {$c$};
                    \node at (90:3cm) [above,yshift=-1] {$e$};
                    \node at (270:3cm) [below,yshift=1] {$\Bar{e}$};
                    \node at (180:3cm) [left] {$d$};
                    \node at (360:3cm) [right] {$\Bar{d}$};
                    \node at (60:3cm) [right] {$\Bar{f}$};

\node at (150:3cm) [left] {$a$};
                    \node at (210:3cm) [left] {$\Bar{b}$};
                    \node at (330:3cm) [right] {$\Bar{a}$};
                    \node at (300:3cm) [right,xshift=-1] {$\Bar{c}$};
                    \node at (240:3cm) [left] {$f$};
                    \node at (30:3cm) [right] {$b$};

                    \begin{scope}[xshift=8cm]
                      \draw (90:3cm) -- (120:3cm) -- (150:3cm) -- (180:3cm) -- (210:3cm) -- (240:3cm) -- (270:3cm) -- (360:3cm)  --  cycle;
                    %\draw (90:3cm) -- node[midway, above, xshift=-1mm, yshift=-1mm] {} (180:3cm);
                    %\draw (90:3cm) -- node[midway, above left,xshift=1.2mm] {} (225:3cm);
                    \draw (90:3cm) -- node[midway, above left,xshift=1mm] {} (270:3cm);
                   % \draw (90:3cm) -- node[midway, above right, xshift=-1mm] {} (315:3cm);
                    %\draw (90:3cm) -- node[midway, above right, xshift=-1mm] {} (360:3cm);
                    \draw[red, line width=0.3mm] (150:3cm) -- (360:3cm); 
                     \draw[cyan, line width=0.3mm] (120:3cm) -- (180:3cm);
\draw[red, line width=0.3mm] (210:3cm) -- (270:3cm);

                  \node at (120:3cm) [left] {$c$};
                    \node at (90:3cm) [above,yshift=-1] {$e$};
                    \node at (270:3cm) [below,yshift=1] {$\Bar{e}$};
                    \node at (180:3cm) [left] {$d$};
                    \node at (360:3cm) [right] {$\ast$};

\node at (150:3cm) [left] {$a$};
                    \node at (210:3cm) [left] {$\Bar{b}$};
                    \node at (240:3cm) [left] {$f$};                     
                    \end{scope}
                \end{tikzpicture}
                               \caption{On the left, the two $\theta$-orbits $[a,b]$ and $[c,d]$. On the right, their rotated restrictions.}
                               \label{fig:enter-label}
                           \end{figure}
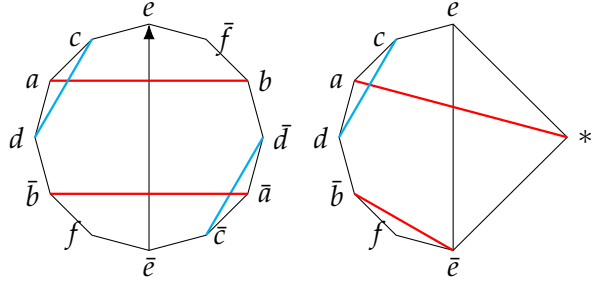

Then, by Lemma \ref{lemma_schiffler}, $(a,\ast) \in \mu_{i_1}\cdots \mu_{i_s}(\Bar{T})$. Therefore, the $\bold{c}$-vector corresponding to the exchange between $[a,b]$ and $[c,d]$ is the bottom part of the $i_1$-th column of $\mu_{i_2}\cdots \mu_{i_s}(\left[\begin{matrix}
    B(\Bar{T})D \\ I
\end{matrix} \right])$. By Lemma \ref{lemma_c_vect_typeC} i), this is equal to $C^{i_1}$, where $C^{i_1}$ is the bottom part of the $i_1$-th column of $\mu_{i_2}\cdots \mu_{i_s}(\left[\begin{matrix}
    B(\Bar{T}) \\ I
\end{matrix} \right])$, which is given by Proposition \ref{up:skein1}.

We have the following exchange relation:
\begin{equation}\label{thm2_eqn3}
u_{i_1}x_{ab}=\bold{y}^{\bold{d}_{ac,d\ast}}x_{ad}x_{bc}+\bold{y}^{\bold{d}_{ad,c\ast}}x_{ac}x_{bd}.
\end{equation}
 
It follows from \ref{thm2_eqn3} that
\begin{equation}
F_{ab}=\bold{y}^{\bold{d}_{ac,d\ast}}F_{bc}+\bold{y}^{\bold{d}_{ad,c\ast}}F_{bd},
\end{equation}
where we have used that $F_{ad}=F_{ac}=1$, since $[a,d]$ and $[a,c]$ must be either boundary edges or pairs of diagonals of $T$.

By inductive hypothesis and repeated applications of Proposition \ref{up:skein1},

    $F_{ab}=\bold{y}^{\bold{d}_{ac,d\ast}}(F_{(c,\ast)}F_{(\Bar{b},\Bar{e})}-\bold{y}^{\Tilde{\text{Res}}(\bold{d}_{\Bar{b}\ast,f\Bar{e}}+\bold{d}_{c\Bar{e},e\ast})}F_{(c,e)}F_{(\Bar{b},f)})+\bold{y}^{\bold{d}_{ad,c\ast}}(F_{(d,\ast)}F_{(\Bar{b},\Bar{e})}$
    
    $-\bold{y}^{\Tilde{\text{Res}}(\bold{d}_{\Bar{b}\ast,f\Bar{e}}+\bold{d}_{d\Bar{e},e\ast})}F_{(d,e)}F_{(\Bar{b},f)})$
    $=F_{(a,\ast)}F_{(\Bar{b},\Bar{e})}-\bold{y}^{\Tilde{\text{Res}}(\bold{d}_{\Bar{b}\ast,f\Bar{e}}+\bold{d}_{a\Bar{e},e\ast})}F_{(a,e)}F_{(\Bar{b},f)}=F_{ab}^C$.
\begin{comment}
    It also follows from \ref{thm2_eqn3} that
\begin{equation}
    \bold{e}_{i_1}+\bold{g}_{ab}=\begin{cases}
        \bold{g}_{ad}+\bold{g}_{bc}\hspace{1cm} \text{if $\bold{y}^{\bold{d}_{ac,d\ast}}=1$}\\
        \bold{g}_{ac}+\bold{g}_{bd}\hspace{1cm} \text{otherwise}.
    \end{cases}
\end{equation}

By inductive hypothesis and Proposition \ref{up:skein1},

\begin{align*}
    \bold{g}_{ab}&=\begin{cases}
        -\bold{e}_{i_1} + \bold{g}_{(a,d)} + \bold{g}_{(c,\ast)} (+\bold{e}_i) + \bold{g}_{(\Bar{b},\Bar{e})}  \hspace{0.8cm}\text{if $\bold{y}^{\bold{d}_{ac,d\ast}}=1$}\\
        -\bold{e}_{i_1} + \bold{g}_{(a,c)} + \bold{g}_{(d,\ast)} (+\bold{e}_i) +\bold{g}_{(\Bar{b},\Bar{e})} \hspace{0.8cm} \text{otherwise};
        \end{cases}\\
    &=\bold{g}_{(a,\ast)}(+\bold{e}_i)+\bold{g}_{(\Bar{b},\Bar{e})}=\bold{g}_{ab}^C.
\end{align*}

\end{comment}

    \end{itemize}
Similarly we prove that $\bold{g}_{ab}=\bold{g}_{ab}^C$.
\end{proof}

\addcontentsline{toc}{section}{Acknowledgments}
\section*{Acknowledgments}
The results presented in this paper are part of my Ph.D. thesis under the supervision of Giovanni Cerulli Irelli. I thank him for his support, and helpful comments and suggestions on a previous version of this work. I am also grateful to Ralf Schiffler, Pierre-Guy Plamondon, Salvatore Stella, Bernhard Keller and Daniel Labardini-Fragoso for many helpful discussions. Finally, I thank the anonymous referees for their valuable comments.

This work has been partially supported by the "National Group for Algebraic and Geometric Structures, and their Applications" (GNSAGA - INdAM), and by NextGenerationEU under NRRP, Call PRIN 2022  No.~104 of February 2, 2022 of Italian Ministry of University and Research; Project 2022S97PMY \textit{Structures for Quivers, Algebras and Representations (SQUARE)}.

\printbibliography[heading=bibintoc]

@article {Y-systems,
    AUTHOR = {Fomin, S. and Zelevinsky, A.},
     TITLE = {{$Y$}-systems and generalized associahedra},
   JOURNAL = {Ann. of Math. (2)},
  FJOURNAL = {Annals of Mathematics. Second Series},
    VOLUME = {158},
      YEAR = {2003},
    NUMBER = {3},
     PAGES = {977--1018},
      ISSN = {0003-486X},
   MRCLASS = {17B20 (20F55 52B12)},
  MRNUMBER = {2031858},
MRREVIEWER = {Ivan Arzhantsev},
       DOI = {10.4007/annals.2003.158.977},
       URL = {https://doi.org/10.4007/annals.2003.158.977},
}

@article {CAII,
    AUTHOR = {Fomin, S. and Zelevinsky, A.},
     TITLE = {Cluster algebras. {II}. {F}inite type classification},
   JOURNAL = {Invent. Math.},
  FJOURNAL = {Inventiones Mathematicae},
    VOLUME = {154},
      YEAR = {2003},
    NUMBER = {1},
     PAGES = {63--121},
      ISSN = {0020-9910},
   MRCLASS = {17B20 (05E15 16S99 52B12)},
  MRNUMBER = {2004457},
MRREVIEWER = {Eric N. Sommers},
       DOI = {10.1007/s00222-003-0302-y},
       URL = {https://doi.org/10.1007/s00222-003-0302-y},
}

@article {CSI,
    AUTHOR = {Canakci, I. and Schiffler, R.},
     TITLE = {Snake graph calculus and cluster algebras from surfaces},
   JOURNAL = {J. Algebra},
  FJOURNAL = {Journal of Algebra},
    VOLUME = {382},
      YEAR = {2013},
     PAGES = {240--281},
      ISSN = {0021-8693},
   MRCLASS = {13F60 (05Exx)},
  MRNUMBER = {3034481},
MRREVIEWER = {Xueqing Chen},
       DOI = {10.1016/j.jalgebra.2013.02.018},
       URL = {https://doi.org/10.1016/j.jalgebra.2013.02.018},
}

@article {MS,
    AUTHOR = {Musiker, G. and Schiffler, R.},
     TITLE = {Cluster expansion formulas and perfect matchings},
   JOURNAL = {J. Algebraic Combin.},
  FJOURNAL = {Journal of Algebraic Combinatorics. An International Journal},
    VOLUME = {32},
      YEAR = {2010},
    NUMBER = {2},
     PAGES = {187--209},
      ISSN = {0925-9899},
   MRCLASS = {13F60 (05C25 05C70 16S99)},
  MRNUMBER = {2661414},
MRREVIEWER = {Michael S. Barot},
       DOI = {10.1007/s10801-009-0210-3},
       URL = {https://doi.org/10.1007/s10801-009-0210-3},
}

@article {FST,
    AUTHOR = {Fomin, S. and Shapiro, M. and Thurston, D.},
     TITLE = {Cluster algebras and triangulated surfaces. {I}. {C}luster
              complexes},
   JOURNAL = {Acta Math.},
  FJOURNAL = {Acta Mathematica},
    VOLUME = {201},
      YEAR = {2008},
    NUMBER = {1},
     PAGES = {83--146},
      ISSN = {0001-5962},
   MRCLASS = {57Q15 (13F60 32G15 52B70)},
  MRNUMBER = {2448067},
MRREVIEWER = {Christof Gei\ss },
       DOI = {10.1007/s11511-008-0030-7},
       URL = {https://doi.org/10.1007/s11511-008-0030-7},
}

@article {FT,
    AUTHOR = {Fomin, S. and Thurston, D.},
     TITLE = {Cluster algebras and triangulated surfaces {P}art {II}:
              {L}ambda lengths},
   JOURNAL = {Mem. Amer. Math. Soc.},
  FJOURNAL = {Memoirs of the American Mathematical Society},
    VOLUME = {255},
      YEAR = {2018},
    NUMBER = {1223},
     PAGES = {v+97},
      ISSN = {0065-9266},
   MRCLASS = {13F60 (30F60 57M50)},
  MRNUMBER = {3852257},
MRREVIEWER = {Christof Gei\ss },
       DOI = {10.1090/memo/1223},
       URL = {https://doi.org/10.1090/memo/1223},
}

@article {DW,
    AUTHOR = {Derksen, H. and Weyman, J.},
     TITLE = {Generalized quivers associated to reductive groups},
   JOURNAL = {Colloq. Math.},
  FJOURNAL = {Colloquium Mathematicum},
    VOLUME = {94},
      YEAR = {2002},
    NUMBER = {2},
     PAGES = {151--173},
      ISSN = {0010-1354},
   MRCLASS = {16G20 (14L35)},
  MRNUMBER = {1967372},
MRREVIEWER = {Iain G. Gordon},
       DOI = {10.4064/cm94-2-1},
       URL = {https://doi.org/10.4064/cm94-2-1},
}

@book {S_Quiver_rep,
    AUTHOR = {Schiffler, R.},
     TITLE = {Quiver representations},
    SERIES = {CMS Books in Mathematics/Ouvrages de Math\'{e}matiques de la SMC},
 PUBLISHER = {Springer, Cham},
      YEAR = {2014},
     PAGES = {xii+230},
   MRCLASS = {16G20 (16G70)},
  MRNUMBER = {3308668},
MRREVIEWER = {Alex Martsinkovsky},
       DOI = {10.1007/978-3-319-09204-1},
       URL = {https://doi.org/10.1007/978-3-319-09204-1},
}

@misc{boos2021degenerations,
      title={On degenerations and extensions of symplectic and orthogonal quiver representations}, 
      author={M. Boos and G. Cerulli Irelli},
      year={2021},
      eprint={2106.08666},
      archivePrefix={arXiv},
      primaryClass={math.RT}
}

@article {CC,
    AUTHOR = {Caldero, P. and Chapoton, F.},
     TITLE = {Cluster algebras as {H}all algebras of quiver representations},
   JOURNAL = {Comment. Math. Helv.},
  FJOURNAL = {Commentarii Mathematici Helvetici. A Journal of the Swiss
              Mathematical Society},
    VOLUME = {81},
      YEAR = {2006},
    NUMBER = {3},
     PAGES = {595--616},
      ISSN = {0010-2571},
   MRCLASS = {16G20 (18E30)},
  MRNUMBER = {2250855},
       DOI = {10.4171/CMH/65},
       URL = {https://doi.org/10.4171/CMH/65},
}

@article {CEHR,
    AUTHOR = {Cerulli Irelli, G. and Esposito, F. and Franzen,
              H. and Reineke, M.},
     TITLE = {Cell decompositions and algebraicity of cohomology for quiver
              {G}rassmannians},
   JOURNAL = {Adv. Math.},
  FJOURNAL = {Advances in Mathematics},
    VOLUME = {379},
      YEAR = {2021},
     PAGES = {Paper No. 107544, 47},
      ISSN = {0001-8708},
   MRCLASS = {14M15 (13F60 14C15 16G20)},
  MRNUMBER = {4198640},
MRREVIEWER = {Xiaobo Zhuang},
       DOI = {10.1016/j.aim.2020.107544},
       URL = {https://doi.org/10.1016/j.aim.2020.107544},
}

@article {DWZ,
    AUTHOR = {Derksen, H. and Weyman, J. and Zelevinsky, A.},
     TITLE = {Quivers with potentials and their representations {II}:
              applications to cluster algebras},
   JOURNAL = {J. Amer. Math. Soc.},
  FJOURNAL = {Journal of the American Mathematical Society},
    VOLUME = {23},
      YEAR = {2010},
    NUMBER = {3},
     PAGES = {749--790},
      ISSN = {0894-0347},
   MRCLASS = {16G20 (13F60)},
  MRNUMBER = {2629987},
MRREVIEWER = {M\'{a}ty\'{a}s Domokos},
       DOI = {10.1090/S0894-0347-10-00662-4},
       URL = {https://doi.org/10.1090/S0894-0347-10-00662-4},
}

@article {CAI,
    AUTHOR = {Fomin, S. and Zelevinsky, A.},
     TITLE = {Cluster algebras. {I}. {F}oundations},
   JOURNAL = {J. Amer. Math. Soc.},
  FJOURNAL = {Journal of the American Mathematical Society},
    VOLUME = {15},
      YEAR = {2002},
    NUMBER = {2},
     PAGES = {497--529},
      ISSN = {0894-0347},
   MRCLASS = {16S99 (14M99 17B99)},
  MRNUMBER = {1887642},
MRREVIEWER = {Eric N. Sommers},
       DOI = {10.1090/S0894-0347-01-00385-X},
       URL = {https://doi.org/10.1090/S0894-0347-01-00385-X},
}

@article {CAIV,
    AUTHOR = {Fomin, S. and Zelevinsky, A.},
     TITLE = {Cluster algebras. {IV}. {C}oefficients},
   JOURNAL = {Compos. Math.},
  FJOURNAL = {Compositio Mathematica},
    VOLUME = {143},
      YEAR = {2007},
    NUMBER = {1},
     PAGES = {112--164},
      ISSN = {0010-437X},
   MRCLASS = {16S99 (05E15 14M17 22E46)},
  MRNUMBER = {2295199},
MRREVIEWER = {Christof Gei\ss },
       DOI = {10.1112/S0010437X06002521},
       URL = {https://doi.org/10.1112/S0010437X06002521},
}

@article {S_CAUS,
    AUTHOR = {Schiffler, R.},
     TITLE = {On cluster algebras arising from unpunctured surfaces. {II}},
   JOURNAL = {Adv. Math.},
  FJOURNAL = {Advances in Mathematics},
    VOLUME = {223},
      YEAR = {2010},
    NUMBER = {6},
     PAGES = {1885--1923},
      ISSN = {0001-8708},
   MRCLASS = {13F60 (05E15 16G20)},
  MRNUMBER = {2601004},
MRREVIEWER = {Gregoire Dupont},
       DOI = {10.1016/j.aim.2009.10.015},
       URL = {https://doi.org/10.1016/j.aim.2009.10.015},
}

@article {LF,
    AUTHOR = {Labardini-Fragoso, D.},
     TITLE = {Quivers with potentials associated to triangulated surfaces},
   JOURNAL = {Proc. Lond. Math. Soc. (3)},
  FJOURNAL = {Proceedings of the London Mathematical Society. Third Series},
    VOLUME = {98},
      YEAR = {2009},
    NUMBER = {3},
     PAGES = {797--839},
      ISSN = {0024-6115},
   MRCLASS = {16G99 (16S99)},
  MRNUMBER = {2500873},
MRREVIEWER = {Gregoire Dupont},
       DOI = {10.1112/plms/pdn051},
       URL = {https://doi.org/10.1112/plms/pdn051},
}

@article {Dupont,
    AUTHOR = {Dupont, G.},
     TITLE = {An approach to non-simply laced cluster algebras},
   JOURNAL = {J. Algebra},
  FJOURNAL = {Journal of Algebra},
    VOLUME = {320},
      YEAR = {2008},
    NUMBER = {4},
     PAGES = {1626--1661},
      ISSN = {0021-8693},
   MRCLASS = {13A99 (16G20)},
  MRNUMBER = {2431998},
MRREVIEWER = {Laurent Demonet},
       DOI = {10.1016/j.jalgebra.2008.03.018},
       URL = {https://doi.org/10.1016/j.jalgebra.2008.03.018},
}

@article {Demonet,
    AUTHOR = {Demonet, L.},
     TITLE = {Categorification of skew-symmetrizable cluster algebras},
   JOURNAL = {Algebr. Represent. Theory},
  FJOURNAL = {Algebras and Representation Theory},
    VOLUME = {14},
      YEAR = {2011},
    NUMBER = {6},
     PAGES = {1087--1162},
      ISSN = {1386-923X},
   MRCLASS = {13F60 (16G10 16G20 17B35 17B67 18E10)},
  MRNUMBER = {2844757},
MRREVIEWER = {Gregoire Dupont},
       DOI = {10.1007/s10468-010-9228-4},
       URL = {https://doi.org/10.1007/s10468-010-9228-4},
}

@article {M,
    AUTHOR = {Musiker, G.},
     TITLE = {A graph theoretic expansion formula for cluster algebras of
              classical type},
   JOURNAL = {Ann. Comb.},
  FJOURNAL = {Annals of Combinatorics},
    VOLUME = {15},
      YEAR = {2011},
    NUMBER = {1},
     PAGES = {147--184},
      ISSN = {0218-0006},
   MRCLASS = {13F60 (05A19 05E15)},
  MRNUMBER = {2785761},
MRREVIEWER = {Gregoire Dupont},
       DOI = {10.1007/s00026-011-0088-3},
       URL = {https://doi.org/10.1007/s00026-011-0088-3},
}

@article {NS,
    AUTHOR = {Nakanishi, T. and Stella, S.},
     TITLE = {Diagrammatic description of {$c$}-vectors and {$d$}-vectors of
              cluster algebras of finite type},
   JOURNAL = {Electron. J. Combin.},
  FJOURNAL = {Electronic Journal of Combinatorics},
    VOLUME = {21},
      YEAR = {2014},
    NUMBER = {1},
     PAGES = {Paper 1.3, 107},
   MRCLASS = {13F60},
  MRNUMBER = {3177498},
MRREVIEWER = {Calin Chindris},
       DOI = {10.37236/3106},
       URL = {https://doi.org/10.37236/3106},
}

@article {R,
    AUTHOR = {Reading, N.},
     TITLE = {Dominance phenomena: mutation, scattering and cluster
              algebras},
   JOURNAL = {Trans. Amer. Math. Soc.},
  FJOURNAL = {Transactions of the American Mathematical Society},
    VOLUME = {376},
      YEAR = {2023},
    NUMBER = {2},
     PAGES = {773--835},
      ISSN = {0002-9947},
   MRCLASS = {13F60 (05E16 20F55 57Q15)},
  MRNUMBER = {4531662},
       DOI = {10.1090/tran/7888},
       URL = {https://doi.org/10.1090/tran/7888},
}

@article {GLS1,
    AUTHOR = {Geiss, C. and Leclerc, B. and Schr\"{o}er, J.},
     TITLE = {Quivers with relations for symmetrizable {C}artan matrices
              {I}: {F}oundations},
   JOURNAL = {Invent. Math.},
  FJOURNAL = {Inventiones Mathematicae},
    VOLUME = {209},
      YEAR = {2017},
    NUMBER = {1},
     PAGES = {61--158},
      ISSN = {0020-9910},
   MRCLASS = {16G10 (16G20 16G70)},
  MRNUMBER = {3660306},
MRREVIEWER = {Intan Muchtadi-Alamsyah},
       DOI = {10.1007/s00222-016-0705-1},
       URL = {https://doi.org/10.1007/s00222-016-0705-1},
}

@book {libro_blu,
    AUTHOR = {Assem, I. and Simson, D. and Skowro\'{n}ski, A.},
     TITLE = {Elements of the representation theory of associative algebras.
              {V}ol. 1},
    SERIES = {London Mathematical Society Student Texts},
    VOLUME = {65},
      NOTE = {Techniques of representation theory},
 PUBLISHER = {Cambridge University Press, Cambridge},
      YEAR = {2006},
     PAGES = {x+458},
   MRCLASS = {16G10 (16-02)},
  MRNUMBER = {2197389},
MRREVIEWER = {Peter W. Donovan},
       DOI = {10.1017/CBO9780511614309},
       URL = {https://doi.org/10.1017/CBO9780511614309},
}

@book {AR_book,
    AUTHOR = {Auslander, M. and Reiten, I. and Smal\o , S. O.},
     TITLE = {Representation theory of {A}rtin algebras},
    SERIES = {Cambridge Studies in Advanced Mathematics},
    VOLUME = {36},
      NOTE = {Corrected reprint of the 1995 original},
 PUBLISHER = {Cambridge University Press, Cambridge},
      YEAR = {1997},
     PAGES = {xiv+425},
   MRCLASS = {16Gxx (16-02)},
  MRNUMBER = {1476671},
}

@article {FST_finite_type,
    AUTHOR = {Felikson, A. and Shapiro, M. and Tumarkin, P.},
     TITLE = {Cluster algebras of finite mutation type via unfoldings},
   JOURNAL = {Int. Math. Res. Not. IMRN},
  FJOURNAL = {International Mathematics Research Notices. IMRN},
      YEAR = {2012},
    NUMBER = {8},
     PAGES = {1768--1804},
      ISSN = {1073-7928},
   MRCLASS = {13F60},
  MRNUMBER = {2920830},
MRREVIEWER = {Yu Zhou},
       DOI = {10.1093/imrn/rnr072},
       URL = {https://doi.org/10.1093/imrn/rnr072},
}

@article {GeuLF2,
    AUTHOR = {Geuenich, J. and Labardini-Fragoso, D.},
     TITLE = {Species with potential arising from surfaces with orbifold
              points of order 2, {P}art {II}: {A}rbitrary weights},
   JOURNAL = {Int. Math. Res. Not. IMRN},
  FJOURNAL = {International Mathematics Research Notices. IMRN},
      YEAR = {2020},
    NUMBER = {12},
     PAGES = {3649--3752},
      ISSN = {1073-7928},
   MRCLASS = {13F60 (05E14 16G20 18G80)},
  MRNUMBER = {4120307},
MRREVIEWER = {Daniel V. Mathews},
       DOI = {10.1093/imrn/rny090},
       URL = {https://doi.org/10.1093/imrn/rny090},
}

@article {GeuLF1,
    AUTHOR = {Geuenich, J. and Labardini-Fragoso, D.},
     TITLE = {Species with potential arising from surfaces with orbifold
              points of order 2, part {I}: one choice of weights},
   JOURNAL = {Math. Z.},
  FJOURNAL = {Mathematische Zeitschrift},
    VOLUME = {286},
      YEAR = {2017},
    NUMBER = {3-4},
     PAGES = {1065--1143},
      ISSN = {0025-5874},
   MRCLASS = {05E99 (13F60 16G20)},
  MRNUMBER = {3671571},
       DOI = {10.1007/s00209-016-1795-6},
       URL = {https://doi.org/10.1007/s00209-016-1795-6},
}

@article {FeST,
    AUTHOR = {Felikson, A. and Shapiro, M. and Tumarkin, P.},
     TITLE = {Cluster algebras and triangulated orbifolds},
   JOURNAL = {Adv. Math.},
  FJOURNAL = {Advances in Mathematics},
    VOLUME = {231},
      YEAR = {2012},
    NUMBER = {5},
     PAGES = {2953--3002},
      ISSN = {0001-8708},
   MRCLASS = {13F60},
  MRNUMBER = {2970470},
MRREVIEWER = {Kyungyong Lee},
       DOI = {10.1016/j.aim.2012.07.032},
       URL = {https://doi.org/10.1016/j.aim.2012.07.032},
}

@misc{5,
      title={Cluster expansion formulas and perfect matchings for type B and C}, 
      author={A. Ciliberti},
      year={2024},
      eprint={2405.14915},
      archivePrefix={arXiv},
      primaryClass={math.RT}
}

@article {MSW11,
    AUTHOR = {Musiker, G. and Schiffler, R. and Williams, L.},
     TITLE = {Positivity for cluster algebras from surfaces},
   JOURNAL = {Adv. Math.},
  FJOURNAL = {Advances in Mathematics},
    VOLUME = {227},
      YEAR = {2011},
    NUMBER = {6},
     PAGES = {2241--2308},
      ISSN = {0001-8708},
   MRCLASS = {13F60 (05C70 05E15)},
  MRNUMBER = {2807089},
MRREVIEWER = {Gregoire Dupont},
       DOI = {10.1016/j.aim.2011.04.018},
       URL = {https://doi.org/10.1016/j.aim.2011.04.018},
}

@article {FeliksonTumarkin,
    AUTHOR = {Felikson, A. and Tumarkin, P.},
     TITLE = {Bases for cluster algebras from orbifolds},
   JOURNAL = {Adv. Math.},
  FJOURNAL = {Advances in Mathematics},
    VOLUME = {318},
      YEAR = {2017},
     PAGES = {191--232},
      ISSN = {0001-8708},
   MRCLASS = {13F60 (16S99)},
  MRNUMBER = {3689740},
MRREVIEWER = {Linhui Shen},
       DOI = {10.1016/j.aim.2017.07.025},
       URL = {https://doi.org/10.1016/j.aim.2017.07.025},
}

\end{document}